\documentclass{article}
\usepackage{graphicx} 
\usepackage{amsmath,amsthm,amsfonts,amssymb}
\usepackage{xcolor}
\usepackage{fullpage}
\usepackage{mathtools}
\usepackage{authblk}

\usepackage{latexsym,color,amsmath,amsthm,amssymb,amscd,amsfonts}
\usepackage{tikz}
\usetikzlibrary{arrows}
\usepackage{scalefnt}
\usepackage[font=small,labelfont=bf]{caption}
\usepackage{float}
\usepackage{graphicx}
\usepackage{relsize}
\usepackage{amsopn}
\usepackage{mathrsfs}  
\usepackage[hidelinks]{hyperref}
\usepackage{bm}
\usepackage{bbm}
\usepackage[shortlabels]{enumitem}
\usepackage{multicol}
\usepackage{float} 
\usepackage{todonotes}
\usepackage{cite}
\usepackage{soul}
\usepackage{cancel}

\newtheorem{theorem}{Theorem}[section]
\newtheorem{question}[theorem]{Question}

\newtheorem{lemma}[theorem]{Lemma}
\newtheorem{corollary}[theorem]{Corollary}
\newtheorem{proposition}[theorem]{Proposition}
\newtheorem{example}[theorem]{Example}

\newtheorem{remark}[theorem]{Remark}
\newtheorem{definition}[theorem]{Definition}

\newtheorem{mainthm}{Theorem}

\newcommand{\R}{{\mathbb R}}
\newcommand{\sfe}{\mathbb{S}^{n-1}}
\newcommand{\RR}{{\mathbb R}}
\newcommand{\N}{{\mathbb N}}
\newcommand{\sn}{\mathbb{S}^{n-1}}
\newcommand{\fconvx}{{\mbox{\rm Conv}(\R^n)}} 
\newcommand{\fconvf}{{\mbox{\rm Conv}(\R^n; \R)}}
\newcommand{\class}{{\mathcal C}}
\newcommand{\bs}{\mathcal{BS}}
\newcommand{\dom}{{\mbox{\rm dom}}}
\newcommand{\inte}{{\mbox{\rm int}}}

\newcommand{\epi}{{\mbox{\rm epi}}}
\newcommand{\dd}{\mathrm{d}}

    \def\XXint#1#2#3{{\setbox0=\hbox{$#1{#2#3}{\int}$}
         \vcenter{\hbox{$#2#3$}}\kern-.5\wd0}}

    \makeatother

\DeclarePairedDelimiterX{\ip}[2]{\langle}{\rangle}{#1, #2}

\title{On weighted Blaschke--Santal\'o\\ and strong Brascamp--Lieb inequalities}

\author[1]{Andrea Colesanti}
\author[2]{Alexander Kolesnikov}
\author[3]{Galyna Livshyts}
\author[4]{Liran Rotem}
\affil[1]{University of Florence, Italy}
\affil[2]{HSE University, Russian Federation}
\affil[3]{Georgia Institute of Technology, USA}
\affil[4]{Technion - Israel Institute of Technology, Israel}
\date{}                     
\begin{document}

\maketitle

\begin{abstract}
In this paper, we study new extensions of the functional Blaschke-Santal\'o inequalities, and explore applications of such new inequalities beyond the classical setting of the standard Gaussian measure. 
In particular, we study functionals of the type
\begin{equation}
\label{BSfunctional}
\Bigl( \int_{\mathbb{R}^n} e^{-\Phi} dx\Bigr)^{\frac{1}{p}}
\Bigl( \int_{\mathbb{R}^n} e^{-\frac{1}{p-1}\Phi^*(\nabla V)} dx\Bigr)^{\frac{p-1}{p}},
\end{equation}
where $p >1$ and $V$ is a $p$-homogeneous convex even function. The function $\Phi$ is assumed to be convex and even.

In particular, we prove that maximizers of (\ref{BSfunctional}) are $p$-homogeneous and solve an equation of the Monge--Amp\`ere type appearing in the $L^{q}$-Minkowski problem. 
 This gives  a novel mass-transport approach to the Blaschke-Santal\'o inequality, which is of independent interest even in the classical setting. We find sufficient conditions for $\Phi=V$ to be the maximizer of (\ref{BSfunctional}). 
 In particular, these conditions are satisfied if $V = \|x\|_q^{p}$, where $p\ge q \ge 2$ and $\|\cdot\|_{q}$ is the $l^q$-norm.
 We prove that in general   $V$ fails to be the maximizer.

In addition, we prove that any maximizer of (\ref{BSfunctional}) satisfies a strong version of the Brascamb--Lieb inequality on the set of even functions. In particular, if $V$ is the maximizer of (\ref{BSfunctional}), then  $\mu  = \frac{e^{-V}dx}{\int e^{-V}dx}$ satisfies the following  strengthening  of the Brascamp--Lieb inequality 
on the set of even functions:
\begin{equation}
\label{lbl}
{\rm Var}_{\mu} f \le \lambda\int_{\mathbb{R}^n} \langle (D^2 V)^{-1} \nabla f, \nabla f \rangle d\mu
\end{equation} with the sharp value 
$
\lambda = 1 - \frac{1}{p}.
$
We also estimate the best constant $\lambda <1 $ in (\ref{lbl})
for probability measures of the form 
$\mu = C e^{-c\|x\|^{p}_q} dx$ for various values of $p$ and $q$.


\end{abstract}

\section{Introduction}

The Blaschke--Santal\'o inequality 
\begin{equation}\label{BS-bodies}
|K| |K^{o}| \le |B^n_2|^2
\end{equation}
discovered in the beginning of XXth century by W.~Blaschke in dimensions $n=2,3$ and proved later by L.~Santal\'o \cite{santalo} for $n>3$,  is one of the fundamental and most celebrated results in convex geometry. Here $K \subset \mathbb{R}^n$ is a symmetric convex body, $|\cdot|$ is the $n$-dimensional volume, and
$$K^o = \{y: \langle x, y \rangle \le 1 \ \forall\ x \in K \}$$
is the corresponding polar body. The notation $B^n_p$ stands for the unit $l^p$-ball in $\R^n$:
$$
B^n_p = \Bigl\{x: \|x\|_p = \Bigl(\sum_{i=1}^n |x_i|^p \Bigr)^{\frac{1}{p}} \le 1\Bigr\}.
$$
The standard $l^p$-norm of vector $x$ is denoted by $\|x\|_p$. In the case $p=2$ we  omit subscript $2$ and  write $|x|$.

The classical proof of (\ref{BS-bodies})  goes via the Steiner symmetrization and the Brunn--Minkowski inequality, see e.g. Artstein-Avidan, Giannopolous, Milman \cite{AGM}, Campi, Gronchi \cite{campigronchi}, Meyer, Pajor \cite{MePa}. The first symmetrization proof was obtained by Saint Raymond \cite{SaintR}. An important feature of the \emph{volume product} functional $|K| |K^{o}|$, which appears on the left hand side of this inequality, is its linear invariance: for any linear transformation $T$ on $\R^n,$ we have $$|K| |K^{o}|=|T(K)| |(T(K))^{o}|.$$

A remarkable functional form of the Blaschke--Santal\'o inequality was discovered by K.~Ball in 
\cite{Ball1986}. 
Let $\Phi$ be an arbitrary proper even function on $\mathbb{R}^n$ with values in $(-\infty, +\infty]$ and 
$$\Phi^*(y) = \sup_{x \in \mathbb{R}^n} \bigl(\langle x, y \rangle - \Phi(x)\bigr)$$ 
be its Legendre transform.
Then 
\begin{equation}
\label{BSBall-G}
\int e^{-\Phi(x)} dx \int e^{-\Phi^*(y)} dy
\le (2 \pi)^n.
\end{equation}
Here and throughout $\int$ stands for $\int_{\R^n}$. The equality holds if and only if $\Phi = a + \langle A x, x \rangle$ for some symmetric non-degenerate matrix $A$.
This result was later generalized by  Artstein-Avidan, Klartag, Milman \cite{AKM},  Fradelizi, Meyer \cite{FradeliziMeyer2007}, \cite{FradeliziMeyer2008}, \cite{FradeliziMeyer2008(2)}, Lehec \cite{Lehec}, Lin, Leng \cite{LinLeng}, Kolesnikov, Werner \cite{KolesnikovWerner}, Gozlan, Fradelizi, Sadovsky, Zugmeyer \cite{GFSZ}. Stability of the inequality has been studied by  B\"or\"oczky \cite{Bor}. Some $L^p$-extensions of the inequality have been obtained by Lutwak and Zhang \cite{LZ} and Berndtsson, Mastrantonis, Rubinstein\cite{BMR}.  A Fourier analytic proof is presented by Bianchi, Kelly \cite{BK}.
 A novel analytic approach  (semigroups and a monotonicity property) to (\ref{BSBall-G})
is presented in the recent works of Nakamura and Tsuji \cite{NakTs} and Cordero-Erausquin, Fradelizi and Langharst \cite{CEFL}. See also another very recent result of Nakamura and Tsuji \cite{NakTs2} with a partial solution of the  conjecture from
\cite{KolesnikovWerner}. We remark that both (\ref{BS-bodies}) and (\ref{BSBall-G}) have a conjectured reverse form, referred to as Mahler conjecture, which shall not be discussed in detail in this work. However, an intense research activity has been carried out in relation to this conjecture; we refer the interested reader to Fradelizi, Meyer, Zvavitch \cite{vol-prod-survey} and the references therein.

Recall that the standard Gaussian measure $\gamma$ on $\R^n$ is the measure with the density $(2\pi)^{-\frac{n}{2}}e^{-\frac{|x|^2}{2}}$. The inequality (\ref{BSBall-G}) states that the functional 
$$
\int e^{-\Phi(x)} dx \int e^{-\Phi^*(y)} dy
$$ is maximized when $\Phi$ is the potential of a Gaussian measure (in particular, the standard one). Hence the inequality (\ref{BSBall-G}) seems like a purely Gaussian phenomenon. However, Fradelizi and Meyer \cite{FradeliziMeyer2007} extended the formulation of (\ref{BSBall-G}) in a way so that other rotation-invariant measures appear as maximizers of similar inequalities. 

In this paper, we study a different possibility to extend the Blaschke-Santal\'o inequality, with the goal of potentially getting a richer class of measures as maximizers of this type of functionals. Namely, we pose the following questions.

\begin{question}\label{Question-main}

What are the  maximizers for the \emph{generalized weighted
Blaschke--Santal\'o functional}
$$
{\mathcal BS}_{\alpha, \beta, \rho_1,\rho_2}(\Phi) = \Bigl( \int e^{-\alpha \Phi(x)} \rho_1(x) dx\Bigr)^{\frac{1}{\alpha}} \Bigl( \int e^{-\beta \Phi^*(y)} \rho_2(y) dy \Bigr)^{\frac{1}{\beta}}, \ \alpha, \beta \in (0,+\infty),
$$
with symmetric weights $\rho_1, \rho_2$ on the set of even functions? Of interest are existence, uniqueness and characterizations of the maximizers of ${\mathcal BS}_{\alpha,\beta,\rho_1,\rho_2}$, and especially situations  when ${\mathcal BS}_{\alpha,\beta,\rho_1,\rho_2}$ admits a closed-form solution.
\end{question}

\medskip

In the partial case when $\alpha=\beta$ and $\rho_1=\rho_2$ is an even log-concave weight, this question was answered by Klartag \cite{Klartag-marginals}, and the maximizer is $\Phi(x)=\frac{|x|^2}{2}.$ In general, the setting of Question \ref{Question-main} is very general even for the problem of existence of maximizers. In the largest part of this work we deal mainly with a particular case of the generalized Blaschke--Santal\'o functional.
Let $p > 1$ and let $V$ be an even 
strictly convex $p$-homogeneous $C^2$ function on $\R^n$.
We consider the functional
$$
\bs_{p,V}(\Phi) = \int e^{-\Phi(x)} dx \left(\int e^{-\frac{1}{p-1}\Phi^*(\nabla V(x))} dx \right)^{p-1}.
$$
Note that $\bs_{p,V}$
is a particular case of ${\mathcal BS}_{\alpha, \beta, \rho_1,\rho_2}$, because after the 
change of variables  $y = \nabla V(x)$ it can  be rewritten  as
$$
\bs_{p,V}(\Phi) = 
\int e^{-\Phi(x)} dx \left(\int e^{-\frac{1}{p-1}\Phi^*(y)} \det(D^2 V^*(y))dy\right)^{p-1}.$$
A remarkable property of $\bs_{p,V}$ is  the following homogeneity invariance:
$$
\bs_{p,V} (\Phi(tx))
= \bs_{p,V} (\Phi(x))
$$
for all $\Phi$ and $t >0$.
Another important property of $\bs_{p,V}$
is that the function $V$ (more generally $V(tx)$) is a natural candidate to be a maximizer. Indeed, this hint comes from the observation that $V$ satisfies the corresponding Euler--Lagrange equation for the minimization problem defining $\bs_{p,V}$ (see Proposition \ref{Euler-Lagrange}). Thus the following question appears naturally.

\begin{question}\label{Question-main-V}Let $p > 1$ and let $V$ be an even strictly convex $p$-homogeneous function on $\R^n$. Under which conditions do we have
\begin{equation}
\label{main-question-intro}
\int e^{-\Phi(x)} dx \left(\int e^{-\frac{1}{p-1}\Phi^*(\nabla V(x))} dx\right)^{p-1}
\le \left(\int e^{-V(x)}dx\right)^p,
\end{equation}
for an arbitrary proper convex even function $\Phi$ on $\mathbb{R}^n$ with values in $(-\infty, +\infty]$?
\end{question}
\medskip

As we shall explain in Section \ref{set-version} (following the ideas of Artstein-Avidan, Klartag, Milman \cite{AKM}), inequality (\ref{main-question-intro}) of Question \ref{Question-main-V} is equivalent to certain Blaschke--Santal\'o type inequality for sets.

\begin{proposition}
\label{set-func}
    Let $p > 1$ and let $V$ be an even strictly convex $p$-homogeneous $C^2$ function on $\R^n$. 
    Inequality (\ref{main-question-intro}) holds for an arbitrary convex proper function $\Phi$ if and only if inequality 
    \begin{equation}
    \label{KKo}
|K|\cdot|\nabla V^* (K^o)|^{p-1}\leq \left|\left\{V\leq \frac{1}{p}\right\}\right|^p
\end{equation}
holds for an arbitrary compact convex body $K$. 

If inequality (\ref{KKo}) holds, then equality is attained when $K$ is a sublevel set of $V$: $K = \{ V \le \alpha \}$. 
\end{proposition}

\begin{remark}
    Some alternative forms of the Blaschke--Santal\'o inequality generalizing  the classical inequality for sets appeared in the study of the generalized Minkowski problem. See \cite{CW}, \cite{CCL}, \cite{chen2018} and comments in Section 7.
\end{remark}

Unfortunately, one can not hope for the affirmative answer to Question \ref{Question-main-V} in too great a generality. To see this, in  view of Proposition \ref{set-func}, consider $V(x)=\|x\|^p_p$ with $p\in [1,2]$. In this case, (\ref{KKo}) states that the maximizer of the functional 
$$|K|\left(\int_{K^o} \prod_{i=1}^n |x_i|^{\frac{2-p}{p-1}} dx\right)^{p-1}$$
among all symmetric convex sets is $K=\frac{1}{p}B^n_p.$ In the limiting case $p\rightarrow 1,$ this would be equivalent to the inequality
$$|K|\sup_{x\in K^o}\prod_{i=1}^n |x_i|\leq |B^n_1|\sup_{x\in B^n_{\infty}}\prod_{i=1}^n |x_i|=\frac{2^n}{n!}.$$
However, this is false: if we let $v=e_1+...+e_n$, and $K^o_R=[-Rv,Rv]+B^n_2$ (which, for large $R$, is a long needle-like body pointing in the direction of $v$), then 
$$
\lim_{R\rightarrow\infty} \Bigl( |K_R|\sup_{x\in K_R^o}\prod_{i=1}^n |x_i|\Bigr)  = \infty.
$$


In a similar manner, one may see that the lack of linear invariance of the functional
$\bs_{p,V}$
makes (\ref{main-question-intro}) hopeless in many situations. 
See Section \ref{nonexistence} for more details and other counterexamples. Nevertheless, in the following Theorem we were able to obtain some sufficient conditions for (\ref{main-question-intro}).

\begin{mainthm}\label{main-thm}
Let $p > 1$ and let $V$ be an even strictly convex $p$-homogeneous $C^2$ function on $\R^n$. Assume that $V$ is an unconditional function, and that the function
$$
x=(x_1,\dots,x_n)\mapsto V\bigl(x^{\frac{1}{p}}_1,...,x^{\frac{1}{p}}_n\bigr)
$$ 
is concave in 
$$
\R^n_+=\{(x_1,\dots,x_n)\colon x_i\ge 0, \ \ 1 \le i \le n\}.
$$ 
Then inequality (\ref{main-question-intro}) holds for every unconditional convex $\Phi$.

Assume, in addition, that  for every coordinate hyperplane $H$, with unit normal $e$, and for every $x'\in H$, the function $\varphi\colon[0,+\infty)\to\R$ defined by
$$
\varphi(t)=\det D^2 V^*(x'+te)
$$
is decreasing.
 Then inequality 
 (\ref{main-question-intro})
holds for every even convex $\Phi$.
\end{mainthm}

\begin{corollary}\label{mainex} Let $V = c \|x\|_q^{p}$, $c\ge0$. Then inequality (\ref{main-question-intro}) holds in the following cases:
    \begin{enumerate}
        \item For $p \ge q >1$ and unconditional $\Phi$
        \item For $p \ge q \ge 2$ and even $\Phi$.
    \end{enumerate}
\end{corollary}
    
The proof of the ``unconditional" part of Theorem \ref{main-thm} is based on the application of the Pr\'ekopa--Leindler inequality in the unconditional case via a change of variables on $\mathbb{R}_+^n$  (see e.g. Fradelizi, Meyer \cite{FradeliziMeyer2007} for another application of this idea). The general case follows by a symmetrization argument: we show that Steiner symmetrization increases the value of the functional and reduce the problem to the unconditional case.

In Section \ref{section-existence-results} we prove that under various assumptions of additional symmetries (for $V$ and for $\Phi$) the left hand side of our inequality admits a non-degenerate maximizer. See, for instance, Theorem \ref{homog-rad symm} for the existence in the case of rotation-invariant weights while $\Phi$ is even, and Theorem \ref{existence homogeneous case0} for the existence when both the weights and $\Phi$ are assumed to be $1$-symmetric. In particular, in Section \ref{section-symmetrization} we use symmetrization techniques to show existence and various properties of such maximizers.

From our perspective, one of the most exciting aspects of this work is the novel mass transport approach to proving the classical  Blaschke-Santal\'o inequality. We note that the transportation approach to the converse form of the Blaschke--Santal{\'o} inequality was recently developed by Fradelizi, Gozlan, Zugmeyer in \cite{FGZ}. The starting point for our mass transport approach is an observation, going back to \cite{KolesnikovWerner}, that any maximizer $\Phi$ for the Blaschke--Santal\'o functional
is a solution for a nonlinear PDE of the Monge--Amp\`ere type.

\begin{theorem}[\cite{KolesnikovWerner}] Let $\alpha, \beta > 0$, and let $\rho_1, \rho_2$ be positive even functions. Assume that $\Phi$ is a maximum point of the functional
$\mathcal{BS}_{\alpha, \beta, \rho_1, \rho_2}$. Then $\nabla \Phi$ is the optimal transportation pushing forward $\mu$ onto $\nu$, where
$$
d\mu = \frac{ e^{-\alpha \Phi} \rho_1 dx }{\int e^{-\alpha \Phi} \rho_1 dx }, \  d\nu = \frac{ e^{-\beta \Phi^*} \rho_2 dy }{\int e^{-\beta \Phi^*} \rho_2 dy }.
$$
\end{theorem}

This result implies, in particular, that $\Phi$ solves the following nonlinear PDE of the Monge--Amp\`ere type.
\begin{equation}
    \label{MA12-0-intro}
\frac{e^{-\alpha\Phi}}{\int e^{-\alpha\Phi} \rho_1 dx} \rho_1= 
\frac{e^{-\beta\Phi^*(\nabla \Phi)}}{\int e^{-\beta\Phi^*} \rho_2 dy}
\rho_2(\nabla \Phi) \det D^2 \Phi.
\end{equation}

Of course, there is no hope to find a closed-form solution for this equation for general $\rho_1, \rho_1$. However, in some particular cases there exists a natural candidate, as is demonstrated by our results. For the case $\alpha=\beta=1, \rho_1=\rho_2=1$ this equation was already studied in \cite{KolesnikovWerner}. It was proved there that under additional regularity assumptions all solutions to equation (\ref{MA12-0-intro}) are positive quadratic forms: $\Phi(x) = \langle A x, x \rangle$. This would not prove the classical Blaschke--Santal\'o inequality, unless we know {\em a priori} that any solution to 
\begin{equation}
    \label{MA12-00-intro}
\frac{e^{-\alpha\Phi}}{\int e^{-\alpha\Phi} dx} = 
\frac{e^{-\beta\Phi^*(\nabla \Phi)}}{\int e^{-\beta\Phi^*} dy} \det D^2 \Phi
\end{equation}
is sufficiently regular; in this paper we prove that any solution to (\ref{MA12-00-intro}) (understood in the mass transportation sense) is indeed regular, which gives, in particular, a transportation proof of the functional Blaschke--Santal\'o inequality (\ref{BSBall-G})
(see Remark \ref{BSMA}).

When trying to extend these arguments to other cases, we face the difficulty that  Monge--Amp\`ere
equation (\ref{MA12-0-intro}) may have many solutions. Indeed, let us consider the functional $\bs_{p,V}$.
We prove the following result.

\begin{mainthm}\label{MA-homogen-solutions} Let  $V$ be convex and $p$-homogeneous. Then any maximum point $\Phi$ of the functional $\bs_{p,V}$ satisfying $\Phi(0)=0$ is $p$-homogeneous.
\end{mainthm}

Using homogeneity one can reduce equation
for the maximizer of $\bs_{1,V}$ to a Monge--Amp\`ere equation on the unit sphere. More precisely, we get that it is equivalent to the so-called  $L^{q}$-Minkowski problem for some corresponding $q(p,n)$.  The latter is the following non-linear elliptic problem on the sphere: given a measure $\mu = f dx$ on $\mathbb{S}^{n-1}$ solve equation of the type
    \begin{equation}
    \label{bmin}
h^{1-q} \det (h_{ij} + h \delta_{ij}) = f.
    \end{equation}
Uniqueness of solution to (\ref{bmin}) would provide a natural way of establishing affirmative answer to Question \ref{Question-main-V}. Unfortunately, it is known that in general equation (\ref{bmin}) has many  (even infinitely many) solutions for those values of $q$ which are of interest for us. For instance, for $p=2$ we get $q=-n$ and this is the so-called centro-affine Minkowski problem.
We refer to \cite{CW}, \cite{HLW}, \cite{JLW}, \cite{QRLi}  for examples of non-uniqueness.
It is known that uniqueness in the Minkowski problem is closely related to the conjectured $L^p$-Brunn--Minkowski inequality. See, in particular,  seminal paper \cite{BLYZ} about log-Brunn--Minkowski inequality.
Some uniqueness results  via $L^p$-Brunn--Minkowski inequality
can be found in  \cite{KM-LpBMproblem}, \cite{MI1}, \cite{MI2}. 
 More information the reader can find in the recent book  \cite{BFR}. See Subsection \ref{lqbm} for  details.

\medskip

Let us now discuss some possible motivations for studying the aforementioned questions. The inequality (\ref{BSBall-G}) is related to many other important and interesting results, such as the Reverse Log-Sobolev inequality of Artstein-Avidan, Klartag, Sch\"utt, Werner \cite{AKSW} (see also Caglar, Fradelizi, Gozlan, Lehec, Sch\"utt, Werner \cite{CFGLSW}), and the generalized Talagrand's Transport-Entropy inequality due to Fathi \cite{fathi}. Therefore one would hope that our results pave the way to discovering this type of phenomena also beyond the Gaussian setting. 

In order to explain a particularly important motivation, let us recall the connection of the Blaschke-Santal\'o inequality to the sharp symmetric Gaussian Poincar\'e inequality. Inequality (\ref{BSBall-G}) implies that for any $t>0$ and any \emph{even} $f \in C^2(\R^n)$, the function
$$F(t)=\int e^{-(\frac{x^2}{2}+tf(x))} dx \int e^{-(\frac{y^2}{2}+tf(y))^*(y)} dy$$
is maximized at $t=0$. 

One may thus check that $F'(0)=0$ and deduce that $F''(0)\leq 0$. Computing $F''(0)$ we obtain  the ``symmetric Gaussian Poincare inequality'': for any even locally-Lipschitz function $f:\R^n\rightarrow\R,$ one has
\begin{equation}\label{half-Poincare}
\int f^2 d\gamma-\left(\int f d\gamma\right)^2\leq \frac{1}{2}\int |\nabla f|^2 d\gamma.
\end{equation}
(see more details e.g. at \cite{Galyna-notes}). This fact can be also proven using Hermite polynomial decomposition, and in many other ways -- see {\em e.g.} Cordero-Erausquin, Fradelizi, Maurey \cite{CEFM}. 

Using an analogous variational argument, one can get
\begin{proposition}\label{fromBStoBL} Let $\Phi$ be the maximum point of $ \mathcal{BS}_{\alpha, \beta, \rho_1, \rho_2}$. Then $\mu = \frac{ e^{-\alpha \Phi} \rho_1 dx }{\int e^{-\alpha \Phi} \rho_1 dx }$ satisfies
$$
{\rm Var}_{\mu} f \le \frac{1}{\alpha + \beta} \int \langle ( D^2 \Phi)^{-1} \nabla f, \nabla f \rangle d\mu.
$$ 
\end{proposition}

The result of Proposition \ref{fromBStoBL} (for $\rho_1=1$) is an improvement of the Brascamp--Lieb inequality for maximizers of the Blaschke--Santal\'o functional on the set of even functions. Recall that a log-concave probability measure 
$$
d\mu = \frac{e^{-V}dx}{\int e^{-V}dx}
$$ 
satisfies the famous Brascamp--Lieb inequality (see \cite{BrLi})
$$
{\rm Var}_{\mu} f \le \int \langle (D^2 V)^{-1} \nabla f, \nabla f \rangle d\mu
$$
for any smooth (in general, not symmetric) $f$, under the assumption that $V\in C^2(\R^n)$ and $D^2V(x)$ is positive definite for every $x$. The last inequality can be viewed as the infinitesimal version of the Pr\'ekopa--Leindler inequality (see \cite{BL1-BrLib}). The Brascamp--Lieb inequality is sharp with equality case $f = V_{x_i}$.

Can the Brascamp--Lieb inequality be improved on the set of even functions in the spirit of inequality (\ref{half-Poincare})? One of the motivation for this question comes from the attempts to solve the so-called $B$-conjecture (see \cite{CEFM}). In particular, it was shown in \cite{CorRot} that the following inequality is equivalent to the $B$-conjecture:
$$
\frac{1}{\mu(K)} \int_K \langle \nabla V,x\rangle^2 d\mu-\Bigl( \frac{1}{\mu(K)} \int_K \langle \nabla V,x\rangle d\mu\Bigr)^2\leq \frac{1}{\mu(K)} \int_K (\langle \nabla V,x\rangle +\langle \nabla^2 V x,x\rangle) d\mu,
$$
where $K$ is an arbitrary even symmetric convex set. The latter inequality can be viewed as a particular case of a strengthening  of the Brascamp--Lieb inequality for measure $\frac{1}{\mu(K)} 1_{K} \cdot \mu$ for  specific function $f(x) = \langle \nabla V(x), x \rangle$. Here $1_K$ is the indicator function of $K$.  

A reasonable guess of what can be the strengthening of the Brascamp--Lieb inequality is the following: let $V$ be a $p$-homogeneous even convex function and $\mu=\frac{e^{-V} dx}{\int e^{-V}dx}$. Is it true that for every even function $f$ the following holds?
\begin{equation}
\label{p-BL}
{\rm Var}_{\mu} f \le \Bigl( 1 - \frac{1}{p}\Bigr)  \int \langle (D^2 V)^{-1} \nabla f, \nabla f \rangle d\mu.
\end{equation}
Note that the conjectured inequality turns out to be an equality for $f(x) = \langle \nabla V(x), x \rangle = p V(x)$.

As a consequence of Proposition \ref{fromBStoBL} together with Theorem \ref{main-thm}, we derive the following strengthening of the Brascamp-Lieb inequality:

\begin{corollary}\label{thm-BL}
Let $p > 1$ and let $V$ be an even strictly convex $p$-homogeneous $C^2$ function on $\R^n$. Assume that $V$ is an unconditional function, and that the function
$$
x=(x_1,\dots,x_n)\mapsto V\bigl(x^{\frac{1}{p}}_1,...,x^{\frac{1}{p}}_n\bigr)
$$ 
is concave in $\R^n_+$. Then inequality (\ref{p-BL}) holds for every unconditional $f$.

Assume, in addition, that for every coordinate hyperplane $H$, with unit normal $e$ and for every $x'\in H$, the function $\varphi\colon[0,+\infty)\to\R$ defined by
$$
\varphi(t)=\det D^2 V^*(x'+te)
$$
is decreasing. Then inequality (\ref{p-BL}) holds for every even $f$.
\end{corollary}

\medskip

Similarly to generalized Blaschke--Santal\'o inequality, inequality (\ref{p-BL}) fails to hold for an arbitrary convex even $p$-homogeneous $V$. In Section \ref{section-BL} we study a strengthening of the Brascamp--Lieb inequality for a specific family of measures. 

\begin{mainthm}\label{strongBL} Let $p,q >1$. Consider the probability measure $\mu$ given by
$$
d\mu = C {e^{-\frac{1}{p} \|x\|_q^{p}} dx}= C {e^{-\frac{1}{p} \bigl( \sum_{i=1}^n |x_i|^q \bigr)^{\frac{p}{q}}}dx},
$$
where $C > 0$ is a normalizing constant. Assume that
$$
\lambda \ge \max \Bigl( 1 -\frac{1}{p}, 1 - \frac{1}{q}, \frac{1}{2(1 + \frac{q-2}{n})} \Bigr).
$$
Then, the inequality 
\begin{equation}\label{BL1008}
{\rm Var}_{\mu} f \le \lambda  \int \langle (D^2 V)^{-1} \nabla f, \nabla f \rangle d\mu,
\end{equation}
holds on the set of even functions. This inequality is sharp and holds with $\lambda = 1 - \frac{1}{p}$ in the following cases:
\begin{enumerate}
    \item 
    $$
 q \ge 2, \ \ p \ge q,
    $$
    \item 
    $$q \le 2, \ \ p \ge \frac{2(n+q-2)}{n+ 2(q-2)} = 2 - \frac{2(q-2)}{n+2(q-2)}.$$
\end{enumerate}
\end{mainthm}

Note that statement (1) of Theorem \ref{strongBL} follows from Corollary \ref{mainex}
and Proposition \ref{fromBStoBL}. Moreover, we prove the following result (see Subsection \ref{counter-BL}).

\begin{proposition}
Let $1<p<2$ and $\mu = C e^{-\frac{1}{p}\|x\|^p_p} dx$. Then the best value of $\lambda$ in inequality (\ref{BL1008}) satisfies $\lambda > 1 - \frac{1}{p}$. In particular, inequality (\ref{main-question-intro}) fails to hold in this case.
\end{proposition}

\begin{remark}
Let us stress  the presence of a big discrepancy between the unconditional and symmetric cases. Indeed, according to Theorem \ref{main-thm} inequality (\ref{main-question-intro}) holds for $\mu = C e^{-\frac{1}{p}\|x\|^p_p} dx$ and  {\bf all $p >1 $} on the set of unconditional functions.
\end{remark}

To conclude the discussion on the  strong Brascamp--Lieb inequality, let us describe the main steps in the proof of Theorem \ref{strongBL}. First we make the change of variables pushing forward the measure $\mu$ into a measure of the form
$$
\prod_{i=1}^n |y_i|^{\frac{2}{p}-1} \cdot m_0,
$$
where $m_0$ is a rotationally invariant measure. Then using homogeneity we show that our problem is equivalent to the spectral gap problem for the following operator on $\mathbb{S}^{n-1}$:
$$
Lf = \Delta_{\mathbb{S}^{n-1}} f + \Bigl( \frac{2}{p}-1\Bigr) \langle \omega, \nabla_{\mathbb{S}^{n-1}} f
\rangle,
$$
where $\Delta_{\mathbb{S}^{n-1}}, \nabla_{\mathbb{S}^{n-1}}$
are the spherical Laplacian and the spherical gradient, respectively, and
$$
\omega = \Bigl( \frac{1}{y_1}, \cdots, \frac{1}{y_n}\Bigr).
$$
Note that a complete orthogonal system of eigenfunctions for $L$ contains non-elementary functions. This is true even for $n=2$, in this case the eigenfunctions (after appropriate change of variables) belong to the family of  Legendre functions, which are non-elementary in general.
Fortunately, it turns out that the eigenfunctions we need for establishing the sharp spectral gap estimate on the set of even functions are elementary and the expressions  for corresponding eigenvalues have simple algebraic form. 

\medskip

The paper is organized as follows. Section 2 contains preliminary material. In Section 3 we study finiteness and continuity of the Blaschke--Santal\'o functional. Section 4 is devoted to the symmetrization approach and symmetric properties of the maximizers, and we prove Theorem \ref{main-thm} (Theorems \ref{BLunconditional-p-homo} and \ref{thA}). In Section 5 we discuss the existence of maximizers for the Blaschke--Santal\'o functional, we present some counterexamples. In Section 6 we present the reduction to the convex body case. In Section 7 we outline the mass transport approach to the Blaschke-Santal\'o inequality and prove Theorem \ref{MA-homogen-solutions} (Theorem \ref{p-hom-solutions}). Finally, in Section 8 we deal with the strong Brascamp-Lieb type inequalities and we prove Theorem \ref{strongBL} (Theorem \ref{BL-homog-powers}).

\medskip

\noindent\textbf{Acknowledgements.}  The first author was supported by the project {\em Disuguaglianze analitiche e geometriche}, funded by the Gruppo per Analisi Matematica la Probabilit\`a e le loro Applicazioni. 
The second author gratefully acknowledges the support of the Basic Research Program of HSE University.
The third named author is supported by NSF DMS-1753260. The fourth named author is supported by ISF grant no. 2574/24. The third and fourth authors are supported by NSF-BSF DMS-2247834. The authors are grateful to the workshop ``Geometric inequalities, Convexity and Probability'' at BIRS IMAG in Granada, Spain in June 2023. We thank
Emanuel Milman for providing us with references and anonymous referees for careful reading and numerous suggestions. 

\section{Preliminaries}
\label{prem}

As a rule we will omit the domain of integration of an integral if this domain is entire $\R^n$:
$$
\int f dx := \int_{\R^n} f dx.
$$

\subsection{Convex bodies}
By a {\em convex body}, we mean a convex and compact subset of $\R^n$ with non empty interior. Our main reference for properties of convex bodies is the monograph \cite{book4}. Given a convex body $K$ we shall consider its support function $h_K:\R^n\rightarrow\R$ defined by
$$
h_K(x)=\max_{y\in K}\langle x,y\rangle.
$$
We will always assume that $0$ is contained in $K$.
Then $h_K$ is nonnegative, convex and positively $1$-homogeneous. The latter means the following property: $h_K(tx) = t h_K(x)$ for all $t \ge 0$. 

The radial function $\rho_K$ of a convex body $K$ containing the origin, is defined, for $x\in\R^n$, by
$$
\rho_K(x)=\sup\{t>0:\, tx\in K\}.
$$
If $K$ is a compact, convex set with a non-empty interior, then $\rho_K=h^{-1}_{K^o}$. 

We also recall the Minkowski functional of an origin symmetric convex body $K$: 
$$
\|x\|_K=\rho^{-1}_K(x)\quad\forall\ x\in\R^n.
$$

\subsection{Convex functions}
We will consider convex functions $\Phi\colon\R^n\to\R\cup\{+\infty\}$. The space of these functions will be denoted by $\fconvx$. Our general references on convex functions are the monographs \cite{Rockafellar} and \cite{Rockafellar-Wets}. Given a convex function $\Phi$ we define its domain as
$$
\dom(\Phi)=\{x\in\R^n\colon \Phi(x)<+\infty\}.
$$
A convex function $\Phi$ is said {\em proper} if its domain is not empty, and {\em coercive} if 
$$
\lim_{|x|\to\infty}\Phi(x)=+\infty.
$$
On the set of convex functions we fix the topology induced by epi-convergence. A sequence $\Phi_k$, $k\in\N$, epi-converges to $\Phi$ if:
\begin{itemize}
    \item[(i)]
    $$
    \liminf_{k\to+\infty}\Phi_k(x)\ge \Phi(x),
    $$
    for every $x\in\R^n$;
    \item[(ii)] for every $x\in\R^n$ there exists a sequence $x_k$, $k\in\N$, converging to $x$ and such that
    $$
    \liminf_{k\to+\infty}\Phi_k(x_k)=\Phi(x). 
    $$
\end{itemize}

We note that for sequences of finite convex functions, epi-convergence to a finite convex function is 
equivalent to uniform convergence on compact subsets of $\R^n$ -- 
see e.g. \cite[Theorem 7.17]{Rockafellar-Wets}.

We say that a (not necessary convex) function $\Phi$ is $\alpha$-homogeneous with $\alpha > 0$ if
$$
\Phi(tx) = t^{\alpha}\Phi(x)
$$
for all $x \in \mathbb{R}^n$ and $t \ge 0$.

Finally, we say that a (not necessary convex) function $f$ is unconditional with respect to a orthogonal basis $\{e_i\}$ if
$$
f(\pm x_1, \cdots,\pm x_n) = f(x_1, \cdots,x_n).
$$
A set $A$ is unconditional if $I_A$ is unconditional.

\subsection{The Legendre transform}\label{Legtr}
Given a lower semi-continuous convex function $\Phi\colon\R^n\to\R\cup\{+\infty\}$, $\Phi\not\equiv+\infty$, we denote by $\Phi^*$ its conjugate, or Legendre transform, which is defined as follows:
$$
\Phi^*(y)=\sup_{x\in\R^n} \bigl(\langle x,y\rangle-\Phi(x) \bigr),\quad\forall\ y\in\R^n.
$$
We collect some properties of the Legendre transform that will be used in this paper.

\begin{proposition} \label{legendre-basic} The following properties hold for lower semi-continuous convex functions on $\R^n$:
\begin{itemize}
    \item  \cite[Theorem 11.1]{Rockafellar-Wets} $\Phi^\ast$ is also a lower semi-continuous convex function and $\Phi^{\ast\ast} = \Phi$.
    \item \cite[Theorem 11.8(d)]{Rockafellar-Wets} $\Phi$ is finite everywhere if and only if $\Phi^*$ is a super-coercive convex function, that is:
    $$
    \lim_{|x|\to\infty}\frac{\Phi^*(x)}{|x|}=\infty.
    $$
    \item \cite[Theorem 11.34]{Rockafellar-Wets} A sequence $\Phi_k$, $k\in\N$, epi-converges to $\Phi$ if and only if the sequence $\Phi_k^*$ epi-converges to $\Phi^*$. 
\end{itemize}
\end{proposition}

\begin{remark}\label{invariance} Let $A\colon\R\to\R$ be an invertible linear map, and $b\in\R$. Given a finite convex function $\Phi$, consider the function $\bar \Phi\colon\R^n\to\R$ defined by
$$
\bar \Phi(x)=\Phi(Ax)+b.
$$
Clearly $\bar \Phi$ is a finite convex function as well. A direct computation shows that
$$
\bar \Phi^*(y)=\Phi^*(A^{-T}y)-b, \quad y\in\R^n,
$$
where $A^{-T}$ is the inverse of the transpose of $A$.
\end{remark}

\medskip

Further important properties of the Legendre transform are contained in the following statement.

\begin{proposition}[Legendre transform of smooth functions]\label{legendre-prop}
Let $V\in C^2(\R^n)$ be such that $D^2 V(x)$ is positive definite for every $x$. Then the following properties hold, for every $x$.
\begin{enumerate}
\item $V(x) + V^*(\nabla V(x)) = \ip{x}{\nabla V(x)},$    
\item $\nabla V(\nabla V^* (x)) = x$, in other words $\nabla V \circ \nabla V^* = \text{Id}.$
\item $\nabla^2 V^*(\nabla V(x)) = (\nabla^2 V)^{-1}(x)$.
\end{enumerate}
\end{proposition}

We will use throughout the following well-known 
formula:
if
$$
V(x) = \frac{1}{p} \|x\|^p_p = \frac{1}{p} 
 \sum_{i=1}^n|x_i|^p,
$$
then 
$$
V^*(y) = \frac{1}{p^*} \|y\|^{p^*}_{p^*} = \frac{1}{p^*}
 \sum_{i=1}^n|y_i|^{p^*},
$$
where $p^* = \frac{p}{p-1}$. In addition, for every $x$ with non-zero coordinates 
$
D^2 V(x)  
$
is a diagonal matrix with values
$(p-1) |x_i|^{p-2}$ on the diagonal. Similarly for $V^*$.
In particular,
$$
\det D^2 V = (p-1)^n \prod_{i=1}^n |x_i|^{p-2} , \ \  \det D^2 V^* = (p^*-1)^n \prod_{i=1}^n |x_i|^{p^*-2}. 
$$
More generally, we will apply these formulae to  $\tilde{V}(x) = c V(x)$. Then
$
\tilde{V}^*(y) = {c} V^*(\frac{y}{c}),
$
$D^2 \tilde{V}^*(y)  = \frac{1}{c} D^2 V^*(y/c)$
and $\det D^2 \tilde{V}^*(y) = \frac{1}{c^n} \det D^2 V^*(y/c).$

\subsection{Optimal transportation}
Let us consider two probability measures $\mu$ and $\nu$ on $\mathbb{R}^n$ with finite second moments. We assume that both measures are absolutely continuous with respect to the Lebesgue measure, and we denote their respective densities by $\rho_{\mu}$ and $\rho_{\nu}$.
 According to the celebrated Brenier theorem (see \cite{Villani}, \cite{BKS})
there exists a lower semi-continuous convex function $U$ such that $\nabla U$
pushes forward $\mu$ onto $\nu$:
$$
\int f(\nabla U) d\mu = \int f d\nu
$$
for any test function $f$.

The mapping $T : x \to \nabla U(x)$ is known as the optimal transportation mapping. Note that $T$ is well defined almost everywhere with respect to the Lebesgue measure, since $U$ is almost everywhere differentiable as a convex function.

The optimal transport mapping $\nabla U$ is $\mu$-a.e. unique, meaning that if $T_1 =\nabla U_1$ and $T_2 = \nabla U_2$ are
pushing forward $\mu$ onto $\nu$ and $U_1, U_2$ are convex, then
$$
T_1 = T_2
$$
$\mu$-almost everywhere.

The regularity of $U$ is in general a difficult issue. Results  about H{\"o}lder continuity of the optimal transportation can be found in \cite{Figalli} (see Section 4.6).  Despite the fact that many regularity results are known to experts, it is usually hard to find a precise reference for a specific situation.    In this paper we use more refined regularity results in Section 7, where we refer to work \cite{KK-eigenvalues}.
We will give more details in Section 7. 

For many purposes it is sufficient to know only the validity of the following change of variables formula, which holds in the non-smooth setting (see \cite{Villani} and \cite{McCann} for explanations):
\begin{equation}
\label{cvf}
\rho_{\mu}(x) = \rho_{\nu}(\nabla U(x)) \det D^2_a U(x),
\end{equation}
where $D^2_a U$ is the absolutely continuous part of the distributional Hessian $D^2 U$. In particular, $D^2_a U$ is a symmetric and positive semi-definite matrix. Equation \eqref{cvf} holds for $\mu$-almost all $x$. 

Another important instrument related to optimal transportation is the so called Hessian metric. Assume that $\rho_{\mu}, \rho_{\nu}$ are smooth and positive and $U$ is smooth and strictly convex. We  consider the following  Riemannian metric
$$
g(x) = D^2 U(x),
$$
and the corresponding Dirichlet form
$$
\mathcal{E}(f) =\int \langle (D^2 U)^{-1} \nabla f, \nabla f \rangle d\mu.
$$
The generator of $\mathcal{E}$
\begin{equation}
    \label{generator}
L f = {\rm Tr} \bigl[ (D^2 U^{-1} D^2 f \bigr] - \langle \nabla f, \nabla W(\nabla U) \rangle
\end{equation}
is a second-order elliptic differential operator, naturally related to the metric-measure space $\mathbb{R}^n$, equipped with metric $g$ and measure $\mu$.
Here $W$ is the potential of $\nu$: $\nu = e^{-W}dy$.
$L$ is symmetric with respect to $\mu$: if $f, g $ are smooth and supported on compact sets lying inside of ${\rm supp}(\mu)$, then
$$
- \int L f g d \mu = \int \langle (D^2 U)^{-1} \nabla f, \nabla g \rangle d \mu.
$$
This metric and its applications to convex geometry has been studied in 
\cite{Kolesnikov2014}, \cite{KM}, \cite{Klartag-pmm}, 
\cite{Klartag-lcmm}, \cite{KK-eigenvalues},
\cite{KK-curvature}, \cite{KK-extremal},
\cite{Kolesnikov-sphere},
 \cite{Caglar-Kolesnikov-Werner}.
Its counterpart on the sphere together with related elliptic operator (Hilbert operator) is a natural 
instrument for studying Minkowski-type problems (see  \cite{KM-LpBMproblem}, \cite{Kolesnikov-sphere}, \cite{Milman-cageometry}, \cite{Milman-cainequality}, \cite{MI1}, \cite{MI2}).

\subsection{Measures}

The image of a measure $\mu$ under a measurable
mapping $T$ is denoted by $\mu \circ T^{-1}$.

We say that a probability measure $\mu$ on $\mathbb{R}^n$ is isotropic if the coordinate functions $x_i$, $1 \le i \le N$ do
satisfy
$$
\int_{\mathbb{R}^n} x_i d\mu =0, \ 1 \le i \le N,
$$
$$
\int_{\mathbb{R}^n} x_i x_j d\mu = \delta_{ij}, \ 1 \le i,j \le N, 
$$
where $\delta_{ij}$ are the usual Kronecker symbols. For every absolutely continuous probability measure  $\mu$ there exists a non-degenerated linear mapping $A \colon \mathbb{R}^n \to \mathbb{R}^n$ such that the image of $\mu$ under $A$ is isotropic.

A probability measure is called log-concave if for all compact sets $A,B$ and a number $\lambda \in [0,1]$ one has
$$
\mu^{1-\lambda}(A) \mu^{\lambda}(B)
 \le \mu ((1-\lambda) A + \lambda B).$$
 Here $+$ means Minkowski addition of sets, i.e., $A+B = \{x+y: \ x\in A, y \in B\}$
According to a celebrated result of Borell \cite{Borell} every absolutely continuous log-concave measure has the form $d\mu = e^{-\Phi}dx$ for some convex function $\Phi$.

\section{Finiteness and continuity conditions for the Blaschke-Santal\'o functional}
\label{section continuity}

Let $\mu_1$ and $\mu_2$ be non-negative Borel measures on $\R^n$. We will assume that $\mu_1$ and $\mu_2$ are absolutely continuous with respect to the Lebesgue measure, and denote by $\rho_1$ and $\rho_2$ their respective densities. We study the functional
$$
\bs_{ \alpha, \beta, \rho_1, \rho_2}(\Phi) =
\Bigl( \int e^{-\alpha \Phi} \rho_1 dx\Bigr)^{\frac{1}{\alpha}} \Bigl(\int e^{-\beta \Phi^*} \rho_2 dy \Bigr)^{\frac{1}{\beta}}, 
$$
where $\alpha, \beta$ are positive numbers.

\begin{definition}\label{class-C} We define $\class$ as the class of lower semi-continuous, even, convex functions $\Phi\colon\R^n\to\R\cup\{+\infty\}$ such that:
\begin{enumerate}
\item \label{C-int}$\inte(\dom(\Phi))\ne\emptyset$,
where ``$\inte$'' denotes the interior;
\item \label{C-coercive} $\lim_{|x|\to\infty}\Phi(x)=+\infty$.
\end{enumerate}
\end{definition}

\begin{remark} 
Since $\Phi$ is even, property (\ref{C-int}) is equivalent to $0 \in \inte(\dom(\Phi))$. Moreover, by \cite[Corollary 3.27]{Rockafellar-Wets}, property (\ref{C-coercive}) is equivalent to the a priori stronger statement that $\Phi(x) \ge a|x| + b$ for some $a>0$ and $b \in \RR$. It follows from \cite[Theorem 11.8(c)]{Rockafellar-Wets} that $\Phi \in \class$ if and only if $\Phi^\ast \in \class$.
\end{remark}

\begin{definition} Let $\mu$ be a Borel measure on $\R^n$. We say that $\mu$ is admissible if
\begin{enumerate}
\item $\mu$ is absolutely continuous with respect to the Lebesgue measure on $\R^n$, and its density $\rho$ is positive on $\R^n$;
\item there exist positive constants $A, B, p$ such that
$$
\rho(x)\le A+B|x|^p
$$
for every $x\in\R^n$.
\end{enumerate}
\end{definition}

\begin{proposition} Let $\mu_1$ and $\mu_2$ be admissible measures, with density $\rho_1$ and $\rho_2$ respectively. Then for every $\Phi\in\class$,
$$
0<\int_{\R^n}e^{-\Phi}d\mu_1, \int_{\R^n}e^{-\Phi^*}d\mu_2<+\infty.
$$
In particular
$$
0<\bs_{\alpha,\beta,\rho_1,\rho_2}<+\infty.
$$
\end{proposition}

\begin{proof} As $\Phi\in\class$, its domain contains a neighborhood of the origin; hence $e^{-\Phi}$ is continuous (by the continuity of $\Phi$ in the interior of its domain) and strictly positive in a neighborhood of $0$. As $\rho_1$ is  positive everywhere, we obtain
$
0<\int_{\R^n}e^{-\Phi}d\mu_1.
$
The condition 
$
0<\int_{\R^n}e^{-\Phi^*}d\mu_2
$
is obtained via the same argument, as $\Phi^*\in\class$ and $\mu_2$ is admissible.

As $\Phi\in\class$,
$
\lim_{|x|\to\infty}\Phi(x)=+\infty.
$
This implies that there exist $a\in\R$ and $b>0$ such that
$
\Phi(x)\ge a+b|x|
$
for every $x\in\R^n$ (see for instance \cite[Lemma 8]{Colesanti-Ludwig-Mussnig}). By the growth condition satisfied by $\rho_1$, we get
$
\int_{\R^n}e^{-\Phi}d\mu_1<+\infty.
$
In a similar way
$
\int_{\R^n}e^{-\Phi^*}d\mu_2<+\infty
$
can be proved.
\end{proof}

We will now show a continuity property of our functional.

\begin{proposition}\label{continuity} Let $\mu_1$ and $\mu_2$ be admissible measure with density $\rho_1$ and $\rho_2$, respectively. Let $\Phi$, $\Phi_k$, $k\in\N$, belong to $\class$. Assume that $\Phi_k$ epi-converges to $\Phi$. Then
$$
\lim_{k\to\infty}\bs_{\alpha, \beta,\rho_1,\rho_2}(\Phi_k)=\bs_{\alpha, \beta,\rho_1,\rho_2}(\Phi).
$$
\end{proposition}

\begin{proof} As $\Phi_k$, $k\in\N$, and $\Phi$ are coercive, and $\Phi_k$ epi-converges to $\Phi$, there exist $a>0$ and $b\in\R$ such that
\begin{equation}\label{uniform bound}
\Phi_k(x)\ge a|x|+b    
\end{equation}
for every $x\in\R^n$ and for every $k\in\N$, and the same property is verified by $\Phi$ (see for instance \cite[Lemma 8]{Colesanti-Ludwig-Mussnig}). 
Moreover, as $\Phi_k$ epi-converges to $\Phi$, $\Phi_k(x)$ converges to $\Phi(x)$ for every $x$ in the interior of $\dom(\Phi)$. On the other hand, if $x\notin\dom(\Phi)$, then 
$$
\lim_{k\to+\infty} \Phi_k(x)=+\infty=\Phi(x).
$$
As the boundary of $\dom(\Phi^*)$ has zero Lebesgue measure, we conclude that 
$$
\lim_{k\to+\infty}e^{-\Phi_k(x)}=e^{-\Phi(x)}
$$
for almost every $x\in\R^n$. Next\eqref{uniform bound} justifies the use of the dominated convergence theorem, which yields
$$
\lim_{k\to+\infty}\int_{\R^n}e^{-\Phi_k} d\mu_1 =\int_{\R^n}e^{-\Phi} d\mu_1.
$$
Next, note that $\Phi^*_k$ epi-converges to $\Phi^*$; moreover $\Phi^*$ and $\Phi_k^*$, $k\in\N$, belong to $\class$. Therefore, we may repeat the same considerations of the previous part of this proof, and deduce that
$$
\lim_{k\to+\infty}\int_{\R^n}e^{-\Phi_k^*} d\mu_2 =\int_{\R^n}e^{-\Phi^*} d\mu_2.
$$
Finally, as $\Phi$ and $\Phi^*$ belong to $\class$,
$$
0< \int_{\R^n}e^{-\Phi} d\mu_1, \int_{\R^n}e^{-\Phi^*} d\mu_2 <\infty.
$$
This concludes the proof.
\end{proof}

\section{Symmetrization}\label{section-symmetrization}

The symmetrization technique is the main tool for proving inequality of the Blaschke--Santal\'o type; it is therefore not surprising that they can be used also in the functional version of this result.

To symmetrize a function $\Phi$ we apply Steiner symmetrization to its level sets. It is a standard observation that the symmetral $\Phi_H$ (where $H$ is the hyperplane with respect to which we symmetrize) has the same distribution as $\Phi$, thus the value of the integral $\int e^{-\Phi}dx$ is preserved under symmetrization. Moreover, we will show that under some natural assumptions on density $\rho_1$
($\rho_1$ must be in a sense decreasing) the value of
$$
\int e^{-\Phi} \rho_1 dx
$$
is increasing under symmetrization. In addition, the value of 
$$
\int e^{-\Phi^*} \rho_2 dx
$$
is increasing under symmetrization as well, if $\rho_2$ is log-concave and admits appropriate symmetries.

\medskip

Let $H$ be a hyperplane of $\R^n$, passing through the origin. Given a convex set $C\subset\R^n$, we denote by $C_H$ the Steiner symmetral of $C$, with respect to $H$, that is,
$$C_H=\left\{x+y:\,x\in H,\,\, y\perp H,\,\, |y|\leq \frac{b(x)-a(x)}{2}\right\},$$
where for $x\in H$ we define $[a(x),b(x)]=C\cap \{x+y,y\perp H\}$, the interval formed by the intersection of $C$ with the line orthogonal to $H$ at the point $x\in H.$

Let $\Phi$ be a convex and coercive function defined in $\R^n$. We denote by $\Phi_H$ the Steiner symmetral of $\Phi$ with respect to $H$. One way of defining $\Phi_H$ is through its sublevel sets:
$$
\{x\colon\Phi_H(x)\le s\}=(\{x\colon\Phi(x)\le s\})_H, \quad\forall\, s\in\R
$$
(with the convention that the Steiner symmetral of the empty set is the empty set).

\subsection{Monotonicity results for $\int_{\R^n}e^{-\Phi^*}d\mu$}

We discuss some classical properties of Steiner symmetrizations for functions. For more information, see \cite{almut}. The main result of this part is the following theorem.

\begin{theorem}\label{symmetrization} Let $H$ be an hyperplane passing through the origin in $\R^n$. Let $\mu = \rho dx$ be a log-concave measure in $\R^n$. Assume that 
\begin{equation}\label{density symmetry}
\rho(t e+y)=\rho(te-y)
\end{equation}
for every $t\in\R$ and $y\in H$, where $e$ is a normal unit vector to $H$. Then for every even, proper, coercive convex function $\Phi$ defined on $\R^n$, we have
$$
\int_{\R^n}e^{-\Phi^*}d\mu\le\int_{\R^n}e^{-(\Phi_H)^*}d\mu.
$$
\end{theorem}

We proceed with some corollaries. The next statements follow from Theorem \ref{symmetrization} applied to the Lebesgue measure, and more generally to radially symmetric log-concave measures, and to measures with unconditional density.

\begin{corollary}\label{Lebesgue} Let $\Phi\colon\R^n\to(-\infty,\infty]$ be convex, even and coercive. Then, for every hyperplane $H$,
$$
\int_{\R^n}e^{-\Phi^*(x)}dx\le\int_{\R^n}e^{-(\Phi_H)^*(x)}dx.
$$
\end{corollary}

If the measure $\mu$ has  a radially symmetric density $\mu=g(|x|)dx$, then it satisfies (\ref{density symmetry}) because
$$
g(te+y) = g(\sqrt{y^2 +t^2}) = g(te-y).
$$
Thus we get the following corollary.

\begin{corollary}\label{radially symmetric} Let $\mu$ be a log-concave measure on $\R^n$, with a  radially symmetric density with respect to the Lebesgue measure. Let $\Phi\colon\R^n\to(-\infty,\infty]$ be convex, even and coercive. Then, for every hyperplane $H$,
$$
\int_{\R^n}e^{-\Phi^*}d\mu\le\int_{\R^n}e^{-(\Phi_H)^*}d\mu.
$$
\end{corollary}

For the next result, we say that a  function $\rho\colon\R^n\to\R$ is unconditional, if
$$
\rho(x_1,\dots,x_n)=\rho(\pm x_1,\dots,\pm x_n)
$$
for every choice of the signs $+$ and $-$ on the right hand side. This is equivalent to saying that the graph of $\rho$ is symmetric with respect to each coordinate hyperplane. Note that, if $\rho$ is unconditional and $H$ is a coordinate hyperplane, then condition \eqref{density symmetry} is verified.

\begin{corollary}\label{unconditional} Let $\mu$ be a log-concave measure on $\R^n$, with unconditional density with respect to the Lebesgue measure. Let $\Phi\colon\R^n\to(-\infty,\infty]$ be convex, even and coercive. Then, for every coordinate hyperplane $H$,
$$
\int_{\R^n}e^{-\Phi^*}d\mu\le\int_{\R^n}e^{-(\Phi_H)^*}d\mu.
$$
\end{corollary}

\medskip

\subsubsection{Proof of Theorem \ref{symmetrization}}
The idea of the proof of Theorem \ref{symmetrization} is inspired by the argument of Meyer and Pajor in \cite{MePa}.

Given a function $f\colon\R^n\to(-\infty,+\infty]$, we denote by $\epi(f)$ its epigraph:
$$
\epi(f)=\{(x,z)\in\R^{n+1}\colon z\ge f(x)\}.
$$

\begin{lemma}\label{conjugate epigraph} Let $\Phi$ be a proper convex function defined in $\R^n$. Then
$$
\epi(\Phi^*)=\{(y,w)\in\R^n\times\R\colon w+z\ge\langle x,y\rangle,\,\forall\,(x,z)\in\epi(\Phi)\}.
$$
\end{lemma}

\begin{proof} Let 
$$
X=\{(y,w)\in\R^n\times \R\colon w+z\ge\langle x,y\rangle,\,\forall\,(x,z)\in\epi(\Phi)\}.
$$
Let $(y,w)\in\epi(\Phi^*)$; then
$$
w\ge\Phi^*(y)=\sup_{x}\langle y,x\rangle-\Phi(x),
$$
Hence
$$
w+\Phi(x)\ge\langle y,x\rangle\quad\forall\, x.
$$
If $z\ge \Phi(x)$, then
$$
w+z\ge\langle y,x\rangle.
$$
This proves that $(y,w)\in X$. Assume now that $(y,z)\in X$. Then, for every $x\in\R^n$,
$$
w+\Phi(x)\ge\langle y,x\rangle,
$$
so that
$$
w\ge\sup_{x}\langle y,x\rangle-\Phi(x)=\Phi^*(y).
$$
This proves that $(y,z)\in\epi(\Phi^*)$. 
\end{proof}

In the sequel, we choose a coordinate system so that
$$
H=\{x=(x_1,\dots,x_n)\in\R^n\colon x_n=0\}.
$$
The points of $\R^n$ will be written in the form
$$
(X,x)\in H\times\R\quad\mbox{or}\quad(Y,y)\in H\times\R.
$$
Similarly, the points of $\R^{n+1}$ will be written as
$$
(X,x,z)\in H\times\R\times\R\quad\mbox{or}\quad(Y,y,w)\in H\times\R\times\R.
$$
We also set
$$
H'=\{(X,0,z)\colon X\in H,\, z\in\R\} = H \times \R \subset\R^{n+1}.
$$
For a convex set $A \subseteq \R^{n+1}$ and $z \in \R$, write 
$A^z = \{ (X,x) \in \R^n\ :\ (X,x,z) \in A \} \subseteq \R^n $ for the sections of $A$. Note that Steiner 
symmetrization with respect to $H'$ acts independently on each section $\R^n \times \{z\}$, and therefore
$ (A_{H'})^z = (A^z)_H$. In particular for a convex function $\Phi$ we may take $A = \epi(\Phi)$ and obtain 
\begin{align*}
[\epi(\Phi)_{H'}]^z &= [(\epi(\Phi))^z]_H = \left[\{ (X,x) \in \R^n \ :\ \Phi(X,x) \le z \}\right]_H \\
&= \left[\{ (X,x) \in \R^n \ :\ \Phi_H(X,x) \le z \}\right] = [\epi(\Phi_H)]^z,
\end{align*}
so $\epi(\Phi)_{H'} = \epi(\Phi_H)$. 

Given a function $\Phi$ as in the statement of Theorem \ref{symmetrization}, we denote by $P_\Phi$ the orthogonal projection of $\epi(\Phi)$ onto $H'$.  We then have
$$
\epi(\Phi_H)=\left\{(X,x,z)\colon (X,z)\in P_\Phi,\, x=\frac{x_2-x_1}{2},\, (X,x_1,z), (X,x_2,z)\in\epi(\Phi)\right\}.
$$
From Lemma \ref{conjugate epigraph} we know that
$$
\epi(\Phi^*)=\{(Y,y,w)\colon z+w\ge\langle X,Y\rangle+xy,\,\forall\,(X,x,z)\in\epi(\Phi)\}.
$$
Clearly here $\langle X,Y\rangle$ denotes the scalar product in $H=\R^{n-1}$. Moreover
\begin{equation}
\label{give-ref-later}
\epi((\Phi_H)^*)=\left\{(Y,y,w)\colon z+w\ge\langle X,Y\rangle+y\left(\frac{x_2-x_1}{2}\right),\
\forall\,(X,x_i,z)\in\epi(\Phi),\, i=1,2\right\}.
\end{equation}

For a general set $A\subset\R^{n+1}$, and $y\in\R$, we set
$$
A(y)=\{(Y,0,w)\in H'\colon (Y,y,w)\in A\}
$$
for the section of $A$ parallel to $H'$. 

\begin{proposition}\label{proposition MP} In the previous notations, for every $y\in\R$:
$$
\frac12\epi(\Phi^*)(y)+\frac12\epi(\Phi^*)(-y)\subset\epi((\Phi_H)^*)(y).
$$
\end{proposition}
\begin{proof} Let
$$
(Y',w')\in\epi(\Phi^*)(y),\quad (Y'',w'')\in\epi(\Phi^*)(-y),
$$
and let
$$
(X,x_1,z),\, (X,x_2,z)\in\epi(\Phi).
$$
Then
$$
w'+z\ge\langle Y',X\rangle+yx_2,\quad w''+z\ge\langle Y'',X\rangle-yx_1,
$$
whence
$$
\left(\frac{w'+w''}{2}\right) + z\ge \Big\langle\left(\frac{Y'+Y''}{2}\right),X\Big\rangle+y \left(\frac{x_2-x_1}{2}\right).
$$
By (\ref{give-ref-later}), it follows that
$$
\left(\frac{Y'+Y''}{2},\frac{w'+w''}{2}\right)\in\epi((\Phi_H)^*)(y).
$$
\end{proof}

Let $g$ be a function defined in $\R^n=H\times\R$; for every $y\in\R$ we denote by $g_y$ the function defined on $H$ by
$$
g_y(Y)=g(Y,y).
$$
Clearly
$$
\epi(g_y)=(\epi(g))(y).
$$

\begin{proof}[Proof of Theorem \ref{symmetrization}]
By Proposition \ref{proposition MP}, we have
$$
\frac12\epi((\Phi^*)_y)+\frac12\epi((\Phi^*)_{-y})\subset\epi(((\Phi_H)^*)_y)
$$
for every $y\in\R$. Now let $\bar\Phi^*$ be the function defined by:
$$
\bar\Phi^*(Y,y)=\Phi^*(Y,-y).
$$
Note that as $\Phi$ is even, $\Phi^*$ and $\bar\Phi^*$ are even as well. 

Given $(Y,w)\in H'$, we have $(Y,w)\in\epi((\Phi^*)_{-y})$ if and only if
$$
w\ge\Phi^*(Y,-y)=\Phi^*(-Y,y)=\bar\Phi(Y,y),
$$
that is, if and only if $(Y,w)\in\epi((\bar\Phi^*)_y)$. Therefore
\begin{equation}\label{epi-sum}
\frac12\epi((\Phi^*)_y)+\frac12\epi((\bar\Phi^*)_{y})\subset\epi(((\Phi_H)^*)_y)\quad\forall\, y\in\R.
\end{equation}
The set on the left hand side of the previous relation is the graph of the function
\begin{eqnarray*}
H\ni Y\,&\mapsto\,&\sup\Bigl\{\frac12\Phi^*(Y_1,y)+\frac12\bar\Phi^*(Y_2,y)\colon\frac{Y_1+Y_2}{2}=Y\Bigr\}\\
&=&\frac12\cdot(\Phi^*)_y\,\square\,\frac12\cdot(\bar\Phi^*)_y,
\end{eqnarray*}
where $\square$ denotes the sup convolution operation, and $\cdot$ the corresponding product by non-negative coefficients (see \cite{Rockafellar}, Section 16]). 

Hence, by (\ref{epi-sum}),
$$
((\Phi_H)^*)_y\le\frac12\cdot(\Phi^*)_y\,\square\,\frac12\cdot(\bar\Phi^*)_y,
$$
so that
\begin{equation}\label{for PL}
e^{-((\Phi_H)^*)_y}\ge e^{-(\frac12\cdot(\Phi^*)_y\,\square\,\frac12\cdot(\bar\Phi^*)_y)}.
\end{equation}
We recall that $\mu$ has a density $\rho$ with respect to the Lebesgue measure; let $\mu_y$ be the measure on $H$, with density $\rho_y$ (with respect to the $(n-1)$ dimensional Lebesgue measure on $H$). Using the log-concavity of $\mu$, and then of $\mu_y$, the inequality \eqref{for PL} and the Pr\'ekopa-Leindler inequality (see Theorem 1.4.1 in \cite{AGM}), we get
$$
\int_H e^{-((\Phi_H)^*)_y}d\mu_y\ge\left(\int_H e^{-(\Phi^*)_y}d\mu_y\right)^{1/2}\,\left(\int_H e^{-(\bar\Phi^*)_y}d\mu_y\right)^{1/2}.
$$
On the other hand
\begin{eqnarray*}
\int_H e^{-(\bar\Phi^*)_y}d\rho_y&=&\int_H e^{-\Phi^*(-Y,y)}\rho(Y,y)dY\\
&=&\int_H e^{-\Phi^*(-Y,y)}\rho(-Y,y)dY\\
&=&\int_H e^{-\Phi^*(Y,y)}\rho(Y,y)dY.
\end{eqnarray*}
We conclude that
$$
\int_H e^{-((\Phi_H)^*)_y}d\mu_y\ge\int_H e^{-(\Phi^*)_y}d\mu_y
$$
for every $y$. That is
$$
\int_{\R^{n-1}} e^{-(\Phi_H)^*(Y,y)}\rho(Y,y)\ge \int_{\R^{n-1}} e^{-\Phi^*(Y,y)}\rho(Y,y),\quad\forall\, y\in\R^n.
$$
The claim of the theorem follows from Fubini's theorem.
\end{proof}

\subsection{Monotonicity results for $\int_{\R^n}e^{-\Phi}d\mu$}

The following statements are probably well-known within the area of rearrangements of functions; we include their proofs for completeness.

\begin{proposition}\label{normal rearrangement inequality general} Let $\mu$ be a measure on $\R^n$, which is absolutely continuous with respect to the Lebesgue measure, with density $\rho$. Let $H$ be a hyperplane through the origin with unit normal vector $e$. Assume that for every $x'\in H$, the function $\varphi\colon\R\to\R$ defined by
$$
\varphi(t)=\rho(x'+te)
$$
is even (in $\R$), and decreasing in $[0,\infty)$. Then, for every $\Phi\in\class$,
$$
\int_{\R^n}e^{-\Phi(x)} d\mu(x)\le\int_{\R^n}e^{-\Phi_H(x)}d\mu(x).
$$
\end{proposition}

The next two results are consequences of Proposition \ref{normal rearrangement inequality general} (the first one can be obtained also as an application of the Hardy-Littlewood inequality for decreasing rearrangements). 

\begin{proposition}\label{normal rearrangement inequality} Let $\mu$ be a measure on $\R^n$, which is absolutely continuous with respect to the Lebesgue measure, with density $\rho$ of the form
$$
\rho=\rho(x)=\varphi(|x|)
$$
where $\varphi\colon[0,+\infty)\to[0,+\infty)$ is decreasing. Let $\Phi\in\class$ and let $H$ be an hyperplane in $\R^n$, passing though the origin. Then
$$
\int_{\R^n}e^{-\Phi(x)} d\mu(x)\le\int_{\R^n}e^{-\Phi_H(x)}d\mu(x).
$$
\end{proposition}

\begin{proposition}\label{normal rearrangement inequality unconditional} Let $\mu$ be a measure on $\R^n$, which is absolutely continuous with respect to the Lebesgue measure, with unconditional density $\rho$. Assume that for every coordinate hyperplane $H$ with unit normal vector $e$, and for every $x'\in H$, the function $\varphi\colon[0,+\infty)\to\R$ defined by
$$
\varphi(t)=\rho(x'+te)
$$
is decreasing. Let $\Phi\in\class$ and let $H$ be a coordinate hyperplane. Then
$$
\int_{\R^n}e^{-\Phi(x)} d\mu(x)\le\int_{\R^n}e^{-\Phi_H(x)}d\mu(x).
$$
\end{proposition}

\subsubsection{Proof of Proposition \ref{normal rearrangement inequality general}}

We will need the following one-dimensional result.

\begin{lemma}\label{Steiner 1d} Let $\mu$ be a measure on $\R$, absolutely continuous with respect to the Lebesgue measure; assume that the density of $\mu$ is a function $\psi\colon\R\to[0,+\infty)$ which is even and decreasing in $[0,+\infty)$. Then, for every $a,b\in\R$ with $a\le b$,
$$
\mu([a,b])\le\mu\left(
\left[
-\frac{b-a}2,\frac{b-a}2
\right]
\right).
$$
\end{lemma}

\begin{proof} Let, for $x\in\R$,
\[
F(x)=
\begin{cases}
\displaystyle \int_0^x \psi(t)\,dt, & \text{if } x>0,\\[6pt]
0, & \text{if } x=0,\\[6pt]
\displaystyle -\int_0^{-x} \psi(t)\,dt, & \text{if } x<0.
\end{cases}
\]
Note that $F$ is odd in $\R$, and it is concave in $[0,+\infty)$, as $\psi$ is decreasing in $[0,+\infty)$. If
$$
a\le0\le b
$$
then 
$$
\mu([a,b])=\int_a^b\psi(t)dt=F(b)-F(a)=F(b)+F(-a).
$$
On the other hand,
\begin{eqnarray*}
\mu\left(\left[-\frac{b-a}2,\frac{b-a}2 \right]\right)&=&\int_{-\frac{b-a}2}^{\frac{b-a}2}\psi(t)dt=2\int_0^{\frac{b-a}2}\psi(t)dt\\
&=&2F\left(\frac{b-a}2\right).
\end{eqnarray*}
The inequality then follows from the concavity of $F$. The case $0\le a\le b$ can be reduced to the previous one, observing that
\begin{align*}
\mu([a,b]) & = \int_a^b \psi(t)dt = \int_a^{\frac{b+a}{2}} \psi(t)dt 
+ 
 \int_{\frac{b+a}{2}}^b\psi(t)dt =  \\& \le 2 \int_a^{\frac{b+a}{2}} \psi(t)dt 
\le 2 \int_0^{\frac{b-a}{2}} \psi(t)dt  = 
 2\mu\left(\left[0,\frac{b-a}{2}\right]\right).
\end{align*}
as $\psi$ is decreasing in $[0,+\infty)$. The case $a\le b\le 0$ is completely analogous. 
\end{proof}

\begin{lemma}\label{normal rearrangement inequality for sets} Let $\mu$ be a measure on $\R^n$, which is absolutely continuous with respect to the Lebesgue measure, with density $\rho$. Let $H$ be a hyperplane through the origin with unit normal vector $e$. Assume that for every $x'\in H$, the function $\varphi\colon\R\to\R$ defined by
$$
\varphi(t)=\rho(x'+te)
$$
is even, and decreasing in $[0,+\infty)$. Then, for every convex body $K$ in $\R^n$,
$$
\mu(K)\le\mu(K_H).
$$
\end{lemma}

\begin{proof} We may assume that 
$$
H=\{x=(x_1,\dots,x_n)=(x',x_n)\colon x_n=0\}.
$$
Let $K'$ be the orthogonal projection of $K$ onto $H$. Them, by Fubini's theorem
$$
\mu(K)=\int_{K'}\int_{K_{x'}}\rho(x',t)dt dx',
$$
where
$$
K_{x'}=\{t\in\R\colon(x',t)\in K\}.
$$
The assertion follows from Lemma \ref{Steiner 1d}, and the fact that for every $x'\in K'$ the function $t\to\rho(x',t)$ is  even and decreasing for $t\ge0$.
\end{proof}

\begin{proof}[Proof of Proposition \ref{normal rearrangement inequality general}] Use  Lemma \ref{normal rearrangement inequality for sets} and the Layer Cake Principle: for every measurable nonnegative function $f$ one has: 
$$
\int f d\mu = \int_0^\infty \mu(f >t)dt.
$$
\end{proof}

\subsection{Applications to the Blaschke-Santal\'o functional}
Let us turn back to our functional $\bs_{\alpha,\beta,\rho_1, \rho_2}$.
We will derive from the previous results some consequences on the behavior of $\bs_{\alpha,\beta,\rho_1,\rho_2}$ under the action of Steiner symmetrizations, in the radially symmetric and  unconditional cases.

\subsubsection{The radially symmetric case} Let us assume that $\rho_1,\rho_2\colon\R^n\to\R$ satisfy the following assumptions:
\begin{itemize}
\item[(R1)] $\rho_1$ is of the form  
$$
\rho=\rho(x)=\varphi(|x|)
$$
where $\varphi\colon[0,+\infty)\to[0,+\infty)$ is  decreasing;
\item[(R2)] $\rho_2$ is radially symmetric and log-concave. 
\end{itemize}

The next results follows from Corollary \ref{radially symmetric} and Proposition \ref{normal rearrangement inequality}.

\begin{proposition}\label{Steiner BS radial} Assume that $\rho_1$ and $\rho_2$ verify assumptions (R1) and (R2). Then, for every hyperplane $H$ through the origin, and for every $\Phi\in\class$
$$
\bs_{\alpha,\beta,\rho_1, \rho_2}(\Phi)\le\bs_{\alpha,\beta,\rho_1, \rho_2}(\Phi_H).
$$
\end{proposition}

\subsubsection{The unconditional case} Let us assume that $\rho_1,\rho_2\colon\R^n\to\R$ satisfy the following assumptions:
\begin{itemize}
\item[(U1)] $\rho_1$ is unconditional, and for every coordinate hyperplane $H$, with unit normal $e$, for every $x'\in H$, the function $\varphi\colon[0,+\infty)\to\R$ defined by
$$
\varphi(t)=\rho_1(x'+te)
$$
is decreasing;  
\item[(U2)] $\rho_2$ is unconditional and log-concave. 
\end{itemize}

From Corollary \ref{unconditional} and Proposition \ref{normal rearrangement inequality unconditional} we deduce the following statement.

\begin{proposition}\label{Steiner BS unconditional} Assume that $\rho_1$ and $\rho_2$ verify assumptions (U1) and (U2). Then, for every coordinate hyperplane $H$ and for every $\Phi\in\class$
$$
\bs_{\alpha,\beta,\rho_1, \rho_2}(\Phi)\le\bs_{\alpha,\beta,\rho_1, \rho_2}(\Phi_H).
$$
\end{proposition}

\subsection{Blaschke--Santal\'o inequality for unconditional functions: Theorem \ref{main-thm}(2) }\label{section-unconditional}

In this subsection we derive Blaschke--Santal\'o inequality for unconditional functions from the Pr\'ekopa--Leindler inequality. The arguments go back to  \cite{FradeliziMeyer2007}. 

For $x,y\in\R^n_+$ and for $s,t \in \mathbb{R}$ define
$$x^sy^t=(x_1^s y_1^t,...,x_n^s y^t_n).$$

\begin{lemma}
Let $V$ be a function, which is twice continuously differentiable and $p$-homogeneous  on $(0,\infty)^n$, $p \ge 1$. Then the following properties are equivalent:
    \begin{enumerate}
        \item for every $a, b \in (0,\infty)^n$ 
        $$
\langle a, \nabla V(b) \rangle \ge 
\langle a^{\frac{1}{p}} b^{\frac{p-1}{p}}, \nabla V(a^{\frac{1}{p}} b^{\frac{p-1}{p}}) \rangle = p V(a^{\frac{1}{p}} b^{\frac{p-1}{p}}).
        $$
        \item the following function is concave:
        $$
x \to V(x^{\frac{1}{p}}).
        $$
    \end{enumerate}
\end{lemma}
\begin{proof}
    Assume that 1. holds. By $p$-homogeneity $\langle b , \nabla V(b) \rangle = pV(b)$. Let $a = b + \varepsilon c $. One has:
    $$
\langle b + \varepsilon c, \nabla V(b) \rangle \ge 
p V \Bigl(\bigl(b+ \varepsilon c \bigr)^{\frac{1}{p}}b^{\frac{p-1}{p}}\Bigr)
    $$
    with equality at $\varepsilon = 0$.
    Expanding in $\varepsilon$, one gets
    $
\bigl(b+ \varepsilon c \bigr)^{\frac{1}{p}}b^{\frac{p-1}{p}} = b + 
\frac{\varepsilon}{p} c + \frac{(1-p)\varepsilon^2}{2p^2} c^2 b^{-1} + o(\varepsilon^2)
    $
    and 
    $$
p V \Bigl(\bigl(b+ \varepsilon c \bigr)^{\frac{1}{p}}b^{\frac{p-1}{p}}\Bigr)
= p V(b) + \varepsilon \langle c, \nabla V(b) \rangle + \frac{(1-p)\varepsilon^2}{2p} \Big\langle {\rm diag} \Bigl(\frac{V_{x_i}(b)}{b_i}\Bigr) c , c\Big\rangle
+ \frac{\varepsilon^2}{2p} \langle (D^2 V)(b) c,c  \rangle
+ o(\varepsilon^2).
    $$
    Finally, we obtain that $V$ must satisfy
    \begin{equation}
    \label{1pconc}
D^2 V(x) \le (p-1) {\rm diag} \Bigl(\frac{V_{x_i}(x)}{x_i} \Bigr).
    \end{equation}
    To see the equivalence  to concavity of $V_p(x) = V(x^{\frac{1}{p}})$, we note that
    \begin{align*}
D^2 V_p(x) & = \frac{1}{p^2} {\rm diag} \Bigl( x^{\frac{1-p}{p}}_i\Bigr) D^2 V(x^{\frac{1}{p}}) {\rm diag} \Bigl( x^{\frac{1-p}{p}}_i\Bigr)
+ \frac{1-p}{p^2} {\rm diag} \Bigl( x^{\frac{1-2p}{p}}_i  V_{x_i}(x^{\frac{1}{p}})\Bigr)
\\& =  \frac{1}{p^2} {\rm diag} \Bigl( x^{\frac{1-p}{p}}_i\Bigr) \Bigl[
 D^2 V(x^{\frac{1}{p}}) + (1-p)  {\rm diag} \Bigl(  \frac{V_{x_i}(x^{\frac{1}{p}})}{x_i^{1/p}}\Bigr)\Bigr]{\rm diag} \Bigl( x^{\frac{1-p}{p}}_i\Bigr).
    \end{align*}
Thus we get that $D^2 V_p \le 0$ if and only if (\ref{1pconc}) holds.

Let us assume 2. Note that function $f(a) = V(a^{\frac{1}{p}}b^{\frac{p-1}{p}})$ is a composition of $V_p$ with diagonal linear mapping ${\rm diag}(b_i^{\frac{p-1}{p}})$, hence $f$ is concave. In particular, $g(t) = f(b + t (a-b))$ is concave and consequently it satisfies $g(1) \le g(0) + g'(0)$. One can easily compute $g(1) = f(a) = V(a^{\frac{1}{p}}b^{\frac{p-1}{p}})$, $g(0) = f(b) = V(b)$, $$g'(0)=\langle a-b, \nabla f(b) \rangle = \frac{1}{p} \langle a-b, \nabla V(b) \rangle = \frac{1}{p} \langle a, \nabla V(b) \rangle- V(b) .$$ 
In the last equality we use again $p$-homogeneity of $V$: $p V(b) = \langle p, \nabla V(p) \rangle$.
This completes the proof.
\end{proof}

\begin{example}
    Let  $V = \frac{1}{p} \|x\|^p_r$. Then $V$ satisfies assumptions of the previous lemma if $r \le p$.
\end{example}

\begin{theorem}
\label{BLunconditional-p-homo}
    Let $\Phi, V$ be unconditional functions. Let, in addition, $V$ satisfy the following assumptions:
    \begin{enumerate}
        \item $V$ is $p$-homogeneous for some $p > 1;$
        \item the function $x \to V(x^{\frac{1}{p}})$
        is concave.
    \end{enumerate}
    Then
    $$
    \left( \int_{\mathbb{R}^n} e^{-\Phi} dx \right)^{\frac{1}{p}}\left(\int_{\mathbb{R}^n} e^{-\frac{1}{p-1}\Phi^*(\nabla V)} dx\right)^{\frac{p-1}{p}}
     \le  \int_{\mathbb{R}^n}  e^{-V} dx .
    $$
    If, in addition, $V$ is convex, then the inequality is sharp and $\Phi=V$ is the maximum point.
\end{theorem}
\begin{proof}
Using change of variables $x=e^r$ and Pr\'ekopa--Leindler inequality, one gets
\begin{eqnarray*}
&&\int_{\mathbb{R}^n}e^{-\Phi(x)}dx\Bigl[ \int_{\mathbb{R}^n} e^{-\frac{1}{p-1}\Phi^*(\nabla V(x))} dx \Bigr]^{p-1}=2^{np}\int_{\mathbb{R}^n_+} e^{-\Phi(x)} dx \Bigl[ \int_{\mathbb{R}^n_+} e^{-\frac{1}{p-1}\Phi^*(\nabla V(x))} dx \Bigr]^{p-1}
\\&&=2^{np}\int_{\mathbb{R}^n} e^{-\Phi(e^t)} e^{\sum_{i=1}^n t_i}dt \Bigl[ \int_{\mathbb{R}^n} e^{-\frac{1}{p-1}\Phi^*(\nabla V(e^s))} e^{\sum_{i=1}^n s_i} ds \Bigr]^{p-1}
\\&&\le 2^{np} \Bigl[\int_{\mathbb{R}^n}
  e^{- \inf_{r=\frac{t+(p-1)s}{p}}\bigl(  \frac{\Phi(e^t) + \Phi^*(\nabla V(e^s))}{p} -\frac{\sum_{i=1} t_i+(p-1)s_i}{p} \bigr)} dr\Bigr]^{p}.
\end{eqnarray*}
Apply the assumptions on $V$ and the previous lemma with $r = \frac{t}{p} + \frac{p-1}{p}s$:
\begin{eqnarray*}
\Phi(e^t) + \Phi^*(\nabla V(e^s))
&\ge& \langle e^t, \nabla V(e^s) \rangle\\
&\ge& \langle e^{\frac{t+(p-1)s}{p}}, \nabla V(e^{\frac{t+(p-1)s}{p}}) \rangle\\
&=&\langle e^r, \nabla V(e^r) \rangle\\
&=& p V(e^r).
\end{eqnarray*}
Thus
\begin{eqnarray*}
\int_{\mathbb{R}^n} e^{-\Phi(x)} dx \Bigl[ \int_{\mathbb{R}^n} e^{-\frac{1}{p-1}\Phi(\nabla V(x))} dx \Bigr]^{p-1} 
\le 2^{np}\Bigl[ \int_{\mathbb{R}^n}e^{-V(e^r) + \sum_{i=1}^n r_i } dr \Bigr]^{p}=\Bigl[\int_{\mathbb{R}^n} e^{-V(x)} dx \Bigr]^{p}.
\end{eqnarray*}
\end{proof}

\begin{corollary}
\label{1105}
    Let $p>1, r>1$, $r\le p$ and
     $V = \frac{1}{p} \|x\|^p_r$. 
    Then the Blaschke--Santal\'o inequality
    $$
\Bigl( \int e^{-\Phi} d x \Bigr)^{\frac{1}{p}} \Bigl( \int e^{-\frac{1}{p-1}\Phi^*(\nabla V(y))} dy \Bigr)^{1-\frac{1}{p}}
\le \int e^{-V} dx
    $$
    holds for unconditional function $\Phi$.
\end{corollary}

\subsection{ Theorem \ref{main-thm}: the general case}

 In this section we prove some sufficient conditions for $V$ to be the maximizer of $\bs_{p,V}$. Our proof is based on symmetrization arguments and the result of the previous section on maximization of $\bs_{p,V}$ on the set of unconditional functions.

\begin{theorem}\label{thA} Let $p>1$  and $V$ be an even convex function satisfying the following assumptions:
\begin{itemize}
        \item $V$ is $p$-homogeneous;
        \item $V$ is  unconditional and the function 
        $$
        x=(x_1,\dots x_n)\mapsto V\bigl(x^{\frac{1}{p}}_1,...,x^{\frac{1}{p}}_n\bigr)
        $$ 
        is concave in $\R^n_+$;
        \item  
        for every coordinate hyperplane $H$, with unit normal $e$, for every $x'\in H$, the function $\varphi\colon[0,+\infty)\to\R$ defined by
$$
\varphi(t)=\det D^2 V^*(x'+te)
$$
is decreasing.
\end{itemize}
Then inequality (\ref{main-question-intro}) holds for any convex even function $\Phi$.
\end{theorem}

\begin{proof}
    Let $\Phi$ be an even convex function. We observe that 
$$
\bs_{p,V}(\Phi) \le  \bs_{p,V}((\Phi)_{H_k})
$$
for every $1 \le k \le n$, where 
$(\Phi)_{H_k}$ is symmetrization of $\Phi$ with respect to the hyperplane $\{x_k=0\}$. This follows from Proposition \ref{Steiner BS unconditional}. Applying consecutively the symmetrizations $H_1, \cdots, H_n$ to $\Phi$, one obtains an unconditional function $\tilde{\Phi}$ such that
$$
\bs_{p,V}(\Phi)
\le  \bs_{p,V}(\tilde{\Phi}).
$$
On the other hand, inequality (\ref{main-question-intro}) holds for unconditional functions, according to Corollary \ref{1105}; this completes the proof.
\end{proof}

The following computations are classical and we write them for the reader's convenience.

\begin{lemma}
\label{aprconj}
Let $K$ be a  convex symmetric body, $K^{\circ}$ be its dual, and $V(x) = a \|x\|_K^p$, where $\|\cdot\|_K$ is the Minkowski functional of $K$, $p>1$ and $a>0$. Then  
$$
V^*(y) = \frac{1}{(ap)^{\frac{1}{p-1}} p^*} \ \|y\|^{p^*}_{K^{\circ}}.
    $$
    where $p^*= \frac{p}{p-1}$.
    
In particular, if  $p >1, r>1, a>0$ and $V(x) = a\|x\|^p_r$, one has
    $
V^*(y) = \frac{1}{(ap)^{\frac{1}{p-1}} p^*} \ \|y\|^{p^*}_{r^*}
    $.
\end{lemma}
\begin{proof}
Using approximation we may assume that $K$ is strictly convex with $C^2$-boundary. 
    We start with the H{\"o}lder inequality
    \begin{align*}
\langle x, y \rangle & \le \|x\|_K\|y\|_{K^{\circ}}
= \Bigl((ap)^{1/p} \|x\|_K \Bigr)\cdot \frac{\|y\|_{K^{\circ}}}{(ap)^{1/p}}
\le \frac{1}{p}\Bigl((ap)^{1/p} \|x\|_{K} \Bigr)^p
+ \frac{1}{p^*} \Bigl[ \frac{\|y\|_{K^{\circ}}}{(ap)^{1/p}}\Bigr]^{p^*}
\\& =  a\|x\|^p_K + \frac{1}{(ap)^{\frac{1}{p-1}} p^*} \|y\|^{p^*}_{K^{\circ}}.
    \end{align*}
    To complete the proof we have to show that for any $x \ne 0$ there exists another vector $\tilde{y}$ such that
    $$
\langle x, \tilde{y} \rangle =  a\|x\|^p_K + \frac{1}{(ap)^{\frac{1}{p-1}} p^*}  \|\tilde{y}\|^{p^*}_{K^{\circ}}.$$
To this end we take $
\tilde{y} = ap \|x\|_K^{p-1} \nabla \|x\|_{K}.
$
Using $1$-homogeneity of $\|\cdot \|_K$  one gets
$\langle x, \nabla \|x\|_K \rangle = \|x\|_K$, hence
$
\langle x, \tilde{y} \rangle - a\|x\|^p_K
= a(p-1)\|x\|^p_{K}.
$

By another well-known property
$\|\nabla \|x\|_K\|_{K^{\circ}}=1$. Indeed, note that
$\|x\|_{K^{\circ}} =\sup_{v : \|v\|_K = 1} \langle x, v \rangle$, and every function $f$ of the form $f(x) = \langle x, v \rangle$ satisfies $\|\nabla f\|_{K} = \| v\|_K  = 1 $; then we get that the Lipschitz constant of $\|x\|_{K^{\circ}}$ in the norm $\|\cdot \|_K$ is not greater than $1$, hence $\|\nabla \|x\|_K\|_{K^{\circ}} \le 1$. To prove the opposite inequality $\|\nabla \|x\|_K\|_{K^{\circ}} \ge 1$ we observe that
$\|x\|_K = \langle x, \nabla \|x\|_K \rangle \le \|x\| _K \| \nabla \|x\|_K\|_{K^{\circ}}$. The claim follows.

Finally, we observe that
$$
\frac{1}{(ap)^{\frac{1}{p-1}} p^*}  \|\tilde{y}\|^{p^*}_{K^{\circ}} = (ap)^{p^*-\frac{1}{p-1}} \frac{1}{p^*} \|x\|^{p^*(p-1)}_K = a(p-1) \|x\|_K^{p} = \langle x, \tilde{y} \rangle - a\|x\|^p_K,
$$
and we get the claim.
\end{proof}

\begin{theorem}
    Let $p \ge r \ge 2$ and
     $V = \frac{1}{p} \|x\|^p_r$. 
    Then the generalized Blaschke--Santal\'o inequality
    $$
\Bigl( \int e^{-\Phi} d x \Bigr)^{\frac{1}{p}} \Bigl( \int e^{-\frac{1}{p-1}\Phi^*(\nabla V(y))} dy \Bigr)^{1-\frac{1}{p}}
\le \int e^{-V} dx
    $$
    holds on the set of even functions.
\end{theorem}
\begin{proof}
   It is sufficient to check that
$
\det D^2 V^* 
$
satisfies assumption (U1) of Proposition \ref{Steiner BS unconditional}. Indeed, according to Lemma \ref{aprconj}
$$
V^*(x)= \frac{1}{p^*} \|x\|^{p^*}_{r^*},
$$
where $p^* = \frac{p}{p-1}, r^* = \frac{r}{r-1}$. Next we compute (for the sake of simplicity let $x_i>0$):
$$
\nabla V^*(x) = \|x\|^{p^*-r^*}_{r^*}
 \bigl( x_i^{r^*-1}  \bigr), 
$$
and
\begin{align*}
D^2 V^*(x) &= (p^*-r^*) \|x\|^{p^*-2r^*}_{r^*}
 \bigl( x_i^{r^*-1} x_j^{r^*-1}  \bigr) + (r^*-1) \|x\|^{p^*-r^*}_{r^*} x_i^{r^*-2} \delta_{ij}
\\&
= \|x\|^{p^*-2r^*}_{r^*} \Bigl[ (p^*-r^*)  x_i^{r^*-1} x_j^{r^*-1}  + (r^*-1) \|x\|^{r^*}_{r^*} x_i^{r^*-2} \delta_{ij} \Bigr].
\end{align*}
We set $\Lambda = \|x\|^{\frac{r^*}{2}}_{r^*} {\rm diag}(x_i^{\frac{r^*}{2}-1})$. Then
$$
D^2 V^* = \|x\|^{p^*-2r^*}_{r^*} \Lambda \Bigl[ (r^*-1) {\rm Id} + (p^*-r^*) a \otimes a\Bigr] \Lambda,
$$
where $a = \frac{1}{\|x\|_{\frac{r^*}{2}}} (x_i^{\frac{r^*}{2}})$
and
$
(a \otimes a) (u,v) = \langle a, u \rangle \langle a, v\rangle.
$
Thus one has (for all $x = (x_1, \dots, x_n)$, $x_i \in \mathbb{R}\setminus \{0\}$)
$$
\det D^2 V^* = (p^*-1) (r^*-1)^{n-1} \|x\|_{r^*}^{n(q-r^*)} \prod_{i=1}^n {|x_i|}^{r^*-2}.
$$
We obtain that the density $D^2 V^*$ satisfies assumption (U1) provided  $p^* \le r^*$ and $r^* \le 2$. Equivalently, $r \ge 2$ and $p \ge r$. 
\end{proof}

\section{Existence of maximizers}\label{section-existence-results}

In this section we present some results which guarantee the existence of maximizers of the Blaschke-Santal\'o functional $\bs_{\alpha,\beta,\rho_1,\rho_2}$, under specific conditions on $\rho_1,\rho_2,\alpha,\beta$. The picture is completed by examples, showing that existence of maximizers under too general conditions can not be expected.

In Section 7 we give a proof of the classical Blaschke--Santal{\'o} inequality (the functional form) based on the results about existence and uniqueness of the maximizers. 

\medskip

\subsection{Examples of non-existence of maximizers}\label{nonexistence}
\begin{example}
Let $\rho_1 \equiv 1$ and $\rho_2$ be $0$-homogeneous. Let $M =\max_{y \in \mathbb{S}^{n-1}} \rho_2(y)$ and assume that the restriction of $\rho_2$ to $\mathbb{S}^{n-1}$ admits exactly two maximum points: $y_0$ and $-y_0$. For simplicity, we may assume that $y_0 = e_1$. Then for every admissible $\Phi$
$$
\int e^{-\Phi} dx \int e^{-\Phi^*} \rho_2 dy < M 
\int e^{-\Phi} dx \int e^{-\Phi^*}  dy < M (2 \pi)^n .
$$
Let 
$\Phi_{\varepsilon} = \frac{1}{2} x_1^2 + \frac{\varepsilon}{2} \sum_{i=2}^n x_i^2$. 
Then
$
\int e^{-\Phi_\varepsilon} dx = (2 \pi)^{\frac{n}{2}} \varepsilon^{-\frac{n-1}{2}}
$
and
\begin{align*}
 \int e^{-\Phi_{\varepsilon}} dx \int e^{-\Phi^*_{\varepsilon}} \rho_2 dy & = (2 \pi)^{\frac{n}{2}} \varepsilon^{-\frac{n-1}{2}}
\int e^{-\frac{1}{2} y_1^2 - \frac{1}{2 \varepsilon} \sum_{i=2}^n y_i^2 } \rho_2 dy
\\&  = (2 \pi)^{\frac{n}{2}}
\int e^{-\frac{1}{2} |t|^2 } \rho_2(t_1, \varepsilon t_2, \cdots \varepsilon t_n) dt.
\end{align*}
Clearly 
$$
\lim_{\varepsilon \to 0}  \int e^{-\Phi_{\varepsilon}} dx \int e^{-\Phi^*_{\varepsilon}} \rho_2 dy  = (2 \pi)^{n} M.
$$
Thus we get that 
$$
\sup \int e^{-\Phi} dx \int e^{-\Phi^*} \rho_2 dy =(2\pi)^n M,
$$
but 
$
\Phi \to \int e^{-\Phi} dx \int e^{-\Phi^*} \rho_2 dy
$
does not attain the supremum.
\end{example}

The following example is related to Theorem \ref{existence homogeneous case0}.

\begin{example}
Without the compatibility condition $\frac{\alpha}{\beta}=\frac{n+s}{n+t}$ the functional $\bs_{\alpha,\beta,\rho_{1},\rho_{2}}$ clearly is not bounded in general, even in dimension $n=1$ (or equivalently, in dimension $n$ even if we assume all densities and functions are not only 1-symmetric but even rotation invariant). This is obvious from the proof in the previous example, but as a concrete example consider 
\[
\bs_{\alpha,\beta,\rho_{1},\rho_{2}}(\Phi)=\int e^{-\Phi}\dd x\cdot\int e^{-\Phi^{\ast}}x^{2}\dd x
\]
(i.e. $\alpha=\beta=1$, $\rho_{1}=1$, $\rho_{2}=x^{2}$, $s=0$
and $t=2$). Then $\frac{\alpha}{\beta}=1\ne\frac{1+0}{1+2}=\frac{n+s}{n+t}$.
And indeed, choosing e.g. $\Phi_{\lambda}(x)=\lambda\frac{x^{2}}{2}$
we see that 
\[
\bs_{\alpha,\beta,\rho_{1},\rho_{2}}(\Phi_{\lambda})=\int e^{-\lambda\frac{x^{2}}{2}}\dd x\cdot\int e^{-\frac{1}{\lambda}\frac{x^{2}}{2}}x^{2}\dd x=2\pi\lambda\xrightarrow{\lambda\to\infty}\infty.
\]
\end{example}

\begin{example}
If $\rho_{1}$ and $\rho_{2}$ are not assumed to be homogeneous then it is possible for $\bs_{\alpha,\beta,\rho_{1},\rho_{2}}$ to be bounded and still not attain a maximum. For instance, in the example above replace $\rho_{2}$ with $\rho_{2}(x)=e^{-\frac{1}{2}x^{2}}$. Then 
\[
\bs_{\alpha,\beta,\rho_{1},\rho_{2}}(\Phi_{\lambda})=\int e^{-\lambda\frac{x^{2}}{2}}\dd x\cdot\int e^{-\frac{1}{\lambda}\frac{x^{2}}{2}}e^{-\frac{x^{2}}{2}}\dd x=\frac{2\pi}{\sqrt{\lambda+1}}\xrightarrow{\lambda\to 0}2\pi.
\]
On the other hand, since $\rho_{2}(x)<1$ for all $x\ne0$ we have, for all $\Phi\in\class$,
\[
\bs_{\alpha,\beta,\rho_{1},\rho_{2}}(\Phi)=\int e^{-\Phi}\dd x\cdot\int e^{-\Phi^{\ast}}e^{-\frac{x^{2}}{2}}\dd x<\int e^{-\Phi}\dd x\cdot\int e^{-\Phi^{\ast}}\dd x\le2\pi.
\]
This shows that $\sup_{\Phi\in\class}\bs_{\alpha,\beta,\rho_{1},\rho_{2}}(\Phi)=2\pi$, but this supremum is not attained. 
\end{example}

\subsection{Existence of maximizers on the set of symmetric functions}

Using the results on symmetrization established in Section \ref{section-symmetrization}, we prove that under natural geometric assumptions (symmetry, monotonicity and log-concavity) on densities any symmetric $\bs_{\alpha,\beta,\rho_1,\rho_2}$ functional having maximizers must have, in particular, symmetric maximizers. 

\begin{theorem}\label{Steiner BS radial 2-0} Let us assume that $\rho_1,\rho_2\colon\R^n\to\R$ satisfy the following assumptions:
\begin{itemize}
\item[(R1)] $\rho_1$ is of the form  
$$
\rho=\rho(x)=\varphi(|x|)
$$
where $\varphi\colon[0,+\infty)\to[0,+\infty)$ is decreasing;
\item[(R2)] $\rho_2$ is radially symmetric and log-concave. 
\end{itemize}
If the functional $\bs_{\rho_1,\rho_2}$ admits a maximizer in some subset of $\class$, invariant under symmetrizations with respect to hyperplanes through the origin, then it admits a radially symmetric maximizer. 
\end{theorem}

\begin{proof}
It is well known (see \cite{BFR}) that for every $\Phi\in\class$ there exists a sequence of hyperplanes $H_k$, $k\in\N$, such that the sequence of functions $\phi_k$ defined recursively as follows:
$$
\Phi_1=\Phi, \quad\Phi_{k+1}=(\Phi_k)_{H_k}\quad\forall\, k\ge1,
$$
is contained in $\class$ and epi-converges to a radially symmetric function $\Phi_\infty\in\class$. In view of this fact, of Proposition \ref{Steiner BS radial} and of the continuity of $\bs_{\rho_1, \rho_2}$ we conclude the proof.
\end{proof}

\begin{theorem}\label{Steiner BS unconditional 2-0} Assume that $\rho_1$ and $\rho_2$  verify assumptions 
\begin{itemize}
\item[(U1)] $\rho_1$ is unconditional, and for every coordinate hyperplane $H$, with unit normal $e$, for every $x'\in H$, the function $\varphi\colon[0,+\infty)\to\R$ defined by
$$
\varphi(t)=\rho_1(x'+te)
$$
is decreasing; 
\item[(U2)] $\rho_2$ is unconditional and log-concave. 
\end{itemize}
If the functional $\bs_{\rho_1,\rho_2}$ admits a maximizer in some subset of $\class$, invariant under symmetrizations with respect to coordinate hyperplanes, then it admits an unconditional maximizer. 
\end{theorem}

\begin{proof}
Let $H_i$, $i=1,\dots,n$ be the coordinate hyperplanes in $\R^n$. Given $\Phi\in\class$,
the function 
$$
\Phi_u=(\dots((\Phi_{H_1})_{H_2})\dots)_{H_n}
$$
belongs to $\class$ and it is unconditional. The proof is completed using Proposition \ref{Steiner BS unconditional}.
\end{proof}

\begin{remark}
    We observe that the same symmetrization argument as in Theorem \ref{Steiner BS unconditional 2-0} was used independently in \cite{vol-prod-survey} to establish a variant of the Blaschke--Santal\'o inequality for sets and unconditional log-concave measures.
\end{remark}

\subsubsection{$1$-symmetric case}

In this subsection we prove an existence result for a special class of functions based on the compactness arguments. We note that another result of this type is proved in Section 7.

 We recall that set or a function are said to be {\it $1$-symmetric} if they possess all the symmetries of the cube. More precisely, we give the following definition. 

\begin{definition} A convex body $K$ is $1$-symmetric if for every $x=(x_1,\dots,x_n)\in K,$ we have $(\epsilon_1 x_1,...,\epsilon_n x_n)\in K$ for any choice of signs $\epsilon_i\in \{-1,1\}$, and also $(x_{\sigma(1)},...,x_{\sigma(n)})\in K$ for any permutation $\sigma.$

A function $\Phi\in\class$ is $1$-symmetric if
$$
\Phi(x_1,\dots,x_n)=\Phi(\epsilon_1 x_1,...,\epsilon_n x_n)
= \Phi(x_{\sigma(1)},...,x_{\sigma(n)})
$$
for every $x=(x_1,\dots,x_n)\in\R^n$ and every permutation $\sigma$.
\end{definition}

\begin{theorem}\label{existence homogeneous case0}
Assume $\rho_{1}$ is $s$-homogeneous and $\rho_{2}$ is $t$-homogeneous
for $s,t>-n$ such that 
$$
\frac{\alpha}{\beta}=\frac{n+s}{n+t}.
$$
 Then $\bs_{\alpha,\beta,\rho_{1},\rho_{2}}$ attains a maximum on the set of $1$-symmetric functions.
\end{theorem}

The following corollary can be deduced.

\begin{corollary}\label{vpc1}
If $V$ is a $p$-homogeneous convex function then the functional $\bs_{p,V}$
 attains a maximum on the set  of $1$-symmetric functions. 
\end{corollary}

 Note that Corollary \ref{vpc1} does not say anything about the precise form of maximizers. 
\medskip

\begin{proof}[Proof of Theorem \ref{existence homogeneous case0}]
Write $M=\sup_{\Phi\in\class_{1s}}\bs_{\alpha,\beta,\rho_{1},\rho_{2}}(\Phi)$,
and choose a sequence $\left\{ \Phi_{k}\right\} \subseteq\class_{1s}$
such that $\bs_{\alpha,\beta,\rho_{1},\rho_{2}}(\Phi_{k})\to M$. 

First observe that for every $\Phi\in\class$ and every $\lambda\in\RR$
we have $\bs_{\alpha,\beta,\rho_{1},\rho_{2}}(\Phi+\text{\ensuremath{\lambda}})=\bs_{\alpha,\beta,\rho_{1},\rho_{2}}(\Phi)$,
so we may assume without loss of generality that $\Phi_{k}(0)=0$
for all $k$. This means in particular that $\Phi_{k}\ge0$. Next,
for $\Phi\in\class$ and $\lambda>0$ we set $\left(H_{\lambda}\Phi\right)(x)=\Phi(\lambda x)$,
and recall that $\left(H_{\lambda}\Phi\right)^{\ast}=H_{1/\lambda}\Phi^{\ast}$.
It follows that 
\begin{align*}
\bs_{\alpha,\beta,\rho_{1},\rho_{2}}(H_{\lambda}\Phi) & =\left(\int_{\R^n} e^{-\alpha\Phi(\lambda x)}\rho_{1}(x)\dd x\right)^{\frac{1}{\alpha}}\left(\int_{\R^n} e^{-\beta\Phi^{\ast}(y/\lambda)}\rho_{2}(y)\dd y\right)^{\frac{1}{\beta}}\\
 & =\left(\frac{1}{\lambda^{n}}\int_{\R^n} e^{-\alpha\Phi(z)}\rho_{1}\left(\frac{z}{\lambda}\right)\dd z\right)^{\frac{1}{\alpha}}\left(\lambda^{n}\int_{\R^n} e^{-\beta\Phi^{\ast}(w)}\rho_{2}\left(\lambda w\right)\dd w\right)^{\frac{1}{\beta}}\\
 & =\left(\frac{1}{\lambda^{n+s}}\int_{\R^n} e^{-\alpha\Phi(z)}\rho_{1}\left(z\right)\dd z\right)^{\frac{1}{\alpha}}\left(\lambda^{n+t}\int_{\R^n} e^{-\beta\Phi^{\ast}(w)}\rho_{2}\left(w\right)\dd w\right)^{\frac{1}{\beta}}\\
 & =\lambda^{\frac{n+t}{\beta}-\frac{n+s}{\alpha}}\bs_{\alpha,\beta,\rho_{1},\rho_{2}}(\Phi)=\bs_{\alpha,\beta,\rho_{1},\rho_{2}}(\Phi).
\end{align*}
 
For every $k$ the set $\left[\Phi_{k}\le1\right]$ is a $1$-symmetric
convex body. 
We observe that for an arbitrary $1$-symmetric body $L$  the corresponding John ellipsoid $E$ is a ball. Indeed, if $T$ is a linear transformation satisfying $T(L)=L$, then $T(E)$ is the John ellipsoid as well. By uniqueness of $E$ one has $T(E)=E$. Hence $E$ is $1$-symmetric. This means that $E$ is a ball.
Therefore by replacing $\Phi_{k}$ with $H_{\lambda}\Phi_{k}$
for a suitable $\lambda>0$ we may assume that $\left[\Phi_{k}\le1\right]$
is in John position. In particular $B_{2}^{n}\subseteq\left[\Phi_{k}\le1\right]\subseteq\sqrt{n}B_{2}^{n}$. 

By passing to a subsequence we may assume without loss of generality
that $\left\{ \Phi_{k}\right\} $ epi-converges to a lower semi-continuous
convex function $\Phi:\RR^{n}\to[-\infty,\infty]$ (such a converging
subsequence always exists). Clearly $\min\Phi=\Phi(0)=0$
and $\Phi$ is 1-symmetric, so in order to prove that $\Phi\in\class_{1s}$
it is enough to show that $\inte\left(\dom(\Phi)\right)\ne\emptyset$
and that $\lim_{\left|x\right|\to\infty}\Phi\left(x\right)=\infty$. 

We know that $\left[\Phi_{k}\le1\right]\to\left[\Phi\le1\right]$
in the Hausdorff sense (see the proof of Lemma 5 of \cite{Colesanti-Ludwig-Mussnig}), so in particular $B_{2}^{n}\subseteq\left[\Phi\le1\right]\subseteq\sqrt{n}B_{2}^{n}$.
This immediately shows that $\inte\left(\dom(\Phi)\right)\supseteq\inte\left(B_{2}^{n}\right)\ne\emptyset$.
On the other hand for all $x\in\RR^{n}$ with $\left|x\right|\ge\sqrt{n}$
we have 
\[
1\le\Phi\left(\frac{\sqrt{n}x}{\left|x\right|}\right)=\Phi\left(\frac{\sqrt{n}}{\left|x\right|}x+\left(1-\frac{\sqrt{n}}{\left|x\right|}\right)0\right)\le\frac{\sqrt{n}}{\left|x\right|}\Phi(x)+\left(1-\frac{\sqrt{n}}{\left|x\right|}\right)\Phi(0)=\frac{\sqrt{n}}{\left|x\right|}\Phi(x),
\]
 or $\Phi(x)\ge\frac{\left|x\right|}{\sqrt{n}}.$ This shows that
$\lim_{\left|x\right|\to\infty}\Phi\left(x\right)=\infty$ and so
$\Phi\in\class$. 

Now the arguments of Proposition \ref{continuity} show that $\bs_{\alpha,\beta,\rho_{1},\rho_{2}}$
is continuous on $\class$ with respect to epi-convergence, so $\bs_{\alpha,\beta,\rho_{1},\rho_{2}}(\Phi)=M$
as claimed. 
\end{proof}

\begin{proof}[Proof of Corollary \ref{vpc1}] We have $\bs_{p,V}=\bs_{\alpha,\beta,\rho_{1},\rho_{2}}$ where $\alpha=1$, $\beta=\frac{1}{p-1}$, $\rho_{1}=1$ and $\rho_{2}=\det D^{2}V^{\ast}$. 

Since $V$ is $p$-homogeneous we know that $V^{\ast}$ is $p^{\ast}$-homogeneous, with $\frac{1}{p}+\frac{1}{p^{\ast}}=1$. Therefore $D^{2}V^{\ast}$
is $(p^{\ast}-2$)-homogeneous and $\rho_{2}$ is homogeneous of degree
\[
t=n(p^{\ast}-2)=\frac{n}{p-1}-n.
\]
Of course $\rho_{1}$ is homogeneous of degree $s=0$. Since 
\[
\frac{n+s}{n+t}=\frac{n}{n/(p-1)}=p-1=\frac{\alpha}{\beta}
\]
 the previous theorem applies and a maximizer exists in $\class_{1s}$. 
\end{proof}

\subsubsection{Radially symmetric case: existence of maximizers}

Our next result concerns the existence of maximizers of $\bs_{\alpha,\beta,\rho_{1},\rho_{2}}$, when $\rho_1$ and $\rho_2$ are homogeneous, and radially symmetric, and it is an application of our results on symmetrization, and in particular of Proposition \ref{Steiner BS radial}. In order to apply the latter result, we would need to assume that $\rho_2$ is log-concave. As the unique log-concave and homogeneous functions defined on $\R$ are (positive) constant functions, we assume that
$$
\rho_2\equiv1,
$$
that is, $\mu_2$ is the Lebesgue measure. Next, we choose
$$
\rho_1=\rho_2(x)=\varphi(|x|)
$$
with $\varphi\colon(0,+\infty)\to\R$ defined by
$$
\varphi(r)=\frac1{r^\gamma},\quad 0\le\gamma<n.
$$
In particular the measure with density $\rho_1$ is locally finite (which can be seen using polar coordinates).

\begin{theorem} 
\label{homog-rad symm}
Let $\gamma\in[0,n)$ and let $\rho\colon\R^n\setminus\{0\}\to\R$ be defined by
$
\rho(x)=\frac1{|x|^\gamma}.
$
Let $\beta>0$ and set
$
\alpha=\beta\,\frac{n-\gamma}{n}.
$
Then the functional defined by
$$
\bs_{\alpha,\beta,\rho,1}(\Phi)=\left(\int e^{-\alpha\Phi}\rho\dd x\right)^{\frac{1}{\alpha}}\left(\int e^{-\beta\Phi^{\ast}}\dd y\right)^{\frac{1}{\beta}},
$$
admits a radially symmetric maximizer.
\end{theorem}

\begin{remark} Theorem \ref{homog-rad symm} is a particular case of Theorem \ref{main-thm}. However, the proof given below is of independent interest, because it is based only on symmetrization techniques and does not use the Pr\'ekopa--Leindler theorem. 
\end{remark}

\begin{proof}[Proof of Theorem \ref{homog-rad symm}] By Proposition \ref{Steiner BS radial} (see also the proof of the following Theorem \ref{Steiner BS radial 2-0}) we may assume that there exists a maximizing sequence $\Phi_k$, $k\in\N$, such that $\Phi_k$ is radially symmetric for every $k$. This means that 
$$
\sup_{\class}\bs_{\alpha,\beta,\rho,1}=\sup_{\class_{1,s}}\bs_{\alpha,\beta,\rho,1}.
$$
The proof can be completed applying Theorem \ref{existence homogeneous case0}.
\end{proof}

\subsection{Inequalities for radially symmetric measures}
\label{radsym}

In this subsection we compute the radially symmetric maximizer of the functional appearing in Theorem \ref{homog-rad symm}. Note that the result is a particular case of Theorem \ref{main-thm}, which will be proved at the end of this section by reduction to the unconditional case.

We assume that $\beta \ge \alpha > 0$ and set $\lambda = \frac{\alpha + \beta}{\beta}$. One can easily verify  that convex potential
$$
U(y) = \frac{1}{\lambda} |x|^{\lambda} 
= \frac{\beta}{\alpha+ \beta} |x|^{\frac{\alpha + \beta}{\beta}}
$$ pushes forward the Lebesgue measure onto the measure with density
$$
\det D^2 U = \frac{\lambda-1} {|x|^{\gamma}} = \frac{\alpha}{\beta} |x|^{-\gamma},
$$
where 
$$
\gamma = n \frac{\beta-\alpha}{\beta}.
$$
We remind the reader that $U^*(y) = \frac{1}{\lambda^*} |y|^{\lambda^*}$, where $\lambda^* = \frac{\alpha + \beta}{\alpha}$.
According to Theorem \ref{homog-rad symm}, the functional
\begin{align*}
\Bigl( \int e^{-\alpha \Phi(\nabla U^*(y)) } dy\Bigr)^{\frac{1}{\alpha}}
\Bigl( \int e^{-\beta \Phi^*(y)} dy\Bigr)^{\frac{1}{\beta}}
& =
\Bigl( \int e^{-\alpha \Phi(|y|^{\frac{\beta}{\alpha}-1} y) } dy\Bigr)^{\frac{1}{\alpha}}
\Bigl( \int e^{-\beta \Phi^*(y)} dy\Bigr)^{\frac{1}{\beta}}
\\& = \Bigl( \frac{\alpha}{\beta} \Bigr)^{\frac{1}{\alpha}}
\Bigl( \int \frac{e^{-\alpha \Phi(x)}}{|x|^{\gamma}} dx\Bigr)^{\frac{1}{\alpha}}
\Bigl( \int e^{-\beta \Phi^*(y)} dy\Bigr)^{\frac{1}{\beta}}
\end{align*}
admits a radially symmetric maximizer, which is $\lambda$-homogeneous, by {Theorem \ref{p-hom-solutions}} (note that we need the assumption $\alpha \le \beta$ to ensure that the weight $|x|^{-\gamma}$ is decreasing). Thus the choice  $\Phi = U$  is optimal because the only radially symmetric $\lambda$-homogeneous function is $x \to |x|^{\lambda}$ and we get the following result.


\begin{corollary}
Let $\beta \ge \alpha \ge 0$. Then for every convex and even function $\Phi$ the following inequality holds:
$$
\Bigl( \int e^{-\alpha \Phi(|y|^{\frac{\beta}{\alpha}-1} y) } dy\Bigr)^{\frac{1}{\alpha}}
\Bigl( \int e^{-\beta \Phi^*(y)} dy\Bigr)^{\frac{1}{\beta}}
\le 
\Bigl( \int e^{-\frac{\alpha\beta}{\alpha + \beta} |y|^{1+\frac{\beta}{\alpha}}} dy\Bigr)^{\frac{1}{\alpha}+\frac{1}{\beta}}.
$$
The equivalent weighted version is:
$$
\Bigl( \int \frac{e^{-\alpha \Phi(x)}}{|x|^{\gamma}} dx\Bigr)^{\frac{1}{\alpha}}
\Bigl( \int e^{-\beta \Phi^*(y)} dy\Bigr)^{\frac{1}{\beta}}
\le \Bigl( \int \frac{e^{-\frac{\alpha\beta}{\alpha + \beta} |x|^{\frac{\alpha + \beta}{\beta}}}}{|x|^{\gamma}} dx\Bigr)^{\frac{1}{\alpha}} \Bigl( \int e^{-\frac{\alpha\beta}{\alpha + \beta} |y|^{\frac{\alpha + \beta}{\alpha}}} dy\Bigr)^{\frac{1}{\beta}}, \ \gamma = n \Bigl( 1 - \frac{\alpha}{\beta}\Bigr). 
$$
\end{corollary}
Taking, in particular, $p \ge 2$  
$$
\alpha = \frac{1}{p-1}, \ \beta=1
$$
(by homogeneity the general case can be reduced to this situation), one gets the following result
\begin{equation}
\label{3107}
\Bigl( \int e^{-\frac{1}{p-1} \Phi(|x|^{p-2} x)} dx \Bigr)^{p-1}
\Bigl( \int e^{-\Phi^*(y)} dy
\Bigr)
\le \Bigl( \int e^{-\frac{1}{p}|x|^p} dx \Bigr)^p.
\end{equation}

\begin{remark} One may ask whether (\ref{3107}) holds also for $1 < p <2$. We will see that this is not true. Indeed, Proposition \ref{BL1906} implies that if it is the case, then the probability measure $\mu = C e^{-\frac{1}{p}|x|^p} dx$ must satisfy the strong Brascamb--Lieb inequality with constant $1-\frac{1}{p}$. As we will see in Section \ref{section-BL}, this is not true for $p<2$.
\end{remark}

\begin{remark} The following inequality of Blaschke--Santal\'o type, with a radially symmetric maximizer, has been proved by Fradelizi and Meyer in \cite{FradeliziMeyer2007}. Given a decreasing function $\rho$ on $\mathbb{R}_+$ and positive even functions $f,g$ satisfying
$$
f(x) g(y) \le \rho^2(\langle x, y \rangle),
$$ for all $x,y$ such that $\langle x, y \rangle \ge 0$,
one has
\begin{equation}
    \label{fminequality}
\int f dx \int g dy \le \Bigl( \int \rho(|x|^2) dx\Bigr)^2.
\end{equation}
Though in this paper we do not analyze relations between our result and inequality of Fradelizi and Meyer, we observe that the strong Brascamp--Lieb inequality deduced from (\ref{3107}) (see Proposition \ref{BL1906}) is weaker than the infinitesimal version of (\ref{fminequality}). The latter coincides with an improvement of the Brascamp--Lieb inequality obtained by Cordero-Erausquin and Rotem in  \cite{CorRot}. The relation between  (\ref{fminequality}) and the result of  Cordero-Erausquin and Rotem was noticed in \cite{{GFSZ}} (see Theorem 4.1).
\end{remark}

 \section{The geometric approach: an equivalence between functional inequalities and inequalities about convex bodies}\label{set-version}

In this section we will see that for a given $p$-homogeneous convex function $V$, the inequality
$$ 
\int_{\R^n} e^{-\Phi(x)} dx\cdot \left(\int_{\R^n} e^{-\frac{1}{p-1}\Phi^*(\nabla V)} dx\right)^{p-1}\leq \left(\int_{\R^n} e^{-V} dx\right)^{p}
$$
for arbitrary convex functions $\Phi$, is equivalent to the following geometric inequality 
$$
|K| \cdot |\nabla V^*(K^{\circ})|^{p-1}
\le \Bigl| \Bigl\{ V \le \frac{1}{p} \Bigr\}\Bigr|^{p}.
$$

In particular, in the next subsection, we prove Proposition \ref{set-func}.

\subsection{The reduction of the functional inequality to a geometric one} For $p\geq 1$, consider $V=\|x\|^p_M/p$, where $M$ is a symmetric convex body, and $|\cdot|_M$ is the associated Minkowski functional. Then $V$ is $p$-homogeneous and all of its level sets are homothetic to $M.$ For a set $K$ in $\R^n$ we shall use the notation
$$
\nabla V(K)=\{\nabla V(x):\,x\in K\}.
$$
Note that the volume of the set $\nabla V(K)$ is given by
$$
|\nabla V(K)|=\int_K\det(\nabla^2 V(x))dx.
$$

\begin{proposition}\label{implication-p} Let $p \ge 1$. Fix a symmetric convex body $M$ and let $V=\|x\|^p_M/p$. Suppose that for any symmetric convex body $K,$ one has:
$$
|K|\cdot|\nabla V^* (K^\circ)|^{p-1}\leq \left|\left\{V\leq \frac{1}{p}\right\}\right|^p=|M|^p,
$$
with equality when $K=M.$ Then for any even strictly convex $\Phi:\R^n\rightarrow \R$ we have
$$
\int_{\R^n} e^{-\Phi(x)} dx\cdot \left(\int_{\R^n} e^{-\frac{1}{p-1}\Phi^*(\nabla V)} dx\right)^{p-1}\leq \left(\int_{\R^n} e^{-V} dx\right)^{p}.
$$
\end{proposition}
\begin{proof} We follow the scheme of Artstein-Avidan, Klartag and Milman \cite{AKM} and Keith Ball \cite{Ball1986}. Note (in view of the definition of the Legendre transform) that for any $x\in \{\Phi^*(\nabla V)\leq s\}$ and any $y\in \{\Phi(y)\leq t\}$ one has $\langle \nabla V(x),y\rangle \leq s+t$; therefore 
$$
\nabla V\left(\{\Phi^*(\nabla V)\leq s\}\right)\subset (s+t)\{\Phi\leq t\}^o.
$$
By the $\frac{1}{p-1}$-homogeneity of $\nabla V^*$ and the relation $\nabla V=(\nabla V^*)^{-1}$, the above is equivalent to
\begin{equation}\label{eq-polars}
\{\Phi^*(\nabla V)\leq s\}\subset \nabla V^*\left((s+t)\{\Phi\leq t\}^o\right)=(s+t)^{\frac{1}{p-1}}\nabla V^*\left(\{\Phi\leq t\}^o\right).
\end{equation}
Using the ``layer-cake'' representation, we write
$$
\int_{\R^n} e^{-\Phi(x)} dx\cdot \left(\int_{\R^n} e^{-\frac{1}{p-1}\Phi^*(\nabla V)} dx\right)^{p-1}=$$
\begin{equation}\label{start-p}
\int_0^{\infty} e^{-t}|\{\Phi\leq t\}| dt \cdot \left(\int_0^{\infty} e^{-s}|\{\Phi^*(\nabla V)\leq s(p-1)\}|ds\right)^{p-1}.
\end{equation}
Consider the functions
$$
f(t)=e^{-t}|\{\Phi\leq t\}|,\quad g(s)=e^{-s}|\{\Phi^*(\nabla V)\leq s(p-1)\}|\quad\mbox{and}\quad h(\tau)=e^{-\tau}|\{V\leq \tau\}|.
$$ 
Letting $K_t=\{\Phi\leq t\}$, by (\ref{eq-polars}) and the assumption of the Proposition, we get
\begin{eqnarray*}
|K_t|^{\frac{1}{p}}\cdot|\{\Phi^*(\nabla V)\leq s(p-1)\}|^{\frac{p-1}{p}}&\leq& 
(s(p-1)+t)^{\frac{n}{p}}|K_t|^{\frac{1}{p}}\cdot|\nabla V^* (K_t^o)|^{\frac{p-1}{p}}\\
&\leq& (s(p-1)+t)^{\frac{n}{p}}\left|\left\{V\leq \frac{1}{p}\right\}\right|\\
&=&\left|\left\{V\leq \frac{s(p-1)+t}{p}\right\}\right|.
\end{eqnarray*}
Therefore,
$$h\left(\frac{1}{p}t+\frac{p-1}{p}s\right)\geq f(t)^{\frac{1}{p}}g(s)^{\frac{p-1}{p}},$$
and the conclusion follows by (\ref{start-p}) combined with the Pr\'ekopa-Leindler inequality.
\end{proof}

\begin{remark}
Note that the assumption of proposition \ref{implication-p} is equivalent to the inequality
$$
|K|\cdot \left(\int_{K^o}\det D^2 V^*(x) dx\right)^{p-1}\leq \left|\left\{V\leq \frac{1}{p}\right\}\right|^p.
$$
\end{remark}

Let us finally prove that inequality (\ref{KKo}) follows from our generalized weighted functional Blaschke--Santal\'o inequality.

\begin{lemma}
\label{170724}
Let $\Phi$ be a convex $p$-homogeneous function:
$
\Phi = \|x\|^p_K.
$
Here  $K$ is a convex symmetric set.
Then 
$$
\int_{\R^n} e^{-\Phi(x)} dx\cdot \left(\int_{\R^n} e^{-\frac{1}{p-1}\Phi^*(\nabla V(x))} dx\right)^{p-1}
= c(n,p) |K|\cdot|\nabla V^* (K^\circ)|^{p-1}
$$
for some constant $c(n,p)$ depending only on $n$ and $p$.
\end{lemma}

\begin{proof}
    The proof is based on direct computations. First we apply polar coordinates:
    \begin{align*}
\int_{\R^n} e^{-\Phi} dx =
    \int_{\R^n} e^{-\|x\|^p_K} dx
    & = \int_0^{\infty}\Bigl( \int_{\mathbb{S}^{n-1}} e^{-r^p \|y\|^p_K} d\sigma(y) \Bigr) r^{n-1} dr 
    = \int_{\mathbb{S}^{n-1}} \Bigl( \int_0^{\infty} e^{- r^p \|y\|^p_K} r^{n-1} dr \Bigr) d\sigma(y) 
    \\& = \int_{\mathbb{S}^{n-1}} \frac{1}{\|y\|^n_K} \Bigl( \int_0^{\infty} e^{-s^p } s^{n-1} ds \Bigr) d\sigma(y) =
    n {\rm Vol}(K) \cdot {\int_0^{\infty} e^{- s^p } s^{n-1} ds}.
\end{align*}
By Lemma \ref{aprconj})
$$
\Phi^*(y) = \frac{p^{1-p*}}{p^*} \|y\|^{p*}_{K^{\circ}},
$$
where $p^* = \frac{p}{p-1}$.
Hence
$$
\frac{1}{p-1} \Phi^*(\nabla V(x))
=  \frac{p^{1-p^*}}{p} \|\nabla V(x)\|^{p^*}_{K^{\circ}}.
$$
Applying definition  of the Minkowski functional and homogeneity of $V$, one gets
\begin{align*}
\|\nabla V(x)\|_{K^{\circ}}
& = \inf \{ t: \nabla V(x)  \in t K^{\circ} \}
= \inf \{ t: x  \in \nabla V^*(t K^{\circ}) \}
= \inf \{ t: x  \in  t^{\frac{1}{p-1}}\nabla V^*( K^{\circ}) \}
\\& =\inf \{ s^{p-1}: x  \in  s \nabla V^*( K^{\circ}) \} = \|x\|^{p-1}_{\nabla V^*(K^{\circ})}.
\end{align*}
Thus we get
$$
\frac{1}{p-1} \Phi^*(\nabla V(x))
=
p^{-p^*} \|x\|^{p^*(p-1)}_{\nabla V^*(K^{\circ})} = p^{-p^*} \|x\|^p_{\nabla V^*(K^{\circ})}.
$$
 Applying polar coordinates again we deduce: 
$$
\int_{\R^n} e^{-\frac{1}{p-1} \Phi^*(\nabla V(x))} dx=\int_{\R^n} e^{- p^{-p^*} \|x\|^p_{\nabla V^*(K^{\circ})}} dx=n {\rm Vol}(\nabla V^*(K^{\circ}))\cdot p^{\frac{np^*}{p}} \int_0^{\infty} e^{- s^p } s^{n-1} ds.
$$
Finally
$$
\Bigl( \int e^{-\frac{1}{p-1} \Phi^*(\nabla V(x))} dx \Bigr)^{p-1}=\Bigl( n {\rm Vol}(\nabla V^*(K^{\circ})) \Bigr)^{p-1}\cdot p^n \Bigl( \int_0^{\infty} e^{- s^p } s^{n-1} ds \Bigr)^{p-1}.
$$
\end{proof}

\begin{corollary}
Let $V$ be a $p$-homogeneous strictly convex symmetric function. Inequality (\ref{main-question-intro}) holds for an arbitrary convex symmetric $\Phi$ if and only if inequality (\ref{KKo}) holds for an arbitrary symmetric convex body $K$.
\end{corollary}
\begin{proof}
Implication $(\ref{KKo}) \Longrightarrow (\ref{main-question-intro})$ was proved in Proposition \ref{implication-p}. To prove $(\ref{main-question-intro}) \Longrightarrow (\ref{KKo})$ let us take  a symmetric convex body $K$ and define $\Phi=\|x\|_K^p$. One has
    $$
\bs_{p,V}(\Phi) \le 
\bs_{p,V}(V) = \bs_{p,V}(\alpha V),
    $$
where $\alpha$ is an arbitrary positive constant. Note that in the last equality we used the invariance of $\bs_{p,V}$ with respect to homotheties and the homogeneity of $V$. Applying Lemma \ref{170724} one gets
    $$
|K| \cdot |\nabla V^*(K^{\circ})|^{p-1}
\le |K_{\alpha}| \cdot |\nabla V^*(K^{\circ}_{\alpha})|^{p-1}
    $$
    where $K_{\alpha}=\bigl\{ V \le \frac{1}{\alpha} \bigr\}$.
    It remains to prove that 
    $
\nabla V^*(K^{\circ}_{p})
= K_p.
    $
    Indeed, one has $V(x) = \|x\|^p_{K_1}$, hence
    $V^*(y) = \frac{p-1}{p^{p^*}} \|y\|^{p^*}_{K_1^{\circ}}$ (see Lemma \ref{aprconj}). From the relation $K_p =\Big\{ x: \|x\|_{K_1} \le \frac{1}{p^{\frac{1}{p}}} \Big\}$ we infer 
    $$
K_p^{\circ} = \Big\{ y : \|y\|_{K_1^\circ} \le p^{\frac{1}{p}}\Big\} = \Big\{ V^* \le \frac{1}{p^*}\Big\}.
    $$

    Let us show that $\nabla V^*(K^{\circ}_p) \subset K_p$.
    Indeed, by $p$-homogeneity of $V$ one has
    $$V(\nabla V^*(y))= \langle y, \nabla V^*(y)\rangle - V^*(y) = (p^*-1)V^*(y) =\frac{1}{p-1}V^*(y) = \frac{1}{p^{p^*}} \|y\|^{p^*}_{K_1^{\circ}}.$$ Then for every $y \in K_p^\circ$ one has $\|y\|_{K^{\circ}_1} \le p^{\frac{1}{p}}$ and
    $
V(\nabla V^*(y)) \le p^{\frac{p^*}{p} -p^*} = \frac{1}{p}.$ This proves that  $\nabla V^*(K^{\circ}_p) \subset K_p$.

    The reverse inclusion $K_p \subset \nabla V^*(K^{\circ}_p) $ is equivalent to $\nabla V(K_p) \subset K_p^\circ$. This can be proved by the same argument, where we replace $V$ by $V^*$, $p$ by $p^*$ and vice versa.
\end{proof}

\subsection{The case of rotation-invariant measures revisited}

In this subsection we show that Theorem \ref{main-thm} in the case of rotationally invariant measures follows from the classical Blaschke--Santal{\'o} inequality.

Suppose $V=\frac{|x|^p}{p}$ and $V^*=\frac{|x|^{p^*}}{p^*},$ with $p$ and $p^*$ conjugate to each other. Then 
$$
\nabla^2 V^*(x)=|x|^{p^*-2}{\rm Id}+(p^*-2)|x|^{p^*-4}x\otimes x,
$$ 
where $x \otimes x$ is a symmetric matrix defined by the relation
$(x \otimes x)(a,b) = \langle a,x \rangle \langle b,x \rangle$,
and thus $\det(\nabla^2 V^*(x))=(p^*-1)|x|^{n(p^*-2)}$. Therefore, the condition of Proposition \ref{implication-p} in the case $V(x)=\frac{|x|^p}{p}$ is: for any symmetric convex $K$ and for $p,p^*\geq 1,$ with $\frac{1}{p}+\frac{1}{p^*}=1,$
$$(p^*-1)^{p-1}|K|\cdot \left(\int_{K^o} |x|^{n(p^*-2)}dx\right)^{p-1}\leq |B^n_2|^{p},$$
which becomes, in view of the fact that $p^*=\frac{p}{p-1}:$
$$(p^*-1)^{\frac{1}{p^*-1}}|K|\cdot \left(\int_{K^o} |x|^{n(p^*-2)}dx\right)^{\frac{1}{p^*-1}}\leq |B^n_2|^{\frac{p^*}{p^*-1}}.$$
Therefore, the rotationally invariant case  of Theorem \ref{main-thm} follows immediately from Proposition \ref{implication-p} and the following result.

\begin{proposition}\label{rot-inv-shad-syst}
For any symmetric convex $K$ and any $p^*\in(1,2],$ we have
$$(p^*-1)^{\frac{1}{p^*-1}}|K|\cdot \left(\int_{K^o} |x|^{n(p^*-2)}dx\right)^{\frac{1}{p^*-1}}\leq |B^n_2|^{\frac{p^*}{p^*-1}}.$$
\end{proposition}

\begin{remark} Note that the condition $p\geq 2$ in (1) of Theorem \ref{main-thm} corresponds to the assumption $p^*\in [1,2]$, since $p$ and $p^*$ are conjugate.
\end{remark}

We show that  Proposition \ref{rot-inv-shad-syst} follows immediately from the classical Blaschke--Santal\'o inequality. 

\begin{lemma}\label{expression}
For a convex body $K$ we have
$$\int_{K^o} |x|^{n(p^*-2)} dx=\frac{1}{(p^*-1)n}\int_{\sfe} h^{(1-p^*)n}_K(\theta) d\sigma.$$ 
\end{lemma}
\begin{proof} Using the polar coordinates, we write
$$
\int_{K^o} |x|^{n(p^*-2)} dx=\int_{\sfe}\int_0^{\rho_{K^o}} t^{n-1+n(p^*-2)}dt d\sigma,
$$
and the equality follows from the fact that $h^{-1}_K=\rho_{K^o}.$
\end{proof}

\begin{proof} [Proof of Proposition \ref{rot-inv-shad-syst}] Using Lemma \ref{expression} combined with H\"older's inequality, we write
\begin{eqnarray*}
(p^*-1)^{\frac{1}{p^*-1}}|K|\cdot \left(\int_{K^o} |x|^{n(p^*-2)}dx\right)^{\frac{1}{p^*-1}}&=&(p^*-1)^{\frac{1}{p^*-1}}|K|\cdot \left(\frac{1}{(p^*-1)n}\int_{\sfe} h^{(1-p^*)n}_K(\theta) d\sigma(\theta)\right)^{\frac{1}{p^*-1}}\\
&\leq&C(n,p^*) |K|\cdot\int_{\sfe} h^{-n}_K(\theta) d\sigma(\theta)\\
&=&C'(n,p^*)|B^n_2|^2,
\end{eqnarray*}
where in the last step we used the Blaschke-Santal\'o inequality and identity $|K^{o}|=\frac{1}{n} \int_{\mathbb{S}^{n-1}} h_K^{-n}(\theta) d \sigma(\theta)$. Here $C(n,p^*)$ and $C'(n,p^*)$ are appropriate constants depending on $n$ and $p$, such that equality is attained in all the inequalities above when $K$ is $B^n_2$. This completes the proof of the Proposition.
\end{proof}

\section{The mass transport approach to the Blaschke-Santal\'o-type inequalities}

In this section we consider the classical case: $\mu_1$ and $\mu_2$ both coincide with the Lebesgue measure.
Our aim is to get a mass transportation proof of the classical Blaschke--Santal{\'o} inequality.
First we prove that functional $\bs_{1,1,1,1}$ is bounded from above and  that maximizers exist. Then we prove uniqueness (up to linear transformation) of the maximizers. Our proof is based on the analysis of the Monge--Amp\`ere equation satisfied by the maximizers. In addition, we establish a general relation between the Blaschke--Santal{\'o} and the Brascamp--Lieb (variance) inequality.

\subsection{Existence of maximizers for the classical Blaschke--Santal\'o functional}

For the sake of brevity we will write $\mathcal{BS}$ for the classical Blaschke--Santal\'o functional $\bs_{1,1,1,1}$. 

\begin{theorem}\label{existence theorem}
There exists $\Phi\in\class$ such that 
$$
\mathcal{BS}(\Phi) = M:= \sup_{\Psi\in\class}\bs(\Psi).
$$
Making appropriate linear change of variables and normalization, one can assume that the measure $\mu$ with density $e^{-\Phi}$ (with respect to the Lebesgue measure) is an isotropic probability measure and $\Phi^*$ satisfies the inequality
\begin{equation}\label{vstarbound}
\Phi^*(y) \ge \frac{1}{2n} \sum_{i=1}^{n} |y_i| - c(n)
\end{equation}
for some constant $c(n)$ depending only on dimension.
\end{theorem}

We recall that a probability measure $\mu$ on $\R^n$ is said to be isotropic if 
$$
\int_{\R^n} \langle x, \theta \rangle^2 d \mu=1
$$
for every unit vector $\theta$. 
Also we recall that in this subsection $\bs$ is the classical 
Blaschke--Santal\'o functional:
$$
\bs(\Phi) = \int_{\R^n} e^{-\Phi} dx \int_{\R^n} e^{-\Phi^*} dy.
$$

We start with a remark. Let $\Phi\in\class$. Without changing the value of $\bs(\Phi)$, we can assume that the measure $\mu$ with density $e^{-\Phi}$ is an isotropic probability measure. This is possible because, by Remark \ref{invariance}, $\mathcal{BS}$ is invariant with respect to transformations of the form $\Phi(x) \to \Phi(Ax) + b$, where $b\in\R$ and $A \colon \mathbb{R}^n \to \mathbb{R}^n$ is an invertible linear transformation. So, we first reduce to a probability density adding a suitable constant, and then, taking a linear image, we make the measure isotropic.

\begin{lemma}\label{vbound} Let $\Phi\in\fconvx$ be coercive and assume that the measure $\mu$ with density $e^{-\Phi}$ is an isotropic probability measure. There exists a constant $ c = c(n) >0$ depending on $n$ such that
$$
\Phi(x)\le c(n)
$$
for every $x\in B = \{ x: |x| \le \frac{1}{2} \}$.
\end{lemma}
\begin{proof} 
First we observe that $\Phi < \infty$ on $B$. Indeed, if it is not the case, take $x \in B \cap \partial\{\Phi < \infty\}$ and a normal unit vector $\theta$  of a support hyperplane  $H$ for $\{x: \Phi(x) < \infty\}$ at $x$. Since $\Phi$ is even, we conclude that $\{\Phi < \infty\} \subset \{y:| \langle y, \theta\rangle| \le  |\langle x,\theta  \rangle|\}$, but this contradicts to isotropicity, because $|\langle x, \theta \rangle| \le \frac{1}{2}$, hence the linear function $\langle y, \theta \rangle$ is bounded by $1/2$ on $B$.

As $\Phi$ is convex and finite  on $B$, by continuity of convex functions there exists $x_0\in\partial B$ such that 
$$
m = \min_{x \in B} e^{-\Phi(x)}
= e^{-\max_{x \in B} \Phi(x)}=e^{-\Phi(x_0)}.
$$ 
One has $x_0 =  \frac{\theta}{2}$, for some $\theta\in\sn$. Let $L = \{ x\in\R^n: |\langle x, \theta \rangle| \le  \frac{1}{2} \}$. We observe that, $e^{-\Phi} \le m$ on $ \mathbb{R}^n \setminus L$. 
Indeed, if $\Phi(x_0)=\Phi(0)$, then $\Phi(x) \ge \Phi(x_0)$ for all $x \in \mathbb{R}^n$, because $0$ is the global minimum point of $\Phi$. If  $\Phi(0) < \Phi(x_0)$, then $\Phi(tx_0)> \Phi(x_0)$ and $x_0$ belongs to  the boundary of $D= \{ y: \Phi(y) \le \Phi(x_0)\}$. Take a unit supporting vector $\tau$ of $D$ at $x_0$:
$$
D \subset \{x : \langle \tau, x -x_0 \rangle \le 0 \}.
$$ 
Since $B \subset D$ by definition of $x_0$, $\tau$ is also a supporting vector to $B$. But this means that $\tau = \theta$. Hence for every $x \in D$ one has
$\langle \tau, x -x_0 \rangle \le 0$, equivalently 
$
\langle x, \theta \rangle \le \big\langle \theta, \frac{\theta}{2}\big\rangle =\frac{1}{2} . 
$
Thus $\langle x, \theta\rangle$ is bounded by $1/2$ on $D$ and by the symmetry assumption  $|\langle x, \theta\rangle| \le \frac{1}{2}$ on $D$. This proves that $D \subset L$, hence $\Phi\ge \Phi_0$ and $e^{-\Phi} \le m$ on $\mathbb{R}^n \setminus L$. 

Since
$$
1 = \int_L  \langle x, \theta \rangle^2 d \mu
+ \int_{ \mathbb{R}^n \setminus L}  \langle x, \theta \rangle^2 d \mu
\le \frac{1}{4} + \int_{ \mathbb{R}^n \setminus L} \langle x, \theta \rangle^2 d \mu,
$$
one has
\begin{align*}
    \frac{3}{4} \le \int_{ \mathbb{R}^n \setminus L} \langle x, \theta \rangle^2 d \mu
    \le \int_{ \mathbb{R}^n \setminus L} |x|^2 d \mu
    & \le \Bigl( \int_{ \mathbb{R}^n \setminus L} |x|^{n+5} d \mu  \Bigr)^{\frac{1}{2}} \Bigl( \int_{ \mathbb{R}^n \setminus L} \frac{1}{|x|^{n+1}} d \mu \Bigr)^{\frac{1}{2}} \\&
   \le  \Bigl( \int_{ \mathbb{R}^n} |x|^{n+5} d \mu  \Bigr)^{\frac{1}{2}}  \Bigl( \int_{ \mathbb{R}^n \setminus L} \frac{dx}{|x|^{n+1}} \Bigr)^{\frac{1}{2}} \sqrt{m}.
\end{align*}
Using the well-known moment equivalence result for log-concave measures (see Theorem 3.5.11 in \cite{AGM}) and the isotropicity of $\mu$, one gets that
$  \int_{ \mathbb{R}^n} |x|^{n+5} d \mu  $ is bounded by a number depending on $n$. Clearly,  since $B \subset L$, one has $\int_{ \mathbb{R}^n \setminus L} \frac{dx}{|x|^{n+1}} \le \int_{ \mathbb{R}^n \setminus B} \frac{dx}{|x|^{n+1}} \le  C(n)<\infty$. Thus we get
$m \ge c(n) >0$. This completes the proof.
\end{proof}

\begin{corollary} Let $\Phi\in\fconvx$ be coercive, and assume that the measure with density $e^{-\Phi}$ is an isotropic probability measure. Then $\Phi^*$ satisfies the inequality
\begin{equation*}
\Phi^*(y) \ge \frac{1}{2n} \sum_{i=1}^{n} |y_i| - c(n)
\end{equation*}
for every $y=(y_1,\dots,y_n)\in\R^n$.
\end{corollary}
\begin{proof} Let $\{e_i\colon i=1,\ldots,n\}$ be the standard orthonormal basis in $\R^n$. Apply inequality
$$
\Phi^*(y) \ge \langle x, y \rangle - \Phi(x)
$$
to $x = \pm \frac{1}{2} e_i$. The previous lemma implies $\Phi^*(y) \ge \frac{1}{2} |y_i| - c(n)$. The result immediately follows.
\end{proof}

\begin{proof}[Proof of Theorem \ref{existence theorem}]
Let $\Phi_k$, $k\in\N$, be a sequence of coercive functions in $\fconvx$ such that
$$
\lim_{k\to+\infty}\mathcal{BS}(\Phi_k)=\sup\{\bs(U)\colon U\in\fconvf\}.
$$ 
As already remarked, we may assume that for every $k\in\N$, the measure $\mu_k$ with density $e^{-\Phi_k}$ is a probability measure and it is isotropic. 

By Chebyshev inequality and isotropicity, we have
$$
\int_{|x| > R} e^{-\Phi_k}dx \le \frac{1}{R^2}\int_{|x| > R} |x|^2 e^{-\Phi_k}dx \le \frac{n}{R^2}
$$
for every $R>0$. From the last inequality we infer that the sequence of measures $\mu_k$ is tight, and then Prokhorov theorem can be applied (see, for instance, \cite{Bogachev}, Theorem 8.6.2). Therefore $\mu_k$ admits a subsequence which is weakly convergent to a probability measure $\mu$.   

As $\mu_k$ is log-concave and isotropic for every $k$, it is easy to see that $\mu$ is log-concave and isotropic, as well. By a well-known theorem of Borell (see \cite{Borell}, Theorem 3.2), the log-concavity of $\mu$ implies that it is absolutely continuous with respect to the Lebesgue measure, and its density is of the form $e^{-\Phi}$, where $\Phi\in\fconvx$. As $\mu$ is a probability measure, $\Phi$ is coercive; moreover, by Lemma \ref{vbound}, $\Phi(x)<\infty$ for every $x$ such that $|x|\le\frac12$. 

The weak convergence of $\mu_k$ to $\mu$ implies that
$$
\lim \Phi_k(x)=\Phi(x),\quad\forall\, x\in\R^n\setminus\partial\dom(\Phi).
$$
By \cite[Theorem 7.17]{Rockafellar-Wets}, $\Phi_k$ epi-converges to $\Phi$. By Proposition \ref{continuity},
$$
\bs(\Phi)=\lim_{k\to\infty}\bs(\Phi_k)=M.
$$
Note that $\Phi^*$ verifies \eqref{vstarbound}, so that, in particular, $M<+\infty$. 
\end{proof}

\begin{remark}
     Theorem \ref{existence theorem} is used in the proof of the classical Blaschke--Santal\'o  inequality by transportation method without  symmetrization arguments (see Remark \ref{classicalBS}).    
\end{remark}

{ \subsection{The Euler-Lagrange equation of $\bs$}
In this subsection we derive the Euler-Lagrange equation of the Blaschke-Santal\'o functional. 
We realize that this is an equation of the Monge--Amp\`ere type. In the following subsection we prove a kind of more precise statement: a monotonicity property of our functional, which also leads to this equation.

The following lemma is known to experts; we include it for the reader's convenience.

\begin{lemma}\label{v*inf}
Let $V\in C^2(\R^n)$ be such that $D^2 V(x)$ is positive definite for every $x$, and let $f\in C^2(\R^n)$ be compactly supported. For $\varepsilon >0$, let
$$
V_{\varepsilon} = V + \varepsilon f.
$$
Then, $V^{*}_{\varepsilon}$ can be expanded in the following way:
$$
V^{*}_{\varepsilon} = V^* - \varepsilon f(\nabla V^*) + \frac{\varepsilon^2}{2} \langle D^2 V^* \nabla f(\nabla V^*), \nabla f(\nabla V^*) \rangle + o(\varepsilon^2). 
$$
The dependence of the term $o(\varepsilon^2)$ on $x$ in this expansion is uniform on $\nabla V({\rm supp}(f))$.
\end{lemma}

\begin{proof}
    Expand $V^*_{\varepsilon}$ in the variable $\varepsilon$ to get:
    $$
    V^*_{\varepsilon}(y) = V^*(y) + \varepsilon a(y) + \frac{\varepsilon^2}{2} b(y) + o(\varepsilon^2),
    $$  where $a,b$ are functions on $\mathbb{R}^n$. 
    
    We observe that for sufficiently small $\varepsilon$ function    
    $V_\varepsilon$ remains convex with a positive hessian at every point. 
    Thus we can use the relation $ V^*_{\varepsilon} (\nabla V_{\varepsilon}(x)) = \langle x, \nabla V_{\varepsilon}(x) \rangle - V_{\varepsilon}(x). $
    In this way we obtain
    \begin{align}
    \label{vstarexp}
     V^*(\nabla V_{\varepsilon}(x)) + \varepsilon a(\nabla V_{\varepsilon}(x)) + \frac{\varepsilon^2}{2} b(\nabla V_{\varepsilon}(x)) + o(\varepsilon^2)
&   \nonumber  = \langle x, \nabla V_{\varepsilon}(x) \rangle - V_{\varepsilon}(x) \\& = V^*(\nabla V(x))
    + \varepsilon ( \langle x, \nabla f(x) \rangle - f(x) ).
    \end{align}
    Consider the following Taylor series expansions as functions on $\mathbb{R}^n$:
$$
V^*(\nabla V_{\varepsilon}(x)) = V^*(\nabla V(x)) + \varepsilon \langle x, \nabla f(x) \rangle + \frac{\varepsilon^2}{2} \langle D^2 V^*(\nabla V(x)) \nabla f(x), \nabla f(x) \rangle  + o(\varepsilon^2).
$$
Plugging this expression into (\ref{vstarexp}) one gets
$$
\varepsilon \langle x, \nabla f(x) \rangle +  \frac{\varepsilon^2}{2} \langle D^2 V^*(\nabla V(x)) \nabla f(x), \nabla f(x) \rangle  + \varepsilon a(\nabla V_{\varepsilon}(x)) + \frac{\varepsilon^2}{2} b(\nabla V_{\varepsilon}(x)) +o(\varepsilon^2) = \varepsilon ( \langle x, \nabla f(x) \rangle - f(x) ).
$$
Equivalently
\begin{equation}
    \label{afveps}
 a(\nabla V_{\varepsilon}(x)) + f(x)  + \frac{\varepsilon}{2} \Bigl(\langle D^2 V^*(\nabla V(x)) \nabla f(x), \nabla f(x) \rangle  +   b(\nabla V_{\varepsilon}(x)) \Bigr) = o(\varepsilon).
\end{equation}
From this we immediately get $a(\nabla V(x))= - f(x)$. Thus
$a(y)= - f(\nabla V^*(y))$ and $$a(\nabla V_{\varepsilon}(x)) = - f(\nabla V^*(\nabla V_{\varepsilon}(x))) = - f(x)  - \varepsilon \langle D^2 V^*(\nabla V(x)) \nabla f(x), \nabla f(x) \rangle + o(\varepsilon).$$
Plugging this into (\ref{afveps}) we immediately get
$b(\nabla V(x)) =  \frac{1}{2} \langle D^2 V^*(\nabla V(x)) \nabla f(x), \nabla f(x) \rangle$. The proof is complete.
\end{proof}

\begin{proposition}\label{Euler-Lagrange} Let $\mu_1$ and $\mu_2$ be admissible measures. Let $\Phi$ be a maximizer of $\bs_{\alpha,\beta,\rho_1,\rho_2}$ in $\class$, and assume that 
$$
\Phi\in C^2(\R^n)\quad\mbox{and}\quad D^2\Phi(x)>0\quad\forall\, x\in\R^n.
$$
Then 
$$
\frac{e^{-\alpha \Phi}\rho_1}{\int_{\R^n} e^{-\alpha\Phi}\rho_1dx}=\frac{e^{-\beta\Phi^*(\nabla\Phi)}\det(D^2\Phi)\rho_2(\nabla\Phi)}{\int_{\R^n} e^{-\beta\Phi^*}\rho_2dy}.
$$
\end{proposition}

\begin{proof}
Let $\tau$  be a smooth compactly supported function. For $\varepsilon>0$ sufficiently small in absolute value, the function
$$
\Phi_\varepsilon=\Phi+\varepsilon\tau
$$
belongs to $\class$ (here we are using the assumption that $D^2\Phi>0$ and that $\tau$ has compact support). Hence the function
$$
\varepsilon\,\rightarrow\,\bs_{\alpha,\beta,\rho_1,\rho_2}(\Phi_\varepsilon)
$$
has a maximum for $\varepsilon=0$. Therefore
$$
\left.\frac{d}{d\varepsilon}\bs_{\alpha,\beta,\rho_1,\rho_2}(\Phi_\varepsilon)\right|_{\varepsilon=0}=0.
$$
On the other hand
\begin{eqnarray*}
\left.\frac{d}{d\varepsilon}\bs_{\alpha,\beta,\rho_1,\rho_2}(\Phi_\varepsilon)\right|_{\varepsilon=0}&=&
\left.\frac{d}{d\varepsilon} \Bigl(\int_{\R^n}e^{-\alpha\Phi_\varepsilon}\rho_1dx\Bigr)^{\frac{1}{\alpha}}\right|_{\varepsilon=0}
\cdot \Bigl(\int_{\R^n}e^{-\beta\Phi^*}\rho_2dx\Bigr)^{\frac{1}{\beta}}+\\
&& \Bigl(\int_{\R^n}e^{-\alpha\Phi}\rho_1dx \Bigr)^{\frac{1}{\alpha}}\cdot\left.\frac{d}{d\varepsilon}\Bigl(\int_{\R^n}e^{-\beta\Phi^*_\varepsilon}\rho_2dx\Bigr)^{\frac{1}{\beta}}\right|_{\varepsilon=0}\\
&=&-\int_{\R^n}\tau e^{-\alpha\Phi}\rho_1dx  \Bigl(\int_{\R^n}e^{-\alpha\Phi_\varepsilon}\rho_1dx\Bigr)^{\frac{1}{\alpha}-1}\Bigl(\int_{\R^n}e^{-\beta\Phi^*}\rho_2dx\Bigr)^{\frac{1}{\beta}}+\\
&&\int_{\R^n}\tau(\nabla\Phi^*)e^{-\beta\Phi^*}\rho_2dy
\Bigl(\int_{\R^n}e^{-\alpha\Phi_\varepsilon}\rho_1dx\Bigr)^{\frac{1}{\alpha}}\Bigl(\int_{\R^n}e^{-\beta\Phi^*}\rho_2dx\Bigr)^{\frac{1}{\beta}-1}
\end{eqnarray*}
where we have used Lemma \ref{v*inf}. By the change of variable $\nabla\Phi^*(y)=x$ we get
$$
\int_{\R^n}\tau(\nabla\Phi^*)e^{-\beta\Phi^*}\rho_2dy=\int_{\R^n}\tau e^{-\beta\Phi^*(\nabla\Phi)}\det(D^2\Phi)\rho_2(\nabla\Phi)dx
$$
We deduce that
$$
\int_{\R^n}\tau\left[I_2 e^{-\alpha\Phi}\rho_1-I_1 e^{-\beta\Phi^*(\nabla\Phi)}\det(D^2\Phi)\rho_2(\nabla\Phi)
\right]dx=0,
$$
where
$$
I_1=\int_{\mathbb{R}^n} e^{-\alpha\Phi}\rho_1dx,\quad
I_2=\int_{\mathbb{R}^n} e^{-\beta\Phi^*}\rho_2dx.
$$
As $\tau$ is arbitrary, the conclusion follows.
\end{proof}}

\subsection{A comparison result for Blaschke--Santal\'o functional via optimal transportation} 

In this section we present a comparison result for the Blaschke-Santal\'o functional, involving optimal transportation.
We will show that the functional
$$
 \mathcal{BS}_{\alpha, \beta, \rho_1, \rho_2}(\Phi) = \Bigl( \int e^{-\alpha \Phi} \rho_1 dx \Bigr)^{\frac{1}{\alpha}} 
 \Bigl( \int e^{-\beta \Phi^*} \rho_2 dy \Bigr)^{\frac{1}{\beta}}
 $$ 
admits a remarkable  monotonicity property related to optimal transportation. This property was already mentioned in \cite{KolesnikovWerner} for $\alpha=\beta=1$. The idea of the proof is essentially the same. However, the statement about maximum points of the Blaschke--Santal\'o functional that we prove here is more precise even for values $\alpha=\beta=1$.  
 
Let $\mu_1,\mu_2$ be non-negative Borel measures on $\R^n$, absolutely continuous with respect to the Lebesgue measure, with strictly positive densities $\rho_1$ and $\rho_2$, respectively. Let $\Phi\colon\R^n\to\R\cup\{\infty\}$ be a convex function; assume that
\begin{equation}
\label{omegaomegastar}
\Omega:=\inte(\dom(\Phi))\quad\mbox{and}\quad\Omega^*:=\inte(\dom(\Phi^*))
\end{equation}
are non-empty, and 
$$
0<\int_{\R^n}e^{-\alpha\Phi}d\mu_1,\int_{\R^n}e^{-\beta \Phi^*}d\mu_2 < \infty.
$$
Let $\alpha>0, \beta>0$ and $\nabla U \colon\Omega\to\Omega^*$ be the optimal transportation of the measure $\mu$ with density 
$$
\frac1{\int_{\R^n}e^{-\alpha \Phi} d\mu_1}\,e^{-\alpha \Phi} \rho_1
$$
onto the measure $\nu$ with density
$$
\frac1{\int_{\R^n}e^{-\beta\Phi^*} d\mu_2}\,e^{-\beta\Phi^* }\rho_2.
$$ 
Let us assume that
\begin{enumerate}
\item $U(0)=\Phi(0)$,
\item $\Phi$ and $U$ are lower semi-continuous,
\item $\Phi=+\infty$ on $\{ U = + \infty\}$.
\end{enumerate}

\begin{remark}
    It is easy to verify that there exists a unique function $U$ satisfying assumptions 1)-3). Note that uniqueness for optimal transportation  guarantees only that $T = \nabla U$ is uniquely determined $\mu$-a.e. In particular, potentials $U_i$ giving the same mapping $T$ can be different outside of support $\mu$. But in our case the support of $\mu$ is convex and $U$ is supposed to take infinite values outside of it. This implies the uniqueness of $U$.
\end{remark}

We set
$$
\nabla U(\Omega) = \left\{\nabla U(x)\colon x \in \Omega,\,\partial U(x)  \ \mbox{contains a unique element}\right\}.
$$ 

\begin{proposition}\label{key comparison} In the previous assumptions
\begin{equation}\label{comparison}
\bs_{\alpha,\beta,\rho_1,\rho_2}(\Phi)\le\bs_{\alpha,\beta,\rho_1,\rho_2}(U).
\end{equation}
Moreover, equality holds if and only if
$$
\Phi=U\quad\mbox{in $\Omega$,}
$$
and $\nabla U(\Omega)$ 
coincides with $\{U^*<\infty\}$ up to a set of the Lebesgue measure zero.
\end{proposition}

\begin{proof} By definition (\ref{omegaomegastar}) of $\Omega$ and $\Omega^*$ 
$$
\int_{\R^n}e^{-\alpha \Phi} \rho_1 dx=\int_{\Omega}e^{-\alpha\Phi} \rho_1 dx,\quad
\int_{\R^n}e^{-\Phi^*} \rho_2 dy=\int_{\Omega^*}e^{-\beta\Phi^*} \rho_2 dy.
$$
By convexity, $\Phi$ and $\Phi^*$ are continuous on $\Omega$ and $\Omega^*$, respectively. By the change of variables formula (\ref{cvf})
$$
\frac{e^{-\alpha \Phi(x)} \rho_1(x)}{\int_{\Omega} e^{-\alpha\Phi} \rho_1 dx} 
= \frac{e^{-\beta\Phi^*(\nabla U(x))} \rho_2(\nabla U)}{\int_{\Omega^*} e^{-\beta\Phi^*} \rho_2 dy} \det D^2_a U(x)
$$
almost everywhere in $\Omega$.

Take the power $\frac{\alpha}{\alpha + \beta}$ of both sides of the last equality:
$$
\Bigl( \frac{\int_{\Omega^*} e^{-\beta\Phi^*} \rho_2 dy}{\int_{\Omega} e^{-\alpha\Phi}\rho_1 dx} \Bigr)^{\frac{\alpha}{\alpha+\beta}}
e^{-\frac{\alpha^2}{\alpha+\beta}\Phi(x)}
\rho_1^{\frac{\alpha}{\alpha+\beta}}(x)
= e^{-\frac{\alpha \beta}{\alpha + \beta}\Phi^*(\nabla U(x))} \rho^{\frac{\alpha }{\alpha + \beta}}_2(\nabla U) \bigl(\det D^2_a U(x)\bigr)^{\frac{\alpha}{\alpha+\beta}}.
$$
Multiply this identity by $ e^{-\frac{\alpha \beta}{\alpha+\beta}{\Phi}}\rho_1^{\frac{\beta}{\alpha+\beta}}$:
$$
\Bigl( \frac{\int_{\Omega^*} e^{-\beta\Phi^*} \rho_2 dy}{\int_{\Omega} e^{-\alpha\Phi} \rho_1 dx} \Bigr)^{\frac{\alpha}{\alpha+\beta}}
e^{-{\alpha}\Phi(x)}
\rho_1(x)
= e^{-\frac{\alpha \beta}{\alpha + \beta}\bigl[\Phi(x) + \Phi^*(\nabla U(x))\bigr]} \rho_1^{\frac{\beta}{\alpha+\beta}} \rho^{\frac{\alpha }{\alpha + \beta}}_2(\nabla U) \bigl(\det D^2_a U(x)\bigr)^{\frac{\alpha}{\alpha+\beta}}.
$$
Integrating over $\Omega$, we obtain
\begin{eqnarray*}
&&\Bigl( {\int_{\Omega}e^{-\alpha \Phi}\rho_1 dx}\Bigr)^{\frac{\beta}{\alpha+\beta}}
\Bigl({\int_{\Omega^*}e^{-\beta \Phi^*}\rho_2 dy} \Bigr)^{\frac{\alpha}{\alpha+\beta}}\\
&&=\int_{\Omega}e^{-\frac{\alpha \beta}{\alpha + \beta}\bigl[\Phi(x) + \Phi^*(\nabla U(x))\bigr]} \rho_1^{\frac{\beta}{\alpha+\beta}} \rho^{\frac{\alpha }{\alpha + \beta}}_2(\nabla U) \bigl(\det D^2_a U(x)\bigr)^{\frac{\alpha}{\alpha+\beta}} dx.
\end{eqnarray*}
As
$$
\Phi^*(\nabla U(x)) + \Phi(x) \geq \langle x, \nabla U(x) \rangle,
$$
while
$$
U^*(\nabla U(x)) + U(x) = \langle x, \nabla U(x) \rangle,
$$ 
for every $x\in\Omega$, we get
\begin{align*}
\Bigl( {\int_{\Omega}e^{-\alpha \Phi}\rho_1 dx}\Bigr)^{\frac{\beta}{\alpha+\beta}}
\Bigl({\int_{\Omega^*}e^{-\beta \Phi^*}\rho_2 dy} \Bigr)^{\frac{\alpha}{\alpha+\beta}} & \le \int_{\Omega}e^{-\frac{\alpha \beta}{\alpha + \beta}\bigl[U(x) + U^*(\nabla U(x))\bigr]} \rho_1^{\frac{\beta}{\alpha+\beta}} \rho^{\frac{\alpha }{\alpha + \beta}}_2(\nabla U) \bigl(\det D^2_a U(x)\bigr)^{\frac{\alpha}{\alpha+\beta}} dx
 \\& \le \Bigl( \int_{\Omega} e^{-\alpha U} \rho_1 dx\Bigl)^{\frac{\beta}{\alpha+ \beta}}  \Bigl(\int_{\Omega} e^{-\beta U^*(\nabla U)} \rho_2(\nabla U) \det D^2_a U dx \Bigr)^{\frac{\alpha}{\alpha+ \beta}} 
 \\& = \Bigl( \int_{\Omega} e^{-\alpha U} \rho_1 dx \Bigr)^{\frac{\beta}{\alpha+ \beta}} 
\Bigl(\int_{\nabla U(\Omega)} e^{-\beta U^*} \rho_2 dy \Bigr)^{\frac{\alpha}{\alpha+ \beta}} 
\\& \le \Bigl( \int_{\mathbb{R}^n} e^{-\alpha U} \rho_1 dx \Bigr)^{\frac{\beta}{\alpha+ \beta}} 
\Bigl(\int_{\mathbb{R}^n} e^{-\beta U^*} \rho_2 dy \Bigr)^{\frac{\alpha}{\alpha+ \beta}},
\end{align*}
where we have used H\"older inequality in the second step, and the change of variable formula in the third one. This proves \eqref{comparison}.

Assume now that
$$
\bs_{\alpha,\beta,\rho_1,\rho_2}(\Phi) = \bs_{\alpha,\beta,\rho_1,\rho_2}(U).
$$ 
This is possible only if 
\begin{equation}\label{equality}
\Phi^*(\nabla U(x)) + \Phi(x) =\langle x, \nabla U(x) \rangle,
\end{equation}
on $\Omega$. On the other hand, if $x\in\Omega$ and $\Phi$ is differentiable at $x$, then
\begin{equation}\label{equality2bis}
\Phi^*(\nabla \Phi(x)) + \Phi(x)  = \langle x, \nabla \Phi(x) \rangle,
\end{equation}
and $\nabla \Phi(x)$ is the only vector for which \eqref{equality2bis} is valid. Hence
$$
\nabla U (x)=\nabla \Phi(x)
$$
for almost every $x\in\Omega$. Together with assumption $V(0)=\Phi(0)$ and lower semi-continuity, this clearly implies $\Phi=U$. 

In addition, by the argument used in the first part of this proof, equality is possible if and only if $\nabla U(\Omega) = \{U^* < \infty\}$ up to a set of zero measure, this implies the last statement of the theorem.
\end{proof}

 \begin{theorem}
 \label{BSOT}
 Let $\alpha, \beta > 0$ be numbers and $\rho_1, \rho_2$ be positive functions. Assume that $\Phi$ is a maximum point of the functional
 $
 \mathcal{BS}_{\alpha, \beta, \rho_1, \rho_2}$. 
 Then $\nabla \Phi$ is the optimal transportation pushing forward measure $\mu$ onto $\nu$, where
$$
d\mu = \frac{ e^{-\alpha \Phi} \rho_1 dx }{\int e^{-\alpha \Phi} \rho_1 dx }, \  d\nu = \frac{ e^{-\beta \Phi^*} \rho_2 dy }{\int e^{-\beta \Phi^*} \rho_2 dy }.
$$
In addition,
$$
\nabla \Phi(\Omega)=\{ \Phi^{*} < \infty\}
$$ 
up to a set of Lebesgue measure zero.
 \end{theorem}
 \begin{proof}
 The result follows immediately from Proposition \ref{key comparison}.
 \end{proof}

 \begin{remark}
     Another remarkable monotonicity property of the Blaschke--Santal\'o functional in terms of a Gaussian diffusion semigroup has been recently obtained by Nakamura and Tsuji in \cite{NakTs}. 
     This result provides an alternative and purely analytical proof of the classical Blaschke--Santal{\'o} inequality.
 \end{remark}

\subsection{Brascamp--Lieb type inequality for maximizers}
\label{BLformaximizers}

In this subsection we derive a partial differential inequality  for the maximum point of $\bs_{\alpha,\beta,\rho_1,\rho_2}$, and we observe that this is an inequality of the Brascamp--Lieb type.

\begin{proposition}
\label{BL1906}
Let $\Phi$ be the maximum point of $\bs_{\alpha,\beta,\rho_1, \rho_2}$. Assume, in addition, that $\Phi$ is strictly convex and twice continuously differentiable inside of $\{ \Phi < \infty\}$. Then the measure $\mu$ with density
$$ 
\frac{e^{-\alpha \Phi} \rho_1}{\int e^{-\alpha \Phi} \rho_1 dx}
$$ 
satisfies 
$$
{\rm Var}_{\mu} f \le \frac{1}{\alpha+\beta} \int \langle (D^2 \Phi)^{-1} \nabla f, \nabla f \rangle d \mu.
$$
for every function $f\in C^1(\R^n)$. 
\end{proposition}
\begin{proof} Take $f$ satisfying 
\begin{equation}
\label{1903}
    \int f e^{-\alpha \Phi} \rho_1 dx=0.
\end{equation}
    One has
    $$
    \frac{1}{\alpha} \log \int e^{-\alpha \Phi_{\varepsilon}} \rho_1 dx
    + \frac{1}{\beta} \log
\int e^{-\beta \Phi^*_{\varepsilon}} \rho_2  dy
 \le
\frac{1}{\alpha} \log \int e^{-\alpha \Phi} \rho_1 dx
    + \frac{1}{\beta} \log
\int e^{-\beta \Phi^*} \rho_2  dy, 
    $$
    where $\Phi_{\varepsilon} = \Phi + \varepsilon f$.
    By Taylor expansion of $e^x$ and (\ref{1903})
\begin{eqnarray*}
\log\int e^{-\alpha \Phi_{\varepsilon}} \rho_1 dx
&=& \log  \int e^{-\alpha \Phi} \Bigl(1 - \varepsilon \alpha f + \frac{\varepsilon^2 \alpha^2}{2} f^2 + o(\varepsilon^2) \Bigr) \rho_1 dx\\
&=& \log  \int e^{-\alpha\Phi} \Bigl(1  + \frac{\varepsilon^2 \alpha^2}{2} f^2 + o(\varepsilon^2) \Bigr) \rho_1 dx\\
&=& \log  \int e^{-\alpha \Phi}  \rho_1 dx
+ \frac{\varepsilon^2 \alpha^2}{2} \frac{ \int f^2 e^{-\alpha\Phi} \rho_1 dx}{\int e^{-\alpha \Phi} \rho_1 dx} + o(\varepsilon^2),
\end{eqnarray*}
     Lemma \ref{v*inf} implies
 \begin{align*}
\log  & \int e^{-\beta \Phi^*_{\varepsilon}} \rho_2 dx
 = \log  \int e^{-\beta \Phi^*} e^{\beta \bigl( \varepsilon f(\nabla \Phi^*) -\frac{\varepsilon^2}{2} \langle D^2 \Phi^* \nabla f(\nabla \Phi^*), \nabla f(\nabla \Phi^*) \rangle + o(\varepsilon^2) \bigr)} \rho_2 dy
 \\& 
 = \log  \int e^{-\beta\Phi^*} \Bigl( 1+ \varepsilon \beta f(\nabla \Phi^*) -\frac{\varepsilon^2 \beta}{2} \langle D^2 \Phi^* \nabla f(\nabla \Phi^*), \nabla f(\nabla \Phi^*) \rangle + \frac{\varepsilon^2 \beta^2}{2} f^2(\nabla \Phi^*) + o(\varepsilon^2)  \Bigr) \rho_2 dy.
 \end{align*}
 Since $\nabla \Phi^*$ sends the measure $\frac{e^{-\beta\Phi^*}\rho_2 dy}{\int{e^{-\beta \Phi^*}\rho_2}dy}$ to the measure $\frac{e^{-\alpha \Phi}\rho_1 dx}{\int{e^{-\alpha \Phi}\rho_1} dx}$, one gets 
$$
\int f(\nabla \Phi^*)  e^{-\beta \Phi^*} \rho_2 dy = C \int f e^{-\alpha \Phi} \rho_1 dx =0.
$$ 
Thus
 \begin{align*}
\log  & \int e^{-\beta \Phi^*_{\varepsilon}} \rho_2 dy
 = \log  \int e^{-\beta\Phi^*} \Bigl(1 -\frac{\varepsilon^2\beta}{2} \langle D^2 \Phi^* \nabla f(\nabla \Phi^*), \nabla f(\nabla \Phi^*) \rangle + \frac{\varepsilon^2\beta^2}{2} f^2(\nabla \Phi^*) + o(\varepsilon^2)  \Bigr) \rho_2 dx.
 \end{align*}
Let us denote for brevity
$$
Q(f) = - \frac{\beta}{2} \langle D^2 \Phi^* \nabla f(\nabla \Phi^*), \nabla f(\nabla \Phi^*) \rangle + \frac{\beta^2}{2} f^2(\nabla \Phi^*). 
$$
Applying the Taylor expansion of $\log$ one gets
 \begin{align*}
\log  \int e^{-\beta \Phi^*_{\varepsilon}} \rho_2 dy
  & = \log  \int e^{-\beta \Phi^*} \big(1 + \varepsilon^2Q(f)  + o(\varepsilon^2)\big) \rho_2 dy 
  \\& = \log  \int e^{-\beta \Phi^*} \rho_2 dy 
  +  \log  \Bigl[\frac{1}{\int e^{-\beta \Phi^*} \rho_2 dy }\int e^{-\beta \Phi^*} \big(1 + \varepsilon^2Q(f)  + o(\varepsilon^2)\big) \rho_2 dy\Bigr] 
  \\&=
  \log  \int e^{-\beta \Phi^*} \rho_2 dy 
  +  \log  \Bigl[
  1 + 
  \frac{\varepsilon^2}{\int e^{-\beta \Phi^*} \rho_2 dy }\int e^{-\beta \Phi^*}  Q(f)   \rho_2 dy + o(\varepsilon^2)\Bigr] 
  \\& 
=  \log  \int e^{-\beta \Phi^*} \rho_2 dy 
  +   
  \frac{\varepsilon^2}{\int e^{-\beta \Phi^*} \rho_2 dy }\int Q(f) e^{-\beta \Phi^*}    \rho_2 dy + o(\varepsilon^2).
\end{align*}
Changing variables we observe that 
$$
\frac{1}{\int e^{-\beta \Phi^*} \rho_2 dy }\int Q(f)e^{-\beta \Phi^*}     \rho_2 dy 
= \frac{1}{2\int e^{-\alpha\Phi} \rho_1 dx} \int \Big( \beta^2 f^2 - 
  \beta \langle (D^2 \Phi)^{-1} \nabla f, \nabla f \rangle 
 \Bigr) e^{-\alpha\Phi} \rho_1 dx.
$$
Hence
$$
\log  \int e^{-\beta \Phi^*_{\varepsilon}} \rho_2 dy
 = \log \int e^{-\beta\Phi^*} \rho_2 dy +  \frac{\varepsilon^2}{2\int e^{-\alpha\Phi} \rho_1 dx} \int \Big( \beta^2 f^2 - 
  \beta \langle (D^2 \Phi)^{-1} \nabla f, \nabla f \rangle 
 \Bigr) e^{-\alpha\Phi} \rho_1 dx + o(\varepsilon^2). $$
Finally, one gets the relation
 \begin{align*}
  \frac{1}{\alpha} \log \int e^{-\alpha \Phi_{\varepsilon}} \rho_1 dx
    + \frac{1}{\beta} \log
\int e^{-\beta \Phi^*_{\varepsilon}} \rho_2  dy
& =  \frac{1}{\alpha} \log \int e^{-\alpha \Phi} \rho_1 dx
    + \frac{1}{\beta} \log
\int e^{-\beta \Phi^*} \rho_2  dy \\& + \frac{\varepsilon^2}{\int e^{-\alpha\Phi} \rho_1 dx} \Bigl(  \int \Big( \frac{\alpha+\beta}{2}f^2 - 
  \frac{1}{2} \langle (D^2 \Phi)^{-1} \nabla f, \nabla f \rangle 
 \Bigr) e^{-\alpha\Phi} \rho_1 dx \Bigr) + o(\varepsilon^2)
 \end{align*}
 and the claim follows.
\end{proof}

 Let $p>1$ and $V$ be a $p$-homogeneous convex function. The latter means, in particular,
 \begin{equation}
 \label{phomV}
 V = \frac{1}{p-1} V^*(\nabla V).
 \end{equation}
 Consider the functional 
 $$
 \int e^{-\Phi} dx  \cdot \Bigl( \int e^{-{\frac{1}{p-1}}  \Phi^*} \rho dy \Bigr)^{ p-1} =  \int e^{-\Phi} dx  \cdot \Bigl( \int e^{-{\frac{1}{p-1}}  \Phi^*(\nabla V)} dx \Bigr)^{ p-1},
 $$
 where $\rho dy$ is the image of the Lebesgue measure under $\nabla V$, in particular 
 $$
 \rho = \det D^2 V^*.
 $$
 Note that 
 $\nabla V$ is the optimal transportation mapping of
 $
 \frac{e^{-V} dx}{\int e^{-V} dx}
 $
 onto
 $
 \frac{ e^{-{\frac{1}{p-1}}  V^*} \rho dy }{\int e^{-{\frac{1}{p-1}}  V^*} \rho dy }
 $
 This follows from (\ref{phomV})
 and the fact that $\rho dy$ is the image of the Lebesgue measure under $\nabla V$.

Thus we get that $V$ is a natural candidate to maximize $\bs_{p, V}$. The expected inequality reads as:
$$
  \int e^{-\Phi} dx  \cdot \Bigl( \int e^{-{\frac{1}{p-1}}  \Phi^*(\nabla V)} dx \Bigr)^{ p-1}
  \le \Bigl( \int e^{-V} dx \Bigr)^{p}. 
 $$
 
In particular, $V$ satisfies the corresponding Euler--Lagrange equation. The second order condition obtained in Proposition \ref{BL1906} means that if $V$ is a maximizer, it must satisfy the inequality.
  $$
 {\rm Var}_{\mu} f \le \Bigl( 1 - \frac{1}{p} \Bigr)\int \langle ( D^2 V)^{-1} \nabla f, \nabla f \rangle d\mu, 
 $$ 
 where $
 \mu = \frac{ e^{-V}  dx }{\int e^{-V} dx }$.

\subsection{Homogeneity of maximizers of $\bs$-functional}

In this section we consider a smooth (say, $C^3(\R^n)$) and strictly convex function $V\colon\R^n\to\R$; we will refer to $V$ as a {\em potential}.
It will be assumed throughout that 
$V$ is a $p$-homogeneous convex function for some fixed $p>1$. 
The proof of the following properties will be left to the reader as an exercise.

\begin{proposition}\label{properties of p-homogeneous functions} Let $V\in C^3(\R^n)$ be a convex, even $p$-homogeneous function.
Then $V$ verifies the following properties. 
\begin{enumerate}
\item $\langle \nabla V(x), x \rangle = p V(x)$;
\item $(p-1) V = V^*(\nabla V)$;
\item $(p-1) \nabla V(x) = D^2 V(x) \cdot x$;
\item for every vector $e$ one has
$(D^2 V(x))_e \cdot x = (p-2) D^2 V(x) \cdot e$
\end{enumerate}
Here $(\cdot)_e$ indicates partial differentiation along $e$:
$(D^2 V(x))_e$ is the symmetric matrix with entries $\frac{\partial^3 V(x)}{\partial e_i \partial e_j \partial e} $.
\end{proposition}

We consider the Blaschke-Santal\'o functional
$$
\bs_{p,V}(\Phi)=
 \int e^{-\Phi} dx  \cdot \Bigl( \int e^{-{\frac{1}{p-1}}  \Phi^*} \rho dy \Bigr)^{ p-1} =  \int e^{-\Phi} dx  \cdot \Bigl( \int e^{-{\frac{1}{p-1}}  \Phi^*(\nabla V)} dx \Bigr)^{ p-1},
 $$
 where $\rho dy$ is the image of the Lebesgue measure under $\nabla V$. Note that 
 $\nabla V$ is the optimal transportation mapping of
 $
 \frac{e^{-V} dx}{\int e^{-V} dx}
 $ 
 onto
 $
 \frac{ e^{-{\frac{1}{p-1}}  V^*} \rho dy }{\int e^{-{\frac{1}{p-1}}  V^*} \rho dy }
 $
 This follows from (\ref{phomV})
 and the fact that $\rho dy$ is the image of the Lebesgue measure under $\nabla V$.
 
 \begin{theorem}
 \label{p-hom-solutions}
Any even  maximum point  of the functional
 $$
 \bs_{p,V}(\Phi) = \int e^{-\Phi} dx  \cdot \Bigl( \int e^{-{\frac{1}{p-1}}  \Phi^*(\nabla V)} dx \Bigr)^{ p-1}
 $$
 is $p$-homogeneous (up to addition of a constant).
 \end{theorem}

 We start with some preliminary considerations.  If $\Phi$ is a maximum point, Theorem \ref{BSOT} implies that:
\begin{enumerate}
\item
$\nabla \Phi$ is the optimal transportation mapping between the measures
$$
\mu = \frac{e^{-\Phi} dx}{\int_{\mathbb{R}^n} e^{-\Phi} dx}\  \ {\rm and} \ 
\nu = \frac{e^{-\frac{1}{p-1}\Phi^*} \det D^2 V^* dy}{\int_{\mathbb{R}^n} e^{- \frac{1}{p-1}\Phi^*} \det D^2 V^* dy},
$$
\item
$\nabla \Phi(\mathbb{R}^n) = \{ \Phi^* < \infty\}$.
\end{enumerate}

Note that $\mathcal{BS}_{p,V}(\Phi) = \mathcal{BS}_{p,V}(\Phi + a) $ for any $a \in \mathbb{R}$. 
Thus adding appropriate $a$ we can assume
$\Phi(0)=0$. This implies by symmetry and convexity of $\Phi$ that $\Phi \ge 0$. But then
$
\Phi^*(0)=\sup_{y} (-\Phi(y)) =0.
$
Thus
without loss of generality we assume
$$
\Phi(0)=\Phi^*(0)=0.
$$

First we observe that $\Phi, \Phi^*$ are smooth on $\mathbb{R}^n \setminus \{0\}$ provided $V$ is smooth on $\mathbb{R}^n \setminus \{0\}$ (we can not 
assume that $V$ is smooth on entire  $\mathbb{R}^n $, because $V$ is $p$-homogeneous).
Indeed, to prove this we apply local H{\"o}lder estimates for solution to the 
Monge--Amp\`ere equation. We refer to \cite{KK-eigenvalues}, proof of Lemma 5.2. We choose a neighborhood $U$ of a point $y_0 \ne 0$ which does not contain the origin. Using that $\det D^2 V^*$ is smooth inside of $U$, we consider equation
$$
C \det D^2 V^* e^{-\frac{1}{p-1}\Phi^*} = e^{-\Phi(\nabla \Phi^*)} \det D^2 \Phi^*,
$$
where $C$ is the corresponding normalizing constant,
and apply the Forzani--Maldonado estimate (see  \cite{ForMal}) to ensure that $\nabla \Phi^*$ is locally H\"older and then the Trudinger--Wang estimates (see \cite{TrudWang}) to prove higher regularity. See details in \cite{KK-eigenvalues}. Applying standard bootstrapping arguments we can conclude that $\Phi^*$ is smooth on $\Omega^* \setminus \{0\}$
and the similar statement holds for $\Phi$.

Next we  prove a lemma about behavior of $\Phi$ near the boundary of $\{\Phi < \infty\}$. 

We recall that $\Omega, \Omega^*$ are defined in (\ref{1903}).

\begin{lemma}
\label{subd-boundary}
    $\partial \Phi(x)$ is empty for every $x \in \partial \Omega$ and $\partial \Phi^*(x)$ is empty for every $y \in \partial \Omega^*$. 
    In particular, $|\nabla \Phi(x_n)| \to \infty$ for every sequence $x_n \to x$, where $x \in \partial \Omega$.
    \end{lemma}
    \begin{proof}

        Assume that $a \in \partial \Phi(x)$, $x \in \partial \Omega$ (this means in particular, that $\Phi(x) <\infty$). Then  $a + tn \in \partial \Phi(x)$ for every  $t \ge 0$, where  $n$ is the outer normal to $\partial \Omega$ at $x$.
        Indeed, to this end we need to check that for all $y$ one has 
        $$
\Phi(y) - \Phi(x) \ge \langle y-x, a + nt\rangle.
        $$
        If $y \notin \Omega$, then $\Phi(y) = +\infty$ and the claim is trivial.
        If $y \in \Omega$, then 
        $\Phi(y) - \Phi(x) \ge \langle y-x,a \rangle$ because $a \in \partial \Phi(x)$ and $\langle y -x, n \rangle \le 0$ because $n$ is the outer normal. This implies that $a + tn \in \partial \Phi(x)$.

        On one hand, we note that by strict convexity of  $\Phi$ one has:  $L_a = \{ a + tn \in \partial \Phi(x), t >0 \} \subset \mathbb{R}^n \setminus \nabla \Phi(\Omega)$.
        On the other hand, we note that
        $$
        \Phi^*(a + nt) = \langle x, a + nt \rangle - \Phi(x) < \infty.
        $$
        Hence $\Phi^*$ is finite on $L_a$ and by convexity $\Phi^*$ is finite in some neighborhood $\tilde{L}$ of $L_a$. Hence there exists an open set $\tilde L \subset \mathbb{R}^n \setminus \nabla \Phi(\Omega)$ with the property $L \subset \Omega^*$. But this contradicts to the fact that $\nabla \Phi(\Omega) = \Omega^*$ (up to a set of  measure zero), see Theorem \ref{BSOT}. 
    \end{proof}

In what follows we consider function
$$
W(x) = \langle x, \nabla \Phi(x) \rangle - p\Phi(x).
$$
In a sense, $W$ ``measures the homogeneity" of $\Phi$. If $\Phi$ is $p$-homogeneous, then $W=0$ according to Proposition \ref{properties of p-homogeneous functions}.

We work with the operator
$L$
described in Section \ref{prem} (see (\ref{generator})).
We recall that $L$ is the symmetric operator associated with the Riemannian metric $D^2 \Phi(x)$ and measure $\mu$. 
Using (\ref{generator}), the presentation of the potential $W$ of the push-forward measure $\nu = \mu \circ (\nabla \Phi)^{-1}$ is
$$
W = c + \frac{\Phi^*}{p-1} - \log \det D^2 V^*
$$
and the relation $\nabla \Phi^*(\nabla \Phi(x))=x$ we get the following formula for $L$:
$$
L f = {\rm Tr} \Bigl[ (D^2 \Phi)^{-1} D^2 f\Bigr] - \Big\langle \nabla f,  \frac{x}{p-1} - [\nabla \log \det D^2 V^*] \circ \nabla \Phi \Big\rangle.
$$ Since $L$ is symmetric with respect to $\mu$, one has the following integration by parts formula: for every smooth $g$ vanishing in a neighborhood of $\partial\Omega \setminus \{0\}$ one has 
\begin{equation}
\label{ibpforL}
- \int Lf g d \mu = \int \langle (D^2 \Phi)^{-1} \nabla f, \nabla g \rangle d \mu. 
\end{equation}
The following lemma is proved by direct computations (and differentiation of change of variables).

\begin{lemma} \label{LW0} Let $\Phi$ be any even maximum point of $\mathcal{BS}_{p,V}$. The following equation holds for all $x \in \Omega \setminus \{0\}$:
$$
LW(x) =0.
$$
\end{lemma}

\begin{proof}
It was explained in the beginning of the section that $\Phi$ is smooth  inside of $\Omega \setminus \{0\}$. Since $\Phi$  is the optimal transportation of $\mu$ onto $\nu$, the following change of variables formula holds on $\Omega \setminus \{0\}$
\begin{equation}
    \label{1903chvar}
\frac{e^{-\Phi} }{\int_{\mathbb{R}^n} e^{-\Phi} dx} = \frac{e^{-\frac{1}{p-1}\Phi^*(\nabla \Phi)} \det D^2 V^*(\nabla \Phi)}{\int_{\mathbb{R}^n} e^{- \frac{1}{p-1}\Phi^*} \det D^2 V^* dy}.
\end{equation}

One has
$$
\nabla W(x) = D^2 \Phi(x)\cdot x + (1-p) \nabla\Phi(x).
$$
For every $e\in\mathbb{R}^n$
$$
D^2 W \cdot e = (2-p) D^2 \Phi \cdot e + (D^2 \Phi)_e \cdot x,
$$
{(where $(\cdot)_e$ indicates partial differentiation along $e$).} Therefore, for every $e_i,e_j\in\R^n$,
$$
W_{e_i e_j} = \langle D^2 W \cdot e_i, e_j \rangle = (2-p) \Phi_{e_i e_j} + \langle \nabla \Phi_{e_i e_j}, x \rangle = (2-p) \Phi_{e_i e_j}  + \sum_{k=1}^n \Phi_{e_i e_j e_k} x_k.
$$
Thus
$$
LW = \sum_{k=1}^n \langle  (D^2 \Phi)^{-1} (D^2 \Phi)_{e_k} \cdot x, e_k \rangle - \langle D^2 \Phi \cdot x + (1-p) \nabla \Phi, \frac{x}{p-1} - [\nabla \log \det D^2 V^*] \circ \nabla \Phi  \rangle + (2-p)n.
$$
Let $e_k$, $k\in\{1,\dots,n\}$, be a basis of eigenvectors for $D^2 \Phi(x)$ at $x$:
$$
D^2 \Phi(x) e_k = \lambda_k e_k,
$$
where $\lambda_1,\dots,\lambda_n$ are the eigenvalues of $D^2\Phi(x)$. One gets the following expression for $LW$:
\begin{equation}
    \label{lw}
    LW = \sum_{k, j=1}^n \frac{\Phi_{e_k e_k e_j} \cdot x_j}{\lambda_k}   - \langle D^2 \Phi \cdot x + (1-p) \nabla \Phi, \frac{x}{p-1} - [\nabla \log \det D^2 V^*] \circ \nabla \Phi  \rangle + (2-p)n.
\end{equation}
Taking logarithm  of (\ref{1903chvar}) we get
$$
\Phi = \frac{\Phi^*(\nabla \Phi)}{p-1} - \log \det D^2 \Phi -\log \det D^2 V^*(\nabla \Phi) + C
$$
for some constant $C$. Differentiating this equation along the vector field $x$ we get
$$
\nabla \Phi 
=  \frac{1}{p-1} D^2 \Phi \cdot x  - {\sum_{j=1}^n {\rm Tr} (D^2 \Phi)^{-1} (D^2 \Phi)_{e_j} \cdot e_j} -  D^2 \Phi \cdot  [\nabla \log \det D^2 V^*] \circ \nabla \Phi. 
$$
 Here we have used the differentiation formula for determinants:
$
\partial_e \log \det D^2 \Phi = {\rm Tr} (D^2 \Phi)^{-1} (D^2 \Phi)_{e}. 
$

Whence
$$
\langle \nabla \Phi, x \rangle 
=  \frac{1}{p-1}\langle D^2 \Phi \cdot x, x \rangle  - \sum_{k,j=1}^n  \frac{\Phi_{e_k e_k e_j} x_j}{\lambda_k} \cdot x_j  -  \langle D^2 \Phi \cdot  [\nabla \log \det D^2 V^*] \circ \nabla \Phi, x \rangle. 
$$
Putting this expression into (\ref{lw}) one gets
\begin{align*}
LW =& \frac{1}{p-1}\langle D^2 \Phi \cdot x, x \rangle  -  
\langle \nabla \Phi, x \rangle  -  \langle D^2 \Phi \cdot  [\nabla \log \det D^2 V^*] \circ \nabla \Phi, x \rangle
\\& - \langle D^2 \Phi \cdot x + (1-p) \nabla \Phi, \frac{x}{p-1} - [\nabla \log \det D^2 V^*] \circ \nabla \Phi  \rangle
\\=& (1-p)\langle \nabla \Phi,  [\nabla \log \det D^2 V^*] \circ \nabla \Phi  \rangle + n(2-p)
\\=& (1-p)\langle x,  \nabla \log \det D^2 V^* \rangle \circ \nabla \Phi + n(2-p).
    \end{align*}
The claim follows from the observation that $\det D^2 V^*$ is a $n \frac{2-p}{p-1}$-homogeneous function, and therefore  $$\langle x,  \nabla \log \det D^2 V^* \rangle =  \frac{\langle x,  \nabla \det D^2 V^* \rangle}{\det D^2 V^*} 
= \frac{n \frac{2-p}{p-1} \det D^2 V^*}{\det D^2 V^*} =n \frac{2-p}{p-1}.$$
\end{proof}

Before giving a complete proof of Theorem \ref{p-hom-solutions}, let us explain the main idea.
Our aim is to conclude that $W=0$ using the equation $LW=0$ proved in Lemma \ref{LW0}. Applying integration by parts (\ref{ibpforL}) one  can formally write
$$
0 = - \int LW \cdot f(W) d\mu =  -  \int f'(W) \langle (D^2\Phi)^{-1} \nabla W, \nabla W \rangle d\mu.
$$
Here $f$ is any monotone function. From this identity we  formally conclude that $W=0$ (note that $W(0)=0$).
However, an application of (\ref{ibpforL}) with functions without compact support inside 
of $\Omega$ requires accurate  work with the boundary terms.

We overcome this difficulties as follows:
we create auxiliary functions of the form
$$
f(W) g_N(\Phi)
$$
vanishing on the $\partial \Omega$, apply (\ref{ibpforL}) and approximate the constant unit function by $g_N(\Phi)$. In the limit we get the desired relation $\int f'(W) \langle (D^2\Phi)^{-1} \nabla W, \nabla W \rangle d\mu$ for a sufficiently broad class of monotone functions.
Along the way we have to estimate some error terms using geometric properties of $\Phi$ and $\Phi^*$.

\begin{proof}[Proof of Theorem \ref{p-hom-solutions}]
We prove that $W$ is constant.
Set:
$$
B = \{ x \in \partial \Omega, \Phi(x) < \infty\},
$$
$$
B_{\infty} = \{ x \in \partial \Omega,  \Phi (x) = + \infty \}.
$$
Note that $W =+\infty$ on $B$ by Lemma \ref{subd-boundary}, because $|\nabla \Phi|(x_n) \to \infty$ as $x_n \to x \in \partial \Omega$.

In what follows we consider a smooth, convex, non negative and decreasing function $f$, which vanishes on $[a,+\infty)$ for some $a$. In addition, we assume that $-f' \le C$ for some $C>0$.
Using the condition $L W =0$, which comes from the Lemma \ref{LW0}, one gets
\begin{equation}
    \label{LfW}
L f(W) = f^{''} (W) \langle (D^2 \Phi)^{-1} \nabla W, \nabla W \rangle.
\end{equation}
Note that $f(W)$ vanishes on some neighborhood of $B$.

Let $g_{N,\varepsilon}$ be a family of smooth nonnegative functions, satisfying the following assumptions:
$g_{N,\varepsilon}(t) =0$ for $t \le \frac{\varepsilon}{2}$, $g_{N,\varepsilon}(t)$ is increasing
on $[\frac{\varepsilon}{2},\varepsilon]$, $g_{N,\varepsilon}=1$ on $[\varepsilon,N]$, $g_{N,\varepsilon}$
is decreasing on $[N,2N]$, $g_{N,\varepsilon}=0$ on $[2N,\infty)$.

Note that $g_{N,\varepsilon}(\Phi)$ vanishes in  neighborhoods of  $B_{\infty}$ and of the origin.
Hence the product $f(W) g_{N,\varepsilon}(\Phi)$ vanishes on $\partial \Omega$ and has compact support not containing the origin. Thus one can apply equations (\ref{LfW}) and (\ref{ibpforL}):
\begin{align*}
\int f^{''}(W) g_{N,\varepsilon}(\Phi) \langle (D^2 \Phi)^{-1} \nabla W, \nabla W \rangle d\mu
& =  \int L f(W) g_{N,\varepsilon}(\Phi) d \mu 
= - \int \langle (D^2 \Phi)^{-1} \nabla [f(W)], \nabla[g_{N,\varepsilon}(\Phi)] \rangle d\mu
\\&- \int f'(W) g'_{N,\varepsilon}(\Phi) \langle (D^2 \Phi)^{-1} \nabla W,  \nabla \Phi \rangle d \mu.
\end{align*}
For any given $g_{n,\varepsilon}$ let us consider two monotone functions:
$
g^1_{N,\varepsilon}
$ and $g^2_{N,\varepsilon}$ such that
$$g^1_{N,\varepsilon} = g_{N,\varepsilon} \ {\rm on} \ (-\infty,N] \ \ {\rm and} \ \ g^1_{N,\varepsilon}=1 \ {\rm on} \ [N,+\infty)$$ $$g^2_{N,\varepsilon} = 1 \ {\rm on} \ (-\infty,N] \ \ {\rm and} \ \  g^2_{N,\varepsilon}=g_{N,\varepsilon}  \ {\rm on} \ [N,+\infty).$$ Note that $g^1_{N,\varepsilon}$ is non-decreasing and $g^2_{N,\varepsilon}$ is non-increasing.
Note that 
$$
g'_{N,\varepsilon}
= (g^1_{N,\varepsilon})' + (g^2_{N,\varepsilon})'.
$$
In particular, the following  holds
\begin{align}
\label{2308}
\int f^{''}(W) g_{N,\varepsilon}(\Phi) \langle (D^2 \Phi)^{-1} \nabla W, \nabla W \rangle d\mu
=&   - \int f'(W) (g^1_{N,\varepsilon})'(\Phi) \langle (D^2 \Phi)^{-1} \nabla W,  \nabla \Phi \rangle d \mu
\\& 
\nonumber
 + \int f(W)  L \bigl( g^2_{N,\varepsilon}(\Phi) \bigr)  d \mu.
\end{align}
We will choose $g_{N,\varepsilon}$ in such a way that $\lim_{\varepsilon, N} g_{N,\varepsilon}=1$ and the limit of the right-hand side of (\ref{2308}) is not positive.

We estimate the first term of the right-hand side of (\ref{2308})
Using that $f$ is decreasing and $g^1_{N,\varepsilon}$ is increasing we observe
\begin{align*}
     - \int f'(W) (g^1_{N,\varepsilon})'(\Phi) \langle (D^2 \Phi)^{-1} \nabla W,  \nabla \Phi \rangle d \mu 
    =&  - \int f'(W) (g^1_{N,\varepsilon})'(\Phi) \langle x,  \nabla \Phi(x) \rangle d \mu
    \\& + (p-1) \int f'(W) (g^1_{N,\varepsilon})'(\Phi) \langle (D^2 \Phi)^{-1} \nabla \Phi ,  \nabla \Phi \rangle d \mu
    \\ \le& - \int f'(W) (g^1_{N,\varepsilon})'(\Phi) \langle x,  \nabla \Phi(x) \rangle d \mu
    \\ \le& C \int (g^1_{N,\varepsilon})'(\Phi) \langle x,  \nabla \Phi(x) \rangle d \mu.
\end{align*}
Using approximations one can relax smoothness assumption and suppose that
$g^1_{N,\varepsilon}(t) = \frac{2}{\varepsilon} t - 1$ on $[\varepsilon/2, \varepsilon]$ for all $N$. Thus
$$
- \int f'(W) (g^1_{N,\varepsilon})'(\Phi) \langle (D^2 \Phi)^{-1} \nabla W,  \nabla \Phi \rangle d \mu  
\le \frac{2C}{\varepsilon \int e^{-\Phi}dx}
\int_{\{\Phi \le \varepsilon\}} 
\langle x,  \nabla \Phi(x) \rangle  e^{-\Phi}dx.
$$
Then we estimate using integration by parts of $\{\Phi \le \varepsilon\}$ (here $\mathcal{H}^{n-1}$ is the $n-1$-dimensional Haussdorff measure
\begin{align*}
\frac{1}{\varepsilon} \int_{\{\Phi \le \varepsilon\}} 
\langle x,  \nabla \Phi(x) \rangle  e^{-\Phi}dx
& \le 
\frac{1}{\varepsilon} \int_{\{\Phi \le \varepsilon\}} 
\langle x,  \nabla \Phi(x) \rangle dx 
= \frac{1}{\varepsilon} \Bigl( - n \int_{\{\Phi \le \varepsilon\}} 
\Phi(x)  dx 
+ \int_{\Phi=\varepsilon}
\Phi \langle x, \nu \rangle 
 d\mathcal{H}^{n-1} \Bigr)
\\&  \le \int_{\Phi=\varepsilon}
 \langle x, \nu \rangle 
 d\mathcal{H}^{n-1}
 = n |\{\Phi \le \varepsilon\}|.
\end{align*}
Thus we get
$$\lim_{\varepsilon \to 0} \bigl( - \int f'(W) (g^1_{N,\varepsilon})'(\Phi) \langle (D^2 \Phi)^{-1} \nabla W,  \nabla \Phi \rangle d \mu \bigr) \le 0.$$

 Let us analyse the second term 
 $$
 \int f(W)  L \bigl( g^2_{N,\varepsilon}(\Phi) \bigr)  d \mu
 $$ of the right-hand side of (\ref{2308}).

 The function $g^2_{N,\varepsilon}$ will not depend on $\varepsilon$ and we will write 
 $$
g_N = g^2_{N,\varepsilon}.
 $$
Recall that  $g_N$ is supposed to be non-increasing and satisfying
 $$
0 \le g_N \le 1, \ g_N|_{(-\infty,N]} =1, \ g_N|_{[2N,+\infty)} =0.
$$
We fix a decreasing smooth function $\psi$ satisfying $\psi(t)=1$ for $t\le 0$ and $\psi(t)=0$ for $t\ge 1$.
Then we set
$$
g_N(t) =  \psi\bigl( \frac{t}{N}-1 \bigr) 
$$
for $t \ge N$.

We observe that for some $C>0$ one has
$$
|g^{'}_N(t)| \le \frac{C}{N}, \ \ |g_N^{''}(t)| \le \frac{C}{N^2} e^{-\frac{t}{N}}.
$$

One has
$$
L g_N(\Phi) = g'_N(\Phi) \Bigl( n(p-1) - \frac{1}{p-1}\langle \nabla \Phi(x), x \rangle \Bigr)  + g^{''}_N(\Phi) \langle (D^2 \Phi)^{-1} \nabla \Phi, \nabla \Phi \rangle : = I + II.
$$
On the set  $\{x: f(W(x))>0\}$, one has $\langle x, \nabla \Phi(x) \rangle \le a + p\Phi(x)$. 
Finally, we obtain
$$
I = f(W) |g'_N(\Phi)| \Bigl| n(p-1) - \frac{1}{p-1}\langle \nabla \Phi(x), x \rangle \Bigr|  \le \frac{c_1}{N} (1 + |W|) \le \frac{c_2}{N} (1 + \Phi) .
$$
for some $c_1, c_2 >0$.
Since $\Phi \in L^1(\mu)$ we immediately conclude that 
$\int f(W) I e^{-\Phi} dx \to 0$
as $N \to \infty$.

Next we use the bound $
|g^{''}_N(t)| \le \frac{C}{N^2} e^{-\frac{t}{N}}$:
$$
f(W) |g^{''}_N(\Phi)| \le \frac{c}{N^2}  e^{-\frac{\langle x, \nabla \Phi(x) \rangle}{pN}}
$$
and
\begin{align*}
\int  f(W) II e^{-\Phi} dx =& \int f(W) |g^{''}_N(\Phi)|   \big\langle \big(D^2 \Phi\big)^{-1} \nabla \Phi, \nabla \Phi \big\rangle e^{-\Phi} dx\\
\le& \frac{c}{N^2}  \int  \big\langle \big(D^2 \Phi\big)^{-1} \nabla \Phi, \nabla \Phi \big\rangle e^{-\frac{\langle x, \nabla \Phi(x) \rangle}{pN}} e^{-\Phi(x)} dx
\\  =&  \frac{\int e^{-\Phi} dx }{\int e^{-\frac{\Phi^*}{p-1}} dy } \frac{c}{N^2}  \int  \langle D^2 \Phi^* y, y \rangle e^{-\frac{\langle y, \nabla \Phi^*(y) \rangle}{pN}} e^{-\frac{\Phi^*(y)}{p-1}} dy.
\end{align*}
Finally, we get that for some constant $d$
$$
\int f(W) |g^{''}_N(\Phi)|   \langle (D^2 \Phi)^{-1} \nabla \Phi, \nabla \Phi \rangle e^{-\Phi} dx
\le - \frac{d}{N} \int  \big\langle \nabla   e^{-\frac{\langle y, \nabla \Phi^*(y) \rangle}{pN}}, y \big\rangle  e^{-\frac{\Phi^*(y)}{p-1}} dy.
$$
Integrating by parts we get 
\begin{align*}
\int & f(W) |g^{''}_N(\Phi)|   \big\langle (D^2 \Phi)^{-1} \nabla \Phi, \nabla \Phi \big\rangle e^{-\Phi} dx \le - \frac{d}{N} \int_{\Omega^*}  \big\langle \nabla   e^{-\frac{\langle y, \nabla \Phi^*(y) \rangle}{pN}}, y \big\rangle  e^{-\frac{\Phi^*(y)}{p-1}} dy
\\&
= \frac{d}{N}  \int_{\Omega^*} \Bigl(n - \frac{\langle \nabla \Phi^*(y), y \rangle}{p-1} \Bigr)  e^{-\frac{\langle y, \nabla \Phi^*(y) \rangle}{pN}}  e^{-\frac{\Phi^*(y)}{p-1}} dy
- \frac{d}{N} \int_{\partial \Omega^*} \langle \nu, y \rangle e^{-\frac{\langle y, \nabla \Phi^*(y) \rangle}{pN}}  e^{-\frac{\Phi^*(y)}{p-1}} d \mathcal{H}_{n-1}.
\end{align*}
Note that by Proposition \ref{subd-boundary} (applied to $\Phi^*$) one has $\langle y, \nabla \Phi^*(y) \rangle =+\infty$ on $\Omega^*$, hence 
$$
\int_{\partial \Omega^*} \langle \nu, y \rangle e^{-\frac{\langle y, \nabla \Phi^*(y) \rangle}{pN}}  e^{-\frac{\Phi^*(y)}{p-1}} d \mathcal{H}_{n-1} =0.
$$
Then we use that $\langle \nabla \Phi^*(y), y \rangle \ge 0$ (because $\Phi^*$ is even) and obtain 
$$
\int  f(W) |g^{''}_N(\Phi)|   \big\langle \big(D^2 \Phi\big)^{-1} \nabla \Phi, \nabla \Phi \big\rangle e^{-\Phi} dx
\le \frac{dn}{N} \int  e^{-\frac{\Phi^*(y)}{p-1}} dy \to 0.
$$
Finally
\begin{align*}
\int f^{''}(W)  & \big\langle \big(D^2 \Phi\big)^{-1} \nabla W, \nabla W \big\rangle d\mu  = \lim_N \int f^{''}(W) g_N(\Phi) \big\langle \big(D^2 \Phi\big)^{-1} \nabla W, \nabla W \big\rangle d\mu
 \le  0.
\end{align*}
From this we get $\nabla W=0$.
The proof is complete.
\end{proof}

\begin{remark}
    \label{classicalBS}
    Theorem  \ref{p-hom-solutions} provides an alternative proof of the classical Blaschke--Santal{\'o} inequality (\ref{BSBall-G}) without application of symmetrization arguments. Indeed, according to Theorem \ref{existence theorem} there exists a  maximizer $\Phi$ of the classical functional $\bs$. According to Theorems \ref{BSOT} and \ref{p-hom-solutions}, $\Phi$ is a $2$-homogeneous solution to the corresponding Monge--Amp\`ere equation. Then following the arguments from \cite{Caglar-Kolesnikov-Werner} one can prove that $\Phi$ is  a quadratic function. This establishes inequality (\ref{BSBall-G}).
\end{remark}

\begin{remark}
\label{BSMA}
     Theorem \ref{p-hom-solutions} 
    can be used to establish the precise form of maximizers in (\ref{main-question-intro}) in the rotationally invariant case (see Subsection \ref{radsym}). Unfortunately, we do not know any other examples of closed-form solutions apart of the radially symmetric cases, where Theorem \ref{p-hom-solutions} can be used. We show in the following subsection that homogeneity of the maximizers allows to reduce the problem to $L^q$-Minkowski problem, which is in general ill-posed.
\end{remark}

\subsection{Relations to $L^{q}$-Minkowski problems}
\label{lqbm}

\begin{lemma}
    Let $\Phi\colon\R^n\to\R$ be $\alpha$-homogeneous and convex, for some $\alpha\ge1$. Then
    \begin{itemize}
    \item $\Phi\ge0$ in $\R^n;$
    \item $\phi=\Phi^{1/\alpha}$ is a 1-homogeneous convex function, i.e. is the support function of a convex body (containing the origin). 
    \end{itemize}
\end{lemma}
\begin{proof}  Let $u$ be a unit vector. The function $\Phi_u\colon\R\to\R$ defined by $\Phi_u(t)=\Phi(tu)$ is a $\alpha$-homogeneous convex function on the real line. Hence it must be of the form $ct^\alpha$, for some $c\ge0$. This proves that $\Phi$ is non-negative. 

Let $\phi=\Phi^{1/\alpha}$. By convexity, for every $s\ge0$
$$
\{\Phi\le s\}
$$
is convex. This proves that for every $\tau\ge0$, the set
$$
\{\phi\le\tau\}
$$
is convex. This implies that
\begin{equation}\label{quasiconvexity2}
\phi((1-t)x_0+tx_1)\le\max\{\phi(x_0),\phi(x_1)\}
\end{equation}
for every $x_0,x_1\in\R^n$, for every $t\in[0,1]$. Let $x_0,x_1\in\R^n$ and let $t\in[0,1]$. Set
$$
\bar x_0=\frac{x_0}{\phi(x_0)},\quad \bar x_1=\frac{x_1}{\phi(x_1)},\quad \bar t=\frac{t\phi(x_1)}{(1-t)\phi(x_0)+t\phi(x_1)}.
$$
Applying \eqref{quasiconvexity2} to $\bar x_0$, $\bar x_1$ and $\bar t$, and using the fact that $\phi$ is 1-homogeneous, one gets:
$$
\phi((1-t)x_0+tx_1)\le (1-t)\phi(x_0)+t\phi(x_1),
$$
that is, $\phi$ is convex.
\end{proof}

The following computational result should be known to the experts, we give the proof for completeness of the picture.

Let us fix a point $\nu \in \mathbb{S}^{n-1}$. We create a local coordinates system $\theta = (\theta_1, \theta_2, \cdots, \theta_{n-1})$
in a  neighborhood $U_0 \subset \mathbb{S}^{n-1}$  of $\nu$ as described below.
Fix any orthogonal frame $(e_1, \cdots, e_{n-1})$ in the hyperplane 
$$L = \{\theta: \theta \bot \nu \} \subset \mathbb{R}^n$$ orthogonal to $\nu$. 
We construct a local chart $\nu(\theta) : U \to \mathbb{S}^{n-1}$, where $U \subset L$ is a neighborhood of the origin in $L$.
Let $\theta = (\theta_1, \cdots, \theta_{n-1})$ be any Euclidean  coordinate system with orthonormal  frame $(e_1, \cdots, e_{n-1})$   in $L$ and  we define the following mapping $\nu(\theta)$ from $U$ to $\mathbb{S}^{n-1}$:
$$
\nu(\theta) = \frac{\nu + \theta}{\sqrt{1 + \theta^2}}.  
$$
Here  $\nu+\theta$ is the addition in $\mathbb{R}^n$ (we consider $L$ as a hyperplane in $\mathbb{R}^n$). Clearly $|\nu(\theta)|=1$ because $\nu$ and $\theta$ are orthogonal. 
Thus $\theta \to \nu(\theta)$ is a parametrization of a neighborhood of $\theta \in \mathbb{S}^{n-1}$
In particular, it is easy to check the following relations: at the point $\nu$ one has
$$
\frac{\partial^2 \nu}{\partial \theta^2_i} \Big|_{\theta=0} = - \nu, 
\ \ 
\frac{\partial^2 \nu}{\partial \theta_{i} \partial \theta_j} \Big|_{\theta=0} = 0. 
$$
This implies, in particular, that 
the Levi-Civita connection $\Gamma^{i}_{jk}$ of $\mathbb{S}^{n-1}$  vanishes at  $\nu$ and, in particular,
the spherical Hessian $\nabla^2_{\mathbb{S}^{n-1}}$ of a function $f : \mathbb{S}^{n-1} \to \mathbb{R}$ at  $\nu$ coincides with the matrix $\partial_{\theta_i \theta_j} f$. 

Next we consider the polar coordinate system $(r,\theta)$ on the cone $C_0 = \mathbb{R}_+ \times U_0$: 
$$
C_0 \ni x = r \cdot \nu(\theta).
$$

\begin{lemma}
    Let $\Phi$ be  a function defined on the neighborhood of $x_0 = r _0\cdot \nu. $ The Euclidean  Hessian of $\Phi$ has the following representation in the frame $(\nu, e_1, e_2, \cdots, e_{n-1})$
    ant the point $x_0$:
    $$
 D^2 \Phi =
  \left[ {\begin{array}{cccc}
    \Phi_{rr} & \frac{\Phi_{r \theta_1}}{r}  - \frac{\Phi_{\theta_1}}{r^2}  & \cdots & \frac{\Phi_{r \theta_{n-1}}}{r}  - \frac{\Phi_{\theta_{n-1}}}{r^2}\\
   \frac{\Phi_{r \theta_1}}{r}  - \frac{\Phi_{\theta_1}}{r^2} & a_{1,1} & \cdots & a_{1, n-1}\\
    \vdots & \vdots & a_{i,j}  & \vdots\\
   \frac{\Phi_{r \theta_{n-1}}}{r}  - \frac{\Phi_{\theta_{n-1}}}{r^2} & a_{n-1, 1} & \cdots & a_{n-1, n-1}\\
  \end{array} } \right],
  $$
where
\[
  A = (a_{i,j})_{(n-1)\times (n-1)}
= \frac{\Phi_r}{r} \delta_{ij} 
+ \frac{1}{r^2} \nabla^2_{\mathbb{S}^{n-1}} \Phi.
\]
and $\nabla^2_{\mathbb{S}^{n-1}}  \Phi$ is the spherical Hessian of the function
$\theta \to \Phi(r_0 \cdot \theta)$.
\end{lemma}
\begin{proof}
 To prove the lemma we perform the following computations at $x =x_0 = r_0 \cdot \nu$:
$$
\partial^2_{rr} \Phi(x) 
= \partial_r (\langle \nabla \Phi(x), \nu \rangle )
= \partial_{\nu} \langle \nabla \Phi(x), \nu \rangle 
= \langle D^2 \Phi(x) \nu, \nu \rangle + 
     \langle \nabla \Phi(x), \partial_{\nu} \nu \rangle  =
     \langle D^2 \Phi(x) \nu, \nu \rangle,
     $$
\begin{align*}
\partial^2_{r \theta_i} \Phi(x) 
& = \partial_{r} \bigl( \partial_{\theta_i} \Phi(x) \bigr)
= \partial_{r} \bigl( r \cdot \partial_{e_i} \Phi(x) \bigr)
= \partial_{e_i} \Phi(x) + r  \partial_{r} \bigl(\partial_{e_i} \Phi(x) \bigr)
= \partial_{e_i} \Phi(x) + r  \partial_{\nu} \langle {e_i}, \nabla \Phi(x) \rangle 
\\& = \partial_{e_i} \Phi(x) + r   \langle \partial_{\nu} {e_i}, \nabla \Phi(x) \rangle  + r \langle  {e_i}, D^2 \Phi(x) \cdot \nu \rangle 
= \partial_{e_i} \Phi(x) +  r \langle  {e_i}, D^2 \Phi(x)  \cdot \nu \rangle \\& = \frac{\partial_{\theta_i} \Phi(x)}{r} +  r \langle  {e_i}, D^2 \Phi(x)  \cdot \nu \rangle,
\end{align*}
\begin{align*}
\partial^2_{\theta_i \theta_j} \Phi(x)
& = \partial_{\theta_i} \bigl(\partial_{\theta_j} \Phi(x) \bigr) = r \partial_{e_i} \bigl( r \partial_{e_j} \Phi(x) ) = r^2 \partial_{e_i} \bigl( \langle e_j, \nabla \Phi(x) \rangle )
 = r^2 \Bigl( \langle \partial_{e_i} e_j, \nabla \Phi(x) \rangle + \langle e_j, D^2 \Phi(x) \cdot e_i\rangle \Bigr)
 \\& =  r^2 \Bigl( - \frac{\partial_{\nu} \Phi(x)}{r}  \delta_{ij}+ \langle e_j, D^2 \Phi(x) \cdot e_i\rangle \Bigr)
 = - r \partial_{r} \Phi(x)
+
r^2 \langle e_j, D^2 \Phi(x) \cdot e_i\rangle. 
\end{align*}     
Using these formulas one can easily get the desired expression for $D^2 \Phi$. We remind the reader that $\nabla^2_{\mathcal{S}^{n-1}} \Phi = (\partial^2_{\theta_i \theta_j} \Phi$).
\end{proof}

\begin{corollary}
\label{010724}
    Let $\Phi = r^{\alpha} \phi^{\alpha}$ be a $\alpha$-homogeneous convex function, where $\phi$ is the restriction of $\Phi^{\frac{1}{\alpha}}$
    onto 
    $\mathbb{S}^{n-1}$. Then
    $$
\det D^2 \Phi = (\alpha-1) \alpha^n r^{n(\alpha-2)} \phi^{(\alpha-1)n +1} \det (\phi \delta_{ij} + \nabla^2_{\mathbb{S}^{n-1}} \phi).
    $$
\end{corollary}
\begin{proof}
    Apply previous Lemma. One has
    $$
\Phi_{rr} = \alpha(\alpha-1) r^{\alpha-2} \phi^{\alpha} 
    $$
      $$
\frac{\Phi_{r \theta_i}}{r} - \frac{\Phi_{ \theta_i}}{r^2}
= \alpha(\alpha-1) r^{\alpha-2} \phi^{\alpha-1} \phi_{\theta_i}
    $$
    $$
A = r^{\alpha-2} \Bigl( \alpha \phi^{\alpha} \delta_{ij}
+ \alpha \phi^{\alpha-1} \nabla^2_{\mathbb{S}^{n-1}} \phi 
+ \alpha(\alpha-1) \phi^{\alpha-2} \nabla_{\mathbb{S}^{n-1}} \phi \otimes
\nabla_{\mathbb{S}^{n-1}} \phi \Bigr).
    $$
  Assume that the frame $(e_1, \cdots, e_{n-1})$ is chosen in such a way that the matrix
  $\phi \delta_{ij} + \nabla^2_{\mathbb{S}^{n-1}} \phi$ is diagonal with eigenvalues $\lambda_i$. Then $D^2 \Phi$ takes the form
  $D^2 \Phi =  \alpha r^{\alpha-2} \phi^{\alpha-1} C$, where
    $$
C = \left[ {\begin{array}{cccc}
    (\alpha-1) \phi &(\alpha-1) \phi_{\theta_1}  & \cdots & (\alpha-1)\phi_{\theta_{n-1}}\\
   (\alpha-1) \phi_{\theta_1} & b_{1,1} & \cdots & b_{1, n-1}\\
    \vdots & \vdots & b_{i,j}  & \vdots\\
   (\alpha-1)\phi_{\theta_{n-1}} & b_{n-1, 1} & \cdots & b_{n-1, n-1}\\
  \end{array} } \right],
    $$
    and 
    $$
b_{i,j} = \lambda_i \cdot \delta_{ij} + \frac{\alpha-1}{\phi} \phi_{\theta_i} \phi_{\theta_j}.
    $$
    Elementary computations give $\det C = (\alpha-1) \phi \prod_{i=1}^{n-1} = 
    (\alpha-1) \phi \det (\phi \delta_{ij} + \nabla^2_{\mathbb{S}^{n-1}} \phi)$. This completes the proof.
\end{proof}

    Now let $p>1$, and assume that the potential $V$ is $p$-homogeneous. In particular, $V$ and $V^*$ have the following forms
    $$
V(x) = |x|^{p} v^{p} \left(\frac{x}{|x|}\right), \ \
V^*(x) = |x|^{p^*} \tilde{v}^{p^*}\left(\frac{x}{|x|}\right),
    $$
    where $$p^* = \frac{p}{p-1}.$$

    According to Theorem \ref{p-hom-solutions}, any symmetric maximum point $\Phi$ of  the functional 
    $\mathcal{BS}_{p, V}$ is $p$-homogeneous. Thus the similar representation holds
    $$
\Phi(x) = |x|^{p} \phi^{p} \left(\frac{x}{|x|}\right), \ \
\Phi^*(x) = |x|^{p^*} \tilde{\phi}^{p^*}\left(\frac{x}{|x|}\right).
    $$
Applying the change of variables formula
$$
\frac{1}{\int e^{-\frac{\Phi^*}{p-1}}dy} e^{-\frac{\Phi^*}{p-1}} \det D^2 V^* = \frac{1}{\int e^{-\Phi}dx} e^{-\Phi(\nabla \Phi^*)} \det D^2 \Phi^*,
$$
the relation $\Phi = \frac{1}{p-1} \Phi^*(\nabla \Phi)$
 and Corollary \ref{010724}, we get the following result.
    
    \begin{theorem}
    The following equation holds
    \begin{equation}
    \label{generalMA}
    \tilde{\phi}^{(p^*-1)n +1} \det (\tilde{\phi} \delta_{ij} + \nabla^2_{\mathbb{S}^{n-1}} \tilde{\phi})
= C
\tilde{v}^{(p^*-1)n +1} \det (\tilde{v} \delta_{ij} + \nabla^2_{\mathbb{S}^{n-1}} \tilde{v}),
\end{equation}
where $C = \frac{\int e^{-\frac{\Phi^*}{p-1}}dy}{\int e^{-\Phi}dx}$.
\end{theorem}

The above Theorem establishes that any maximizer of $\bs_{p,V}$ is a solution to a corresponding $L^{q}$-Minkowski problem. This fact gives in a sense more precise information about the relation between the functional and the set versions of the problem (see Proposition \ref{set-func}). 

In particular, uniqueness of solution to (\ref{generalMA}) (for fixed $C, \tilde{v}$ and unknown $\tilde{\phi}$) would imply an affirmative answer to Question \ref{Question-main-V}. 
Unfortunately, it is known that in general equation (\ref{generalMA}) has many  (possibly, infinitely many) solutions for those values of $p$ which are of interest for us (see \cite{CW}, \cite{HLW}, \cite{JLW}, \cite{QRLi}). 

\begin{remark}
Note that a variational problem related to equation (\ref{bmin}) is known in the literature about $L^{q}$-Minkowski problem (see \cite{CW}). Usually it is stated in the following form: maximize the functional
$$
\int_{\mathbb{S}^{n-1}} h^{q} f dx 
$$
with constraint $|K_h|=1$, where $h$ is the support function of $K_h$. Then the solution satisfies equation (\ref{bmin}). This problem is equivalent to some appropriate form of our maximization problem (\ref{KKo}) for sets.
Some other variants of the Blaschke--Santal\'o inequality for sets can be found in \cite{CCL}, \cite{chen2018}.
\end{remark}

\section{Strong Brascamp-Lieb inequalities}\label{section-BL}

In this section we study the following strengthening of $\bs_{p,V}$ on the set of even functions:
\begin{equation}
    \label{lambdaBL}
{\rm Var}_{\mu} f \le \lambda \int \langle (D^2 V)^{-1} \nabla f, \nabla f \rangle d\mu,
\end{equation}
with $\lambda<1$. Here $\mu = \frac{e^{-V(x)}dx}{\int e^{-V(x)}dx}$. As we have seen in Subsection \ref{BLformaximizers}, the maximizers of the generalized Blaschke--Santal\'o functional satisfy (\ref{lambdaBL}) with $\lambda = 1 -\frac{1}{p}$.  We  estimate the best value of $\lambda$ 
in (\ref{lambdaBL}) for  log-concave measures with potential of the form  $V = c \|x\|^p_q$. We prove that in general (\ref{lambdaBL}) fails to hold with  $\lambda = 1 -\frac{1}{p}$. This proves, in particular, that $V$ is not always maximizer for
$\bs_{p,V}$

\subsection{Powers of $l_q$-norms}

In this subsection we study strong Brascamp--Lieb inequality for measure
$$\mu = \frac{e^{-\frac{\|x\|_q^{p}}{p}} dx}{\int e^{-\frac{\|x\|_q^{p}}{p}} dx}
= \frac{e^{-\frac{1}{p} \bigl( \sum_{i=1}^n |x_i|^q \bigr)^{\frac{p}{q}}}dx}{\int e^{-\frac{\|x\|_q^{p}}{p}} dx}.
$$
All the functions below are assumed to be even.

\begin{remark}
    Due to homogeneity invariance the constant $\lambda$ in (\ref{lambdaBL}) remains  the same  for $V = c \|x\|_q^{p} $ with any $c>0$.
\end{remark}

At the first step we do the following change of variables:
$$
x_i = {\rm sign}(y_i) |y_i|^{\frac{2}{q}}.
$$
The reader can easily verify that the image of $\mu$ under the mapping $x \to y(x)$ coincides with
$$
\nu = \frac{1}{C} {e^{-\frac{1}{p} |y|^{\frac{2p}{q}}}} \Bigl( \prod_{i=1}^n |y_i| \Bigr)^{\frac{2}{q}-1} dy.
$$
The measure $\nu$ can be represented in polar coordinates as follows:
$$
\nu = \gamma(dr)  m(d \theta),
$$
where
$$ \gamma = \frac{e^{-\frac{1}{p} r^{\frac{2p}{q}}} r^{\frac{2n}{q}-1} dr}{\int_0^{\infty} e^{-\frac{1}{p} r^{\frac{2p}{q}}} r^{\frac{2n}{q}-1} dr}$$
and
$m$ is a probability measure on $\mathbb{S}^{n-1}$
which has the form
$$
m = \frac{|y_1 \cdots y_n|^{\frac{2}{q}-1} \cdot \sigma}{\int_{\mathbb{S}^{n-1}}|y_1 \cdots y_n|^{\frac{2}{q}-1} d\sigma}, 
$$
where $\sigma$ is the normalized probability surface measure on $\mathbb{S}^{n-1}$.

We are interested 
in the best estimate 
of $\lambda$
in  the strong  Brascamb--Lieb inequality for $\mu$
\begin{equation}
    \label{1905bl}
{\rm Var}_{\mu} f \le \lambda  \int \langle (D^2 V)^{-1} \nabla f, \nabla f \rangle d\mu,
\end{equation}
where $V= \frac{1}{p} \bigl( \sum_{i=1}^n |x_i|^q \bigr)^{\frac{p}{q}}$. 

\begin{remark}
Note that the best value of $\lambda$ we can hope for is $ 1 - \frac{1}{p}$.
Indeed, let  $V$ be $p$-homogeneous and convex. One can easily prove that 
$f = \langle \nabla V(x),x \rangle$ satisfies equality
${\rm Var}_{\mu} f = \lambda \int \langle (D^2 V)^{-1} \nabla f, \nabla f \rangle d\mu$,
with $\lambda = 1 - \frac{1}{p}$.

Indeed, one has $f = pV$, $\nabla f = p \nabla V$ and $\langle (D^2 V)^{-1} \nabla f, \nabla f \rangle = {p^2}\langle (D^2 V)^{-1} \nabla V, \nabla V \rangle  = \frac{p^2}{(p-1)} \langle x, \nabla V(x) \rangle$.
Integrating by parts one gets
$$
\int f d\mu = \int \langle \nabla V(x),x \rangle d\mu = n,
$$
$$
\int f^2 d\mu = p \int V(x)  \langle x, \nabla V(x) \rangle d \mu = np \int V d\mu + p \int  \langle x, \nabla V(x) \rangle d\mu = (n+p) \int \langle x, \nabla V(x) \rangle d\mu = n(n+p).
$$
Thus ${\rm Var}_{\mu} f = np$. On the other hand
$
\int \langle (D^2 V)^{-1} \nabla f, \nabla f \rangle d\mu = \frac{p^2}{p-1} \int \langle x, \nabla V(x) \rangle d\mu = \frac{np^2}{p-1}.
$ This proves the claim.
\end{remark}

Let $$a = ({\rm sign}(x_i) |x_i|^{q-1}), \ \ \ b = ({\rm sign}(x_i) |x_i|^{\frac{q}{2}}).$$
One has $\nabla V = \bigl( \sum_{i=1}^n |x_i|^q \bigr)^{\frac{p}{q}-1} a$ and
\begin{align*}
D^2 V(x) & =  \bigl( \sum_{i=1}^n |x_i|^q \bigr)^{\frac{p}{q}-1}  \Bigl[(q-1)   {\rm diag} (|x_i|^{q-2}) + (p-q) \frac{a \otimes a}{ \sum_{i=1}^n |x_i|^q }\Bigr]
\\&
=
\bigl( \sum_{i=1}^n |x_i|^q \bigr)^{\frac{p}{q}-1} {\rm diag} (|x_i|^{\frac{q}{2}-1}) \Bigl[(q-1)   I + (p-q) \frac{b \otimes b}{ \sum_{i=1}^n |x_i|^q }\Bigr] {\rm diag} (|x_i|^{\frac{q}{2}-1}). 
\end{align*}

 Let us make the change of variables. Apply inequality (\ref{1905bl}) to $f=g({\rm sign}(x_i) |x_i|^{\frac{q}{2}})$. One has
 $$
 {\rm Var}_{\mu} f = {\rm Var}_{\nu} g
 $$
 and
 $$
\langle (D^2 V)^{-1} \nabla f, \nabla f \rangle
=
\frac{q^2}{4|y|^{2(\frac{p}{q}-1)}}
\Big\langle \Bigl[(q-1)   I + (p-q) \frac{y \otimes y}{|y|^2 }\Bigr]^{-1} \nabla g(y), \nabla g(y) \Big\rangle.
 $$
Thus we get that  (\ref{1905bl}) is equivalent to the following inequality:
 $$
{\rm Var}_{\nu} g \le \frac{\lambda q^2}{4(q-1)} \int \frac{1}{|y|^{2(\frac{p}{q}-1)}}
\Big\langle \left( I +  \bigl(\frac{p-q}{q-1}\bigr) \frac{y \otimes y}{|y|^2 }\right)^{-1} \nabla g(y), \nabla g(y) \Big\rangle d\nu.
 $$
 In what follows we denote by 
$
\nabla_{\mathbb{S}^{n-1}} g 
$
the projection of $\nabla g$ onto $\mathbb{S}^{n-1}$:
$$
\nabla_{\mathbb{S}^{n-1}} g(y) = g - \big\langle \nabla g(y) , \frac{y}{r} \big\rangle  \frac{y}{r}.
$$
Note that
$$
\nabla_{\theta} g = \frac{\nabla_{\mathbb{S}^{n-1}} g}{r} .
$$
Thus $\nabla g = g_r \cdot \frac{y}{r} + {\nabla_{\mathbb{S}^{n-1}} g}$, where $|y|=r$. We get that in polar coordinates the last inequality looks like
 \begin{equation}
     \label{2105}
     {\rm Var}_{\nu} g \le \frac{\lambda q^2}{4(q-1)} \int \int \frac{1}{r^{2(\frac{p}{q}-1)}}
\Bigl( \frac{q-1}{p-1} g^2_r 
+ \frac{|\nabla_{\theta} g|^2}{r^2}\Bigr)
 \gamma(dr) m(d\theta).
 \end{equation}

Let $g^r(\theta) = \int_0^{\infty} g(r,\theta) d\gamma$.
One has
$$
{\rm Var}_{\nu} g =
{\rm Var}_{\nu} (g-g^r(\theta))
+
{\rm Var}_{\nu} (g^r(\theta)).
$$
To estimate the first term we apply the following one-dimensional Poincar{\'e}-type inequality and the Fubini theorem.
 \begin{equation}
     \label{sBL1dim+}
     {\rm Var}_{e^{-u}dr} (g) \le  \int_0^{+\infty} \frac{(g')^2}{u^{''}  + \frac{u'}{r}} e^{-u} dr. 
 \end{equation}

Inequality (\ref{sBL1dim+}) is  the $1$-dimensional case of the result obtained by Cordero-Erausquin and Rotem (see Theorem 3 in \cite{CorRot}). Note that one can extend measures and functions symmetrically to get an equivalent inequality on $\mathbb{R}$, and then apply the result from \cite{CorRot} for $n=1$.

In particular, one gets
 $$
{\rm Var}_{\nu} (g-g^r(\theta))
 \le \frac{q^2}{4p} \int\frac{g^2_r}{r^{2 (\frac{p}{q}-1)}} d\nu.$$
To estimate the second term we apply 
 Poncar{\'e} inequality for measure $m$ for even functions with the best constant $C_m$
\begin{align*}
{\rm Var}_{\nu} (g^r(\theta))
= {\rm Var}_{m} (g^r(\theta))
& \le C_m \int |\nabla_{\theta}
 g^r(\theta)|^2 dm
\\& = C_m \int {r^2} \Bigl|\int \nabla_{\mathbb{S}^{n-1}} g d\gamma \Bigr|^2 dm
\le C_m
\int r^{\frac{2p}{q}} d\gamma   \int \frac{|\nabla_{\mathbb{S}^{n-1}} g |^2 }{
 r^{2 \bigl(\frac{p}{q}-1\bigr)} }d\nu.
\end{align*}
Next we compute
$$ \int r^{\frac{2p}{q}} d\gamma = \frac{\int_0^{\infty}e^{-\frac{1}{p} r^{\frac{2p}{q}}} r^{(\frac{2n}{q}-1) + \frac{2p}{q}} dr}{\int_0^{\infty} e^{-\frac{1}{p} r^{\frac{2p}{q}}} r^{(\frac{2n}{q}-1)} dr}
= - \frac{q}{2}\frac{\int_0^{\infty} \bigl(e^{-\frac{1}{p} r^{\frac{2p}{q}}} \bigr)' r^{\frac{2n}{q} } dr}{\int_0^{\infty} e^{-\frac{1}{p} r^{\frac{2p}{q}}} r^{(\frac{2}{nq}-1)} dr} =n.$$
Finally, we get that every even $g$ satisfies
$$
{\rm Var}_{\nu} g
\le  \frac{q^2}{4p} \int\frac{g^2_r}{r^{2 (\frac{p}{q}-1)}} d\nu
+ n C_m    \int \frac{|\nabla_{\mathbb{S}^{n-1}} g |^2 }{
 r^{2 \bigl(\frac{p}{q}-1\bigr)} }d\nu.
$$
Thus comparing this result with  (\ref{2105}) we obtain get the following.
\begin{proposition}
\label{lcm}
    Assume that
$$ \lambda \ge 1 -\frac{1}{p}, \ 
 C_m \le \frac{\lambda q^2}{4n(q-1)}.
$$
    Then inequality (\ref{1905bl}) holds on the set of even functions.
\end{proposition}

 Thus we have reduced our problem to the following question:
 what is the best Poincar{\'e} inequality on the set of even functions for the measure
$m$?

The associated weighted Laplacian for $m$ has the form
$$
Lf = \Delta_{\mathbb{S}^{n-1}} f + \Bigl( \frac{2}{q}-1\Bigr) \langle \omega, \nabla_{\mathbb{S}^{n-1}} f
\rangle,
$$
where $\Delta_{\mathbb{S}^{n-1}}, \nabla_{\mathbb{S}^{n-1}}$
are the spherical Laplacian and gradient and
$$
\omega = \Bigl( \frac{1}{x_1}, \cdots, \frac{1}{x_n}\Bigr).
$$
Thus we have to find the first non-zero eigenvalue on the domain of even functions for operator $L$.

\begin{remark}
    Making appropriate change of variables one can show that for $n=2$, equation $Lf = - \lambda f$ can be reduced to  the so-called Legendre equation and the corresponding eigenvalue functions are known as Legendre functions. In general, they are not elementary,
\end{remark}

For the sake of simplicity the computations below are done on  $\{x_i > 0\}$. 
Recall that Euclidean and spherical Laplacians are related by
$$
\Delta = \partial^2_r + \frac{n-1}{r} \partial_r + \frac{1}{r^2} \Delta_{\mathbb{S}^{n-1}}.
$$
Using this representation, we get immediately

    $$
    L x_i = -\bigl( \frac{2n}{q}-1 \bigr) x_i + \frac{2-q}{q x_i}
    $$
    
    $$
    L x_i^2 = -\frac{4n}{q}\Bigl (x^2_i -\frac{1}{n} \Bigr)
    $$

    $$
L(x_i x_j) = -\frac{4n}{q} x_i x_j 
+ \Bigl( \frac{2}{q}-1\Bigr)
\Bigl( \frac{x_i}{x_j} + \frac{x_j}{x_i}\Bigr), \ i \ne j.
    $$

\begin{lemma}
    For all $i\ne j$ and all $N$ we have
    $$
L\bigl[ (x_i x_j)^N \bigr]=
- 4N  \Bigl[  \frac{n}{q} + N-1 \Bigr] (x_i x_j)^N
+ N \Bigl[ N +  \frac{2}{q} -2 \Bigr] (x_ix_j)^{N-1} (x^2_i + x^2_j). 
    $$
    In particular
    $$L\left[(x_i^2 x_j^2)^{1 - \frac{1}{q}}\right] 
    = -8 \Bigl( 1 - \frac{1}{q} \Bigr) \Bigl( 1 + \frac{n-2}{q} \Bigr) (x_i^2 x_j^2)^{1 - \frac{1}{q}}. $$
\end{lemma}
Since this is a local computation performed only on $\{ x_i > 0,\ x_j >0 \}$, there are no issues in considering the power $(x_i x_j)^N$ also for non-integer $N$. 
\begin{proof}
    One has
    $$
L\bigl[(x_i x_j)^N \bigr]
= N (x_i x_j)^{N-1} L(x_i x_j)
+ N(N-1) (x_i x_j)^{N-2} |\nabla_{\mathbb{S}^{n-1}} (x_i x_j)|^2.
    $$
    One can easily verify:
    $$
\nabla_{\mathbb{S}^{n-1}} (x_i x_j) = x_j e_i + x_i e_j - 2 x_ix_j \cdot x
    $$
    and
    $$
|\nabla_{\mathbb{S}^{n-1}} (x_i x_j)|^2 = x^2_i + x^2_j - 4 x^2_i x^2_j.
    $$
    Thus
    \begin{align*}
        L\bigl[(x_i x_j)^N \bigr]
& = N (x_i x_j)^{N-1} \Bigl[  -\frac{4n}{q} x_i x_j 
+ \Bigl( \frac{2}{q}-1\Bigr)
\Bigl( \frac{x_i}{x_j} + \frac{x_j}{x_i} \Bigr)\Bigr] + N(N-1) (x_i x_j)^{N-2} \Bigl[  x^2_i + x^2_j - 4 x^2_i x^2_j\Bigr]
\\& 
= - N (x_i x_j)^N \Bigl[ \frac{4n}{q} + 4(N-1) \Bigr]
+ N (x_ix_j)^{N-1}  \Bigl[ \Bigl( \frac{2}{q} -1\Bigr) + N-1 \Bigr] (x^2_i + x^2_j).
    \end{align*}
This completes the proof.
\end{proof}
    
  We observe that $L$ preserves even and unconditional functions. Using this observation and the above computations, we obtain the following corollary.

\begin{corollary}
    The following functions are eigenfunctions of  $L$:
    \begin{itemize}
        \item 
        $$
x^2_i -\frac{1}{n}
        $$
        with eigenvalue $-\frac{4n}{q}$.
        \item 
        $$
        |x_i x_j|^{2 \bigl(1 - \frac{1}{q}\bigr)} 
        $$
        with eigenvalue 
        $
        -8 \Bigl( 1 - \frac{1}{q} \Bigr) \Bigl( 1 + \frac{n-2}{q} \Bigr).
        $
    \end{itemize}
\end{corollary}

\begin{theorem}
\label{mpoin}
    The best constant $C_m$ of measure $m$ in the Poincar{\`e} inequality on the set of even functions satisfies
    $$
C_m = \max \Bigl( \frac{q}{4n}, \frac{q^2}{8(q-1)(n+q-2)}\Bigr).
    $$
\end{theorem}
\begin{proof} 

Given even function $f$ we represent  it as follows
\begin{equation}
\label{ffa}
f = \sum_{a \in \{0,1\}^n} f_a,
\end{equation}
where every function $f_a(x_1, \cdots,x_n)$ is even in $x_i$ if $a_i=0$ and odd in $x_i$ if $a_i=1$. For instance, if all $a_i$ are zero, then $f_a$ is unconditional. Note that 
if $a=(a_1, \cdots,a_n)$ contains odd amount of $1$, then $f_a=0$, because $f$ is even.

To obtain this representation we use the operators 
$$
\sigma_i(x) = (x_1, \cdots, -x_i, x_n)
$$ and
$$
T_i^+ f =\frac{f(x) + f(\sigma_i(x))}{2}, 
\ 
T_i^- f =\frac{f(x) - f(\sigma_i(x))}{2}.
$$
Note that $f(x) = T_i^+ f + T_i^- f$, where $T_i^+ f$ is even in $x_i$ and 
$T_i^- f$ is odd in $x_i$. Consequently applying the operators
$$
T_1^{\pm}, T_2^{\pm}, \cdots, T_n^{\pm},
$$
we obtain representation (\ref{ffa}),
where 
$$
f_a = T_1^{b_1} \cdots T_n^{b_n} f.
$$
Here $b_i =1$, if $a_i = 1$ and $b_i=-1$ if $a_i=0$. 

Next we note that 
$$
{\rm Var}_m f = \sum_{a \in \{0,1\}^n} {\rm Var}_m f_a.
$$
This is because $f_a f_b$ is odd at least in one variable for $a \ne b$, hence $\int f_a f_b dm=0$, because the measure $m$ is unconditional.
Similarly, 
$$
\int_{\mathbb{S}^{n-1}} |\nabla_{\mathbb{S}^{n-1}} f|^2 dm= \sum_{a \in \{0,1\}^n} \int_{\mathbb{S}^{n-1}} |\nabla_{\mathbb{S}^{n-1}} f_a|^2 dm.
$$
Indeed, let $a \ne b$. There exists $j$ such that  
$f_a$ (say) is even in $x_j$
and $f_b$ is odd in $x_j$.
Then
$\partial_{x_i} f_a \cdot \partial_{x_i} f_b$ is odd in $x_j$ for all $i$. Indeed, if $i \ne j$, then $\partial_{x_i} f_a $ is even in $x_j$
and $\partial_{x_i} f_b$ is odd in $x_j$.
 If $i =j$, then $\partial_{x_i} f_a $ is odd in $x_j$
and $\partial_{x_i} f_b$ is even in $x_j$
Finally, $\int \partial_{x_i} f_a \cdot \partial_{x_i} f_b d\mu=0$ and
$$
\int \langle \nabla f_a, \nabla f_b 
\rangle d\mu = 
\sum_{i=1}^n \int \partial_{x_i} f_a \cdot \partial_{x_i} f_b \ d\mu  =0.
$$

Thus we have reduced the  statement to the case of 
$f_a$ for arbitrary $f, a$. 

For an unconditional function $f_0$ we have
\begin{equation}
\label{mf0}
{\rm Var}_m f_0 \le \frac{q}{4n} \int_{\mathbb{S}^{n-1}} |\nabla_{\mathbb{S}^{n-1}} f_0|^2 dm.
\end{equation}
Indeed,  Theorem \ref{BLunconditional-p-homo} and Proposition \ref{fromBStoBL} imply inequality (\ref{lambdaBL}) for unconditional functions for $q>1$ with $\lambda=1- \frac{1}{q}$. Then we deduce (\ref{mf0}) applying  (\ref{lambdaBL}) to homogeneous functions (see computations in the next subsection).

Let $a$ has at least one non-zero coordinate. We have shown above that this is an even number.  For simplicity let us assume that $a_1=a_2=1$. Then 
$f_a=0$ on the  sets $\{x_1=0\}$, $\{x_2=0\}$.
Using this observation and  identity
$L g = - \lambda g$, where $g= |x_i x_j|^{2(1 - \frac{1}{q})}$, $\lambda = - 8 \Bigl( 1 - \frac{1}{q} \Bigr) \Bigl( 1 + \frac{n-2}{q} \Bigr)$
one gets
\begin{align*}
\lambda {\rm Var}_m f_a & = \lambda \int_{\mathbb{S}^{n-1}} f^2_a dm =
- \int_{\mathbb{S}^{n-1}} f^2_a \frac{Lg }{g} dm  = 2\int_{\mathbb{S}^{n-1}} \frac{\langle \nabla_{\mathbb{S}^{n-1}} g, \nabla_{\mathbb{S}^{n-1}} f_a \rangle}{g} f_a dm -  \int_{\mathbb{S}^{n-1}} f^2_a \frac{|\nabla_{\mathbb{S}^{n-1}} g|^2}{g^2} dm 
\\& \le  \int_{\mathbb{S}^{n-1}} |\nabla_{\mathbb{S}^{n-1}} f_a|^2 dm.
\end{align*}
This completes the proof.
\end{proof}

\begin{theorem}
\label{BL-homog-powers}
    Assume that
$$ \lambda \ge \max \Bigl( 1 -\frac{1}{p}, 1 - \frac{1}{q}, \frac{1}{2(1 + \frac{q-2}{n})} \Bigr).
$$
    Then inequality (\ref{1905bl}) holds on the set of even functions.

    In particular, if $p \le q$, then one can take
    $$
\lambda = 1 - \frac{1}{q}, \ {\rm if} \ q \ge 2,
    $$
    $$
\lambda = \frac{1}{2(1 + \frac{q-2}{n})}, \ {\rm if} \ q \le 2.
    $$

    Inequality (\ref{1905bl})  is sharp and holds with $\lambda = 1 - \frac{1}{p}$ if
    $$
p \ge q, \ \ p \ge \frac{2(n+q-2)}{n+ 2(q-2)} = 2 - \frac{2(q-2)}{n+2(q-2)}.
    $$
\end{theorem}
\begin{proof}
    The estimate of $\lambda$ follows from Proposition \ref{lcm} and Theorem \ref{mpoin}. The sharpness result follows from the observation that the value $1 - \frac{1}{p}$ in inequality ${\rm Var}_{\mu} f \le \bigl( 1 - \frac{1}{p}\bigr) \int \langle (D^2 V)^{-1} \nabla f, \nabla f \rangle d\mu$, where $V$ is $p$-homogeneous can not be improved, because ${\rm Var}_{\mu} f = \bigl( 1 - \frac{1}{p}\bigr) \int \langle (D^2 V)^{-1} \nabla f, \nabla f \rangle d\mu$ for $f = \langle \nabla V(x), x \rangle$.
\end{proof}

   \subsection{Counterexamples to the strong Brascamp--Lieb inequality}
\label{counter-BL}

Theorem \ref{BL-homog-powers} gives, in particular,  sharp inequalities
${\rm Var}_{\mu} f \le \bigl( 1 - \frac{1}{p}\bigr) \int \langle (D^2 V)^{-1} \nabla f, \nabla f \rangle d\mu$
for 
$$
V = \frac{1}{p} \|x\|^{p}_q
$$
if
$$
p \ge q \ge 2
$$
and
$$
q< 2, \ p \ge  2 - \frac{2(q-2)}{n+2(q-2)}.
$$
Note that in both cases the best constant in the inequality (\ref{lambdaBL}) is strictly bigger than $\frac{1}{2}$ except the Gaussian case $p=q=2$.

Unfortunately, we can not claim that the values of $\lambda$ in other cases are optimal. However, we are able to answer the following natural questions:
\begin{enumerate}
    \item Is it true that inequality (\ref{lambdaBL}) holds with $p=q<2$ and $\lambda = 1 - \frac{1}{q}$?
\item Is it true that inequality (\ref{lambdaBL}) holds with $p=2$, $\lambda=\frac{1}{2}$, and some $q \ne 2$? 
\end{enumerate}

The answers to both questions are negative.
Indeed, in the first case the inequality in question is equivalent (see (\ref{2105})) to
$$    {\rm Var}_{\nu} g \le \frac{ q}{4} \int \int 
\Bigl(  g^2_r 
+ \frac{|\nabla_{\theta} g|^2}{r^2}\Bigr)
 \gamma(dr) m(d\theta)
$$
with 
$$
\gamma = \frac{e^{-\frac{1}{q} r^2}r^{\frac{2n}{q}-1} dr}{\int_0^{\infty} e^{-\frac{1}{q} r^2}r^{\frac{2n}{q}-1} dr}, \ \ 
m = \frac{|y_1 \cdots y_n|^{\frac{2}{q}-1} \cdot \sigma}{\int_{\mathbb{S}^{n-1}}|y_1 \cdots y_n|^{\frac{2}{q}-1} d\sigma}.
$$
Apply this inequality to $g = r^{2} \omega(\theta)$ with $\int \omega dm=0$. One gets
$$
 {\rm Var}_{\nu} g 
 = \int r^{4} d \gamma \cdot \int \omega^2 dm
$$
$$
\frac{q}{4} \int \int 
\Bigl(  g^2_r 
+ \frac{|\nabla_{\theta} g|^2}{r^2}\Bigr)
 \gamma(dr) m(d\theta)
 = \frac{q}{4}\int r^{2} d\gamma \cdot \Bigl( 4 \int \omega^2 dm
 + \int |\nabla_{\theta}\omega|^2 dm \Bigr)
$$
Applying integration by parts, we get
$$
\int r^4 d\gamma = \frac{q}{2} \Bigl( 2 + \frac{2n}{q} \Bigr) \int r^2 d\gamma.
$$
After rearrangement  of the terms 
we get inequality
$
{\rm Var}_m  \omega \le 
\frac{q}{4n} \int |\nabla_{\theta} \omega|^2 dm
$
for an arbitrary symmetric $\omega$.
But this contradicts Theorem \ref{mpoin}, because the best constant in the Poincar{\`e} inequality is $\frac{q^2}{8(q-1)(n+q-2)} > \frac{q}{4n}$.

To prove that inequality 
$${\rm Var}_{\mu} f \le \frac{1}{2} \int \langle (D^2 V)^{-1} \nabla f, \nabla f \rangle d\mu,$$
where $p=2$,
does not hold for all values of $q$ except $2$ we use the same arguments: apply the corresponding equivalent  inequality (\ref{2105}) to function $f$ with appropriate degree of homogeneity (this is $\frac{4}{q}$) and show that it contradicts the sharp estimate obtained in Theorem \ref{mpoin}.
We omit the computations here.

\end{document}